%% file: main.tex
\documentclass[reqno,11pt,noamsfonts]{amsart}
\usepackage[foot]{amsaddr}         

\usepackage[default,scale=0.92]{opensans}
\usepackage[notext]{stix}

\usepackage{algorithm}  
\usepackage{algpseudocode}  

\usepackage{etoolbox}
\patchcmd{\section}{\scshape}{\bfseries}{}{}
\makeatletter
\renewcommand{\@secnumfont}{\bfseries}
\makeatother

\usepackage{amsthm}
\newtheorem{remark}{Remark}
\newtheorem{proposition}{Proposition}

\usepackage{bm}
\usepackage{mathtools}
\usepackage{xfrac} 
\usepackage{accents} 
\usepackage[svgnames]{xcolor}

\usepackage{graphicx}
\usepackage[labelfont=bf]{caption}
\usepackage{subcaption}
\usepackage{enumitem}
\usepackage{array}    
\usepackage{booktabs}  
\usepackage[numbers,sort&compress]{natbib} 

\calclayout

\begin{document}


\include{def}

\definecolor{Mosss}{RGB}{0,192,81}
\definecolor{Reddd}{RGB}{240,40,40}
\definecolor{Mocca}{RGB}{148,82,0}
\newcommand{\RI}  [1]{#1} 
\newcommand{\RII} [1]{#1} 
\newcommand{\RIII}[1]{#1} 
\newcommand{\RIV} [1]{#1} 
\newcommand{\JS}  [1]{#1}

\DeclarePairedDelimiterXPP\BigOSI[2]%
  {\mathcal{O}}{(}{)}{}%
  {\SI{#1}{#2}}


\title
{
Robust semi-implicit multilevel SDC methods for conservation laws
}
\author{Erik Pfister, Jörg Stiller}
\address{%
  TU Dresden, Institute of Fluid Mechanics, 01062 Dresden, Germany}
\email{erik.pfister@tu-dresden.de}

\begin{abstract}
Semi-implicit multilevel spectral deferred correction (SI-MLSDC) methods provide a promising approach for high-order time integration for nonlinear evolution equations including conservation laws.
However, existing methods lack robustness and often do not achieve the expected advantage over single-level SDC.
This work adopts the novel SI time integrators from \cite{TI_Stiller2024} for enhanced stability and extends the single-level SI-SDC method with a multilevel approach to increase computational efficiency.
The favourable properties of the resulting SI-MLSDC method are shown by linear temporal stability analysis for a convection-diffusion problem.
The robustness and efficiency of the fully discrete method involving a high-order discontinuous Galerkin SEM discretization are demonstrated through numerical experiments for the convection–diffusion, Burgers, Euler and Navier–Stokes equations. 
The method is shown to yield substantial reductions in fine-grid iterations compared to single-level SI-SDC across a broad range of test cases. Finally, current limitations of the SI-MLSDC framework are identified and discussed, providing guidance for future improvements.
\end{abstract}


\maketitle


\input{introduction}


\input{problem}


\input{sdc}


\input{mlsdc}


\input{stability}


\input{spatial-discretization}


\input{mlsdc-basics}


\input{numerical-studies}


\input{conclusion}


\section*{Data Availability}

The datasets generated during the current study are available from the corresponding author on reasonable request.

\section*{Funding}

Funding by German Research Foundation (DFG) in frame of the project STI \RII{157/9-1} is gratefully acknowlged. The authors would like to thank ZIH, TU Dresden, for the provided computational resources. Open Access funding is enabled and organized by Projekt DEAL.

\section*{Ethics declarations}

The authors declare that they has no known competing financial interests or personal relationships that could have appeared to influence the work reported in this paper.


\bibliographystyle{abbrvnat}
\bibliography{BM,FE,FV,SE,MG,TI}


\appendix
\input{appendix}

\end{document}

%% file: def.tex



\newcommand   {\mathsc}[1]{\text{{\rmfamily\scshape #1}}}
\renewcommand {\mathit}[1]{\text{{\rmfamily\itshape #1}}}

\newcommand {\CFL}    {\mathit{CFL}}

\newcommand {\EX}     {\mathrm{ex}}       
\newcommand {\IM}     {\mathrm{im}}       

\newcommand{\vs} {\nu_{\mathrm s}}
\newcommand{\cs} {c_{\mathrm s}}
\newcommand{\ds} {\kappa_{\mathrm s}}

\newcommand {\fc}     {\M f_\mathrm{c}}
\newcommand {\hc}     {\M h_\mathrm{c}}
\newcommand {\fd}     {\M f_\mathrm{d}}
\newcommand {\fs}     {\M f_\mathrm{s}}
\newcommand {\fex}    {\M f_\mathrm{ex}}
\newcommand {\fim}    {\M f_\mathrm{im}}

\newcommand {\Ac}     {\M A_\mathrm{c}}
\newcommand {\Ad}     {\M A_{\mathrm{d}}}
\newcommand {\nud}    {\tilde\nu}
\newcommand {\kappad} {\tilde\kappa}

\newcommand {\PO}  {\mathcal M}
\newcommand {\g}   {\mathcal G}
\newcommand {\F}   {\mathcal F}
\newcommand {\Fc}  {\mathcal F_{\!\mathrm{c}}}
\newcommand {\Fd}  {\mathcal F_{\!\mathrm{d}}}
\newcommand {\Fs}  {\mathcal F_{\!\mathrm{s}}}
\newcommand {\Fex} {\mathcal F_{\!\mathrm{ex}}}
\newcommand {\Fim} {\mathcal F_{\!\mathrm{im}}}

\newcommand {\lbf} {{\ell\!}}

\newcommand {\iu}  {\mathrm{i}}
\newcommand {\zi}  {z_\mathrm{i}}
\newcommand {\zr}  {z_\mathrm{r}}

\newcommand {\V}  [1] {\vec{#1}}                       
\newcommand {\T}  [1] {\bm{#1}}                        
\newcommand {\TS} [1] {\bm{#1}}                        
\newcommand {\TC} [1] {\mathbf{#1}}                    
\newcommand {\M}  [1] {\bm{#1}}                        
\newcommand {\MC} [1] {\mathbf{#1}}                    

\renewcommand {\d}            {\partial}
\newcommand   {\D}            {\mathrm{d}\mspace{1.0mu}}

\newcommand {\NM} [1] {\underaccent{\bar}{#1}}         
\newcommand {\SM} [1] {\underaccent{\tilde}{#1}}       

\providecommand {\abs} [1] {\left\vert#1\right\vert}
\providecommand {\norm}[1] {\left\Vert#1\right\Vert}

\renewcommand {\d}            {\partial}
\newcommand   {\supp}         {\operatorname{supp}}
\newcommand   {\diag}         {\operatorname{diag}}

\newcommand   {\transpose}[1] {#1^{\mathrm{t}}}
\newcommand   {\jmp}      [1] {\left\lBrack #1\right\rBrack}   
\newcommand   {\avg}      [1] {\left\lBrace #1\right\rBrace}   


\newcommand   {\ek}   {\frac{1}{2}v^2}  
\newcommand   {\ei}   {e}               
\newcommand   {\et}   {E}               
\newcommand   {\htot} {H}               

\newcommand{\N} [2][] {N_{\textsc{#2} #1}}

%% file: introduction.tex
\section{Introduction}

The use of higher-order numerical methods for conservation laws in fluid mechanics and related areas has increased significantly due to their superior convergence properties. 
Consequently, there is a corresponding need for time integration schemes that match the convergence rates and efficiency of spatial discretization techniques.
This is essential to prevent an imbalance when solving partial differential equations (PDEs), where spatial and temporal errors should scale appropriately.
The objective is to develop time integration schemes that enable the stable solution of PDEs with arbitrarily high order in both space and time. 
During the past decades, this challenge has led to the development of various high-order methods, including explicit, implicit, and semi-implicit schemes incorporating multiple steps or stages \cite{TI_Gottlieb2016a}, \cite{TI_Hairer1993a}, \cite{TI_Hairer1996a}. 
Explicit time integration methods become increasingly inefficient as stiffness increases, which is often the case in the presence of unresolved waves, diffusion, or local mesh refinement. This issue is particularly pronounced in high-order spatial finite element methods.
To overcome the restrictive time step limitations imposed by stiffness, one strategy is to utilize fully implicit methods, such as diagonally implicit Runge-Kutta methods \cite{TI_Kennedy2016a}, \cite{TI_Pan2021a}, \cite{TI_Pazner2017a}, Rosenbrock methods \cite{TI_Bassi2015a}, Taylor methods \cite{TI_Baeza2020a} or discontinuous Galerkin methods in time \cite{TI_Tavelli2018a}.
Unfortunately, these fully implicit methods frequently result in ill-conditioned algebraic systems that are costly to solve, making it challenging for them to compete with explicit methods in time-resolved simulations.
Semi-implicit methods handle non-stiff convection explicitly and stiff diffusion implicitly, resulting in manageable algebraic systems using standard linear algebra.
Common techniques include implicit-explicit (IMEX) multistep and Runge-Kutta methods, balancing efficiency and stability.
In practice, however, IMEX multistep methods typically do not exceed an order of $3$ \cite{TI_Fehn2017a}, \cite{TI_Karniadakis1991a}, \cite{TI_Klein2015a}. 
A more favourable alternative is the use of IMEX Runge-Kutta methods, which offer improved stability properties and lower error constants \cite{TI_Ascher1997a}, \cite{TI_Cavaglieri2015a}, \cite{TI_Kennedy2003a}, \cite{TI_Pareschi2005a}.\\
Recently, significant advances have been made in the development of higher-order semi-implicit time integration methods through the use of multi-derivative techniques  \cite{TI_Frolkovic2023a}, \cite{TI_Schuetz2021a}, \cite{TI_Schuetz2022} and \cite{TI_Zeifang2023}.
Another promising approach is the use of spectral deferred correction (SDC) methods, which enable the construction of time integration schemes with arbitrarily high order \cite{TI_Dutt2000a}.
Semi-implicit spectral deferred correction (SI-SDC) methods were introduced \cite{TI_Minion2003b} and have since undergone further evolvement \cite{TI_Christlieb2011a}, \cite{TI_Christlieb2015a}, \cite{TI_Layton2005a}, \cite{TI_Stiller2020a}.
The semi-implicit SDC methods proposed in \cite{TI_Stiller2024} even proof to be $L$-stable up to order 11. 
However, a key challenge with SDC methods is the relatively large number of iterations required to achieve convergence, highlighting the need for techniques that enhance efficiency.\\
So-called Ladder methods, as referenced in \cite{TI_Layton2009a}, \cite{TI_Minion2003b}, initiate with fewer substeps in early SDC iterations to align with the solution's accuracy, gradually advancing from a coarse, low-order solution to a refined, high-order one.
These methods employ one or more SDC sweeps at a coarse level, then interpolate the outcome in both time and space to serve as the provisional solution for subsequent correction sweeps.
The authors conclude however, that while ladder methods reduce workload, this benefit is counterbalanced by a decrease in accuracy, resulting in no computational efficiency gains.
In \cite{TI_Layton2009a}, the results indicate that this problem can be overcome, if the ladder methods allow both spatial and temporal coarsening.
\RI{
Similar to the ladder approach, \citet{TI_Micalizzi2024, TI_Micalizzi2025} proposed a systematic analysis  for an iterative arbitrary high order deferred correction (DC) framework, demonstrating enhanced efficiency compared to classical DC through the newly introduced interpolation process for the purely explicit case.
While \cite{TI_Micalizzi2024} covers only the temporal setting, the extension to both space and time is discussed in \cite{TI_Micalizzi2025}.
}
An extension of \RI{the Ladder} concept involves employing a general multilevel strategy, following \citet{MG_Brandt2011} in the context of SDC methods, as \RII{proposed} by \citet{MG_Speck2014}.
A comparable approach is implemented in the Parallel Full Approximation Scheme in Space and Time (PFASST) \cite{TI_Emmett2012a}, \cite{TI_Minion2010a}, \cite{MG_Bolten2016a} or \cite{MG_Schoebel2019a}.
~\citet{MG_Kremling2021} proved the convergence of multilevel SDC (MLSDC) and examined the role of the coarsening strategy.
Furthermore, \citet{MG_Hamon2019a} introduced an IMEX MLSDC method for the shallow water equations with a spatial discretization based on the global Spherical Harmonics transform.
Finally, \citet{MG_Emmett2019} developed an adaptive MLSDC approach applied to the Navier-Stokes equations, representing another step toward \RII{saving} runtime. However, it \RII{employs} an explicit treatment of convection and diffusion and is restricted to an order of four in space and time.\\
This paper develops a novel SI-MLSDC method that extends the \RII{semi-implicit} single-level SDC method proposed in \cite{TI_Stiller2024} and thereby achieves vastly improved stability properties compared to existing SI-MLSDC methods \RII{and at the same time reduces the number of iterations needed by the single-level SDC method}.
Furthermore, earlier studies focus on spatial discretizations based on finite difference or finite volume methods, while high-order discontinuous Galerkin spectral element methods (DG-SEM) have received comparatively little attention in this context.
Additionally, no prior work known to the authors systematically investigates the influence of starting strategies to the MLSDC algorithm and the influence of the transfer operators.
Moreover, none of the aforementioned works utilize embedded interpolation as a \RI{fine-to-coarse} projection method, which is feasible in our approach due to the availability of spatial and temporal basis functions.
In addition, the distinction between incremental and non-incremental MLSDC formulations has also not been explicitly explored. 
Most importantly, none of the existing studies incorporate the novel semi-implicit (SI) time integrators from \cite{TI_Stiller2024} that offer enhanced stability properties.
The approach presented here introduces a highly stable and robust SI-MLSDC scheme combined with DG-SEM for conservation laws, enabling $p$-coarsening in time as well as simultaneous $p$- and $h$-coarsening in space, thereby addressing these gaps - while also critically examining current limitations of MLSDC. 
The paper is structured as followed: Section~\ref{ch:problem} provides a concise overview of the one-dimensional conservation laws under consideration.
Section \ref{ch:sdc} and \ref{ch:mlsdc} present the proposed SI-SDC and SI-MLSDC methods,
followed by an investigation of the stability properties of the methods with a linear stability analysis of the \RIII{convection-diffusion} problem in Section \ref{ch:stability}.
Spatial discretization, including stabilization techniques, is described in Section \ref{ch:spatial_discretization}.
Section \ref{ch:mlsdc_basics} provides preliminary studies for the application of the SI-MLSDC method to partial differential equations and Section \ref{ch:numerical_experiments} extends the studies by onvection-diffusion, Burgers, and compressible flow
problems.
A summary of findings is provided in Section~\ref{ch:conclusion}.

%% file: problem.tex
\section{Problem}
\label{ch:problem}

A system of conservation laws for the variables ${\M u(x,t) \in \mathbb R^d}$, governed by the equation
\begin{align}
\d_t \M u + \d_x \fc(\M u) = \d_x \fd(\M u) + \fs(\M u, x,t)
\end{align}
is considered where ${x \in \Omega \subset \mathbb R}$ and ${t \in \mathbb R^{+}}$. The vector ${\M u(x,t) \in \mathbb R^d}$ represents the conserved quantities, $\fc$ denotes the convective fluxes, $\fd$ represents the diffusive fluxes, and $\fs$ accounts for any source terms. By introducing the combined right-hand side
\begin{align}
\M f = -\d_x \fc + \d_x \fd + \fs \,,
\end{align}
the system can be rewritten in the following evolutionary form:
\begin{align}
\d_t \M u = \M f(\M u, x, t) \,.
\end{align}
It is assumed that the convective fluxes possess the Jacobian ${\Ac = \fc'(\M u)}$ and
the diffusive fluxes can be written in the form ${\fd = \Ad \d_x \M u}$, where
 ${\Ad(\M u,\d_x \M u)}$
is the positive-semidefinite diffusion matrix.
These definitions yield the quasilinear form
\begin{equation}
  \label{eq:conservation-system:quasilinear}
  \d_t \M u = -\Ac \d_x \M u + \d_x \Ad \d_x \M u + \fs
  \,.
\end{equation}

\subsection{Scalar conservation laws}
For these problems, the variable $ \M u $ \RII{consists} only of one scalar component $u$.
Two particular equations are of interest. The first is the convection-diffusion equation
\begin{align}
\label{eq:conv-diff}
\d_t u = -\d_x (vu) + \d_x(\nu \d_x u) \,,
\end{align}
where $v$ is a constant velocity.
The second case is the Burgers equation, which introduces nonlinear convection and includes an additional source term. It is given by
\begin{align}
\label{eq:burgers}
\d_t u = - \d_x \bigl(\tfrac{1}{2}u^2\bigr) + \d_x(\nu \d_x u) + f_{\mathrm s}(x,t) \,,
\end{align}
where $f_{\mathrm s}$ represents a source term which \RII{can} be tailored to match a specific solution.
In both equations ${\nu \ge 0}$ is a constant diffusivity.

\subsection{Euler and  Navier-Stokes equations}
\label{sec:conservation-laws:cns}

For the Navier-Stokes equations, the conservative variables, convective fluxes and diffusive fluxes given by
\begin{equation}
  \M u = \begin{bmatrix}
           \rho \\ \rho v \\ \rho \et
         \end{bmatrix}
  \,, \quad
  \fc = \begin{bmatrix}
          \rho v \\ \rho v^2 + p\\ \rho v \et + pv
        \end{bmatrix}
  \,, \quad
  \fd = \begin{bmatrix}
           0 \\ 
           \frac{4}{3} \eta \d_x v \\ 
           \frac{4}{3} \eta v \d_x v + \lambda \d_x T
        \end{bmatrix}
  \,, 
\end{equation}
where
  $\rho$ 
is the density,
  $v$
the velocity,
  $\et$ 
the total specific energy,
  $\eta$
the dynamic viscosity and
  $\lambda$
the heat conductivity.
Assuming a perfect gas with the gas constant ${R = c_p - c_v}$ and the ratio of the specific heats ${\gamma = c_p/c_v}$ it is possible to calculate the temperature
  ${T = (\et - v^2/2)/c_v}$,
the pressure
  ${p = \rho R T}$,
the total specific enthalpy ${\htot = \et + p/\rho}$
and
the speed of sound ${a = ((\gamma-1) R T)^{1/2}}$.
The convective Jacobian and the diffusion matrix are given by
\begin{equation}
  \label{eq:convective-jacobian:cns}
  \Ac =
      \begin{bmatrix}
        0                                 & 1                     & 0        \\[1mm]
        \frac{1}{2}(\gamma-3)v^2          & (3-\gamma)v           & \gamma-1 \\[1mm]
        \frac{1}{2}(\gamma-1)v^3 - v\htot & \htot - (\gamma-1)v^2 & \gamma v
      \end{bmatrix}
\end{equation}
and
\begin{equation}
  \label{eq:diffusion-matrix:cns}
  \Ad =
     \begin{bmatrix}
          0                                               
       &  0                      
       &  0 
       \\[1mm]
         -\frac{4}{3}\nu v                              
       &  \frac{4}{3}\nu         
       &  0
       \\[1mm]
         -\bigl(\frac{4}{3}\nu - \gamma\kappa\bigr)v^2 - \gamma\kappa \et
       &  \bigl(\frac{4}{3}\nu - \gamma\kappa\bigr)v 
       &  \gamma\kappa
     \end{bmatrix}
     \,,
  \end{equation}
where
${\nu = \eta/\rho}$
is the kinematic viscosity
and
${\kappa = \lambda /(\rho c_p)}$
the thermal diffusivity.
The inviscid case with ${\eta = \lambda = 0}$ leads to the Euler equations for which
${\Ad = \M 0}$.

%% file: sdc.tex


\newcommand{\Restrict}    {\NM{\mathcal R}}
\newcommand{\Project}     {\NM{\mathcal P}}
\newcommand{\Interpolate} {\NM{\mathcal I}}


\section{Single-level collocation and SDC methods}
\label{ch:sdc}
\subsection{Notation}
Before discussing the numerical time integration methods in detail, we introduce the required notation. The fundamental necessity for the methods considered in this section is to subdivide the time interval $[t, t+\Delta t]$ using a set of collocation points $\{t_m\}_{m=1}^{M}$, such that
\begin{align}
t = t_0 \leq t_1 < \dots < t_m < \dots < t_M = t + \Delta t.
\end{align}
Depending on the choice of collocation points, different integration properties can be achieved. The two point distributions used in this study are:
\begin{itemize}
    \item Radau Right nodes: exclude the left endpoint \RI{($t_1 > t_0$)}.
    \item Lobatto nodes: include both endpoints \RI{($t_1 = t_0$)}.
\end{itemize}
The initial value of the solution is denoted by $\M u_0 = \M u(t_0)$, and the numerical solution at the collocation nodes is represented by $\M u_m \approx \M u(t_m)$ for $m = 1, \dots, M$. Those nodal solutions, excluding the initial value $\M u_0$, are grouped into a solution vector:
\begin{equation}
\NM{\M u} = \bigl[ \{\M u_m\}_{m=1}^{M} \bigr].
\end{equation}
\begin{remark}
    $\mathcal{L}(\NM {\M u},t) = \sum_{m=1}^M \M u_m \ell_m(t)$ is the Lagrange interpolant constructed through the collocation nodes $t_m$. The polynomials $\ell_m(t)$ serve as the basis functions.
\end{remark}
Since the methods under consideration are iterative, the solution at iteration $k$ is denoted by $\M u^k_m \approx \M u(t_m)$. This notation will be used consistently throughout the following sections.


\subsection{Incremental collocation method}

The derivation of SDC starts from the Picard integral form of a generic Initial Value Problem. An approximation to this integral can be achieved by utilizing a quadrature formula:
\begin{equation}
  \label{eq:cm:incremental}
  \M u_{m} = \M u_{m-1} + \Delta t \sum_{i=1}^{M} w^{\mathsc{nn}}_{m,i} \M f(\M u_i,t_i),
  \quad
  m = 1,\dots, M,
\end{equation}
where the incremental integration weights are defined as
\RII{
\begin{align}
\label{eq:nn_weights}
w^{\mathsc{nn}}_{m,i} = \frac{1}{\Delta t}\int_{t_{m-1}}^{t_m} \ell_{\! i}(t) \D t.
\end{align}
}
Note that $\M u_1 = \M u_0$ for Lobatto points. 
The above formulation is also referred to as the \textit{node-to-node} or incremental collocation method (CM), indicated by the $\mathsc{nn}$ superscript.
Rewriting the collocation method in an operator form, we obtain:
\begin{equation}
  \label{eq:cm:incremental:op}
  \M F^{\mathsc{nn}}_{m}(\NM{\M u}) 
    \equiv \M u_{m} - \M u_{m-1} - \Delta t \sum_{i=1}^{M} w^{\mathsc{nn}}_{m,i} \M f(\M u_i,t_i)
    = \M 0\,, \quad
  m = 1,\dots M.
\end{equation}
Consequently, the collocation method can be interpreted as a system of equations:
\begin{equation}
  \label{eq:cm:incremental:matrix}
  \NM{\M F}^{\mathsc{nn}}(\NM{\M u})  = \NM{\M 0}.
\end{equation}
This \RII{system} is solved iteratively using the incremental SDC method, which is explained below.


\subsection{Incremental SDC method}

The SDC method serves as an iterative solver for the collocation problem through low-order substepping methods. Refer to \cite{TI_Ong2020a}, \cite{TI_Minion2003b}, \cite{TI_Dutt2000a}, \cite{TI_Stiller2020a} for a detailed derivation.
For the SDC scheme in this work, predictors and correctors are used, which can be one-staged or two-staged. The predictor is a low-order time integrator, stepping through the subinterval points $m$ to provide the $0$-th approximation to the solution $\M u^{0}_{m}$.
The specific time integration methods used in this work are the IMEX Euler and the one- and two-stage semi-implicit time integration methods SI(1) and SI(2) developed in  \cite{TI_Stiller2024}.
The one-stage IMEX Euler predictor reads
\begin{align}
\label{eq:sdc:incremental_predictor_eu}
\M u^{0}_{m}
= \M u^0_{m-1} 
+ \M H^{\mathsc{eu}}(\M u^{0}_{m-1}, \M u^{0}_{m}, \Delta t_m)
 \,,\quad
  m = 1,\dots M.
\end{align}
Note that the superscript $0$ is not an exponent, but an iteration count. Each iteration $k$ is commonly referred to as a sweep.
The \RII{second term on the right side} represents the low-order approximation of the integral $\int_{t_{m-1}}^{t_m} ( -\d_x \fc + \d_x \fd) \D t$  over the $m$-th subinterval.
For the IMEX Euler, it reads
\begin{equation}
\label{eq:H_IMEX}
\M H^{\mathsc{eu}}(\M u_a, \M u_b, \Delta t_m) =  \Delta t_m \left[ -\d_x \fc(\M u_a) + \d_x \Ad(\M u_a) \d_x \M u_b \right].
\end{equation}
The one-stage predictor SI(1) reads
\begin{align}
\label{eq:sdc:incremental_predictor_si1}
\M u^{0}_{m}
= \M u^0_{m-1} 
+ \M H^{\mathsc{si}}(\M u^{0}_{m-1}, \M u^{0}_{m-1}, \M u^{0}_{m}, \Delta t_m)
 \,,\quad
  m = 1,\dots M,
\end{align}
and the two-stage predictor SI(2)
\begin{align}
\begin{split}
\label{eq:sdc:incremental_predictor_si2}
\M u^{(1)}
& = \M u^0_{m-1} 
+ \M H^{\mathsc{si}}(\M u^0_{m-1}, \M u^0_{m-1}, \M u^{(1)}, \Delta t_m)
 \,,\quad
  m = 1,\dots M, \\
\M u^{0}_{m}
& = \M u^0_{m-1} 
+ \M H^{\mathsc{si}}(\M u^{(1)}, \M u^0_{m-1}, \M u^{0}_m, \Delta t_m)
 \,,\quad
  m = 1,\dots M \,, 
\end{split}
\end{align}
where
\begin{equation}
    \label{eq:H_SI(1)}
\M H^{\mathsc{si}}(\M u_a, \M u_b, \M u_c, \Delta t_m) =  \Delta t_m \left[ -\d_x \fc(\M u_{a}) + \d_x \left(\frac{\Delta t_m}{2} \Ac(\M u_{b})^2 + \Ad(\M u_{b})\right) \d_x \M u_c \right].
    \end{equation}
These correctors are motivated by the implicit Lax-Wendroff method and have excellent stability properties \cite{TI_Stiller2024}.
The corrector equation of the SDC method to calculate the $k$-th approximation is given by
\begin{equation}
  \label{eq:sdc:incremental_corrector}
  \M u^{k}_{m} = \M u^{k}_{m-1} 
               + \Delta t \sum_{i=1}^{M} w^{\mathsc{nn}}_{m,i} \M f(\M u^{k-1}_i,t_i)
               + \Delta \M H^k_m
  \,,\quad
  m = 1,\dots M.
\end{equation}
The deltas of the low-order terms  $\M H^k_m$ for the different time integrations schemes are defined as follows. Firstly, the IMEX Euler:
\begin{equation}
\label{eq:DeltaH_EU}
\Delta \M H^{k,\mathsc{eu}}_m = \M H^{\mathsc{eu}}(\M u^{k}_{m-1}, \M u^{k}_{m}, \Delta t_m) - \M H^{\mathsc{eu}}(\M u^{k-1}_{m-1}, \M u^{k-1}_{m}, \Delta t_m).
\end{equation}
The SI(1) and SI(2) deltas are
\begin{equation}
\label{eq:DeltaH_SI1}
       \Delta \M H^{k,\mathsc{si}(1)}_m = \M H^{\mathsc{si}}(\M u^{k}_{m-1}, \M u^{k}_{m-1}, \M u^{k}_{m}, \Delta t_m) - \M H^{\mathsc{si}}(\M u^{k-1}_{m-1}, \M u^{k-1}_{m-1}, \M u^{k-1}_{m}, \Delta t_m)
\end{equation}
and, respectively, 
\begin{equation}
\label{eq:DeltaH_SI2}
       \Delta \M H^{k,\mathsc{si}(2)}_m = \M H^{\mathsc{si}}(\M u^{(1)}, \M u^k_{m-1}, \M u^{k}_m, \Delta t_m) - \M H^{\mathsc{si}}(\M u^{k-1}_{m}, \M u^{k-1}_{m-1}, \M u^{k-1}_{m}, \Delta t_m)
\end{equation}
with
\begin{equation}
\begin{split}
        \M u^{(1)} & = \M u^{k}_{m-1} 
               + \Delta t \sum_{i=1}^{M} w^{\mathsc{nn}}_{m,i} \M f(\M u^{k-1}_i,t_i) \\
               & + \M H^{\mathsc{si}}(\M u^k_{m-1}, \M u^k_{m-1}, \M u^{(1)}, \Delta t_m) - \M H^{\mathsc{si}}(\M u^{k-1}_{m-1}, \M u^{k-1}_{m-1}, \M u^{k-1}_{m}, \Delta t_m)
  \,,\quad
  m = 1,\dots M.
\end{split}
\end{equation}
\begin{remark}
    Implementation note: In practice, the \RI{$\M f$} terms from sweep $k-1$ are stored to ensure the efficient \RI{computation of the $\M H$ terms} in subsequent sweeps.
\end{remark}
For readability, the $\mathsc{eu}$, $\mathsc{si(1)}$, or $\mathsc{si(2)}$ superscripts are omitted in the following discussion, as the SDC formulation remains identical across these time integration schemes.
Upon convergence \,${\Delta \M H^k_m \rightarrow \M 0}$ \,
and the incremental CM \eqref{eq:cm:incremental} is recovered.
Selecting Radau Right points for the collocation leads to the Radau IIA method resulting in a convergence of order ${2M-1}$ (see \cite[Ch. 6]{TI_Deuflhard2002a}) while selecting Lobatto points leads to the Lobatto IIIA method of order ${2M-2}$ (see \cite[Ch. 8]{TI_Deuflhard2002a}).
The residual of the incremental SDC method for the $k$-th approximation reads
\begin{equation}
  \label{eq:sdc:incremental:residual}
  \M r^{k}_{m} = \M u^{k}_{m} - \M u^{k}_{m-1} 
               - \Delta t \sum_{i=1}^{M} w^{\mathsc{nn}}_{m,i} \M f(\M u^{k-1}_i,t_i)
  \,,\quad
  m = 1,\dots M.
\end{equation}


\subsection{Non-incremental collocation method}
Until now, the method discussed has swept through the subintervals from node to node. There exists an equivalent formulation, called \textit{zero-to-node} or non-incremental formulation, for which the collocation method reads
\begin{equation}
  \label{eq:cm:non-incremental}
  \M u_{m} = \M u_{0} + \Delta t \sum_{i=1}^{M} w^{\mathsc{0n}}_{m,i} \M f(\M u_i,t_i)
  \,,\quad
  m = 1,\dots M \,,
\end{equation}
where zero-to-node weights are given by
\begin{equation}
  \label{eq:0n_weights}
  w^{\mathsc{0n}}_{m,i} = \sum_{j=1}^{m} w^{\mathsc{nn}}_{j,i}.
\end{equation}
This form \RII{is identical to} the standard form of the CM and can be achieved by recursive substitution of ${\M u_{m-1}}$ in \eqref{eq:cm:incremental}.
As a system, the non-incremental CM reads
\begin{equation}
  \label{eq:cm:non-incremental:matrix}
  \NM{\M F}^{\mathsc{0n}}(\NM{\M u})  = \NM{\M 0}.
\end{equation}
\begin{remark}
In collocation methods, the coefficients are denoted $a_{ij}$ as, for example, in \cite[Ch. 2]{TI_Hairer1993a}. Here, the notation $w^{\mathsc{0n}}_{m,i}$ similar to the incremental method is chosen for better comparison. 
\end{remark}

\subsection{Non-incremental SDC method}

Similarly, one obtains the non-incremental SDC method
\begin{equation}
  \label{eq:sdc:non-incremental}
  \M u^{k}_{m} = \M u_{0} 
               + \Delta t \sum_{i=1}^{M} w^{\mathsc{0n}}_{m,i} \M f(\M u^{k-1}_i,t_i)
               + \sum_{i=1}^{m}
                \Delta \M H^k_i
  \,,\quad
  m = 1,\dots M.
\end{equation}
Note that the correction terms in \eqref{eq:sdc:non-incremental} are composed of the same deltas as in the incremental form, but are summed up over the current and all preceding subintervals. The convergence of the zero-to-node SDC method recovers the non-incremental collocation method \eqref{eq:cm:non-incremental}.
\begin{remark}
SDC can also be interpreted as a preconditioned Picard iteration for the collocation problem, see e.g. \cite{TI_Akramov2024}, \cite{TI_Speck2017}, \cite{TI_Caklovic2024}.
\end{remark}

%% file: mlsdc.tex



\section{Multilevel collocation and SDC methods}
\label{ch:mlsdc}

All multilevel techniques derived and applied in this section originate from Brandt \cite{MG_Brandt2011} and can be comprehended there in great detail.
Firstly, a sequence of grid levels $l$, from a coarse discretization of the computational domain ($l = 1$) to a fine discretization ($l = L$) is considered.
The variable $\NM{\M u}$ denotes from now on the solution on all levels. 
The basic idea is to utilize a coarser grid solution to accelerate convergence on a finer grid or, conversely, to leverage the finer grid to improve the solution on the coarser grid.
These objectives are achieved with the Full Approximation Scheme (FAS). With FAS there is ideally a fully-fledged solution to the fine grid problem at all levels.


\subsection{Incremental multilevel collocation}

The nodal solution vector on level \( l = 1, \dots, L \) is denoted by \(\NM{\M u}_l = \bigl[ \{\M u_{l,m}\}_{m=1}^{M_l} \bigr]\). 
The FAS-ML (Full Approximation Scheme Multi-Level) framework is designed to approximate solutions across multiple levels.  
The representation on level $l$ is given by
\begin{equation}
  \label{eq:ml-cm:incremental:fas-mg}
  \NM{\M F}^{\mathsc{nn}}_l(\NM{\tilde{\M{u}}}_l) = \NM{\M g}^{\mathsc{nn}}_l(\NM{\tilde{\M{u}}}_{l+1}) 
  \,,\quad
  l = 1, \dots, L \,,
\end{equation}
\RI{where $\NM{\tilde{\M{u}}}_{l+1}$ denotes an approximate solution on level $l+1$ after one or more SDC sweeps.}
The primary distinction between the single-level and multilevel case is the introduction of the FAS right-hand side (RHS), \(\NM{\M g}^{\mathsc{nn}}_l\). This addition ensures that all levels collectively approximate the solution on the finest grid. \RII{This allows a dual point of view (compare Section~8.2 in \cite{MG_Brandt2011}), where (i) the multilevel FAS coupling accelerates fine-level convergence through coarse-level corrections and (ii) the fine level discretization is used to correct the coarse grid solution.
}
The FAS RHS is defined by
\begin{equation}
  \label{eq:ml-cm:incremental:fas-mg:rhs}
  \NM{\M g}^{\mathsc{nn}}_l(\NM{\tilde{\M{u}}}_{l+1}) =
  \begin{cases}
     \NM{\M F}^{\mathsc{nn}}_l(\Project_l\NM{\tilde{\M{u}}}_{l+1})
     + \Restrict_l 
         \bigl( \NM{\M g}^{\mathsc{nn}}_{l+1} 
              - \NM{\M F}^{\mathsc{nn}}_{l+1}(\NM{\tilde{\M{u}}}_{l+1})
         \bigr)
     & l < L \\
     0
     & l = L
  \end{cases}\,\,.
\end{equation}
The term $\NM{\M g}^{\mathsc{nn}}_l(\NM{\tilde{\M{u}}}_{l+1})$ is computed prior to the correction sweep at the start of every multigrid cycle and contains information from the finer grids, detailed further in algorithm \ref{alg:mlsdc:v-cycle}.
Transfer operators are introduced to facilitate interactions between different levels.
To transfer residuals to the next coarser level, we use the residual restriction \(\NM{\M r}_{l} = \Restrict_l \NM{\M r}_{l+1}\). 
Similarly, the solution projection to the next coarser level is represented by \(\NM{\M v}_{l} = \Project_l \NM{\M u}_{l+1}\), 
and the solution interpolation to the next finer level is given by \(\NM{\M u}_{l} = \Interpolate_l \NM{\M u}_{l-1}\).
\begin{remark}
    The transfer operators are algebraic matrices and are applied to the coefficient vectors.
\end{remark}
Convergence \RI{of this} FAS-ML \RI{formulation} is characterized by several key behaviors. In addition to the convergence of the solution \(\NM{\tilde{\M{u}}}_l \rightarrow \NM{\M{u}}_l\) at each level,
the residual for levels $l < L$ converge $\NM{\M{r}}^{\mathsc{nn}}_l = \NM{\M g}^{\mathsc{nn}}_l - \NM{\M F}^{\mathsc{nn}}_{l+1}(\NM{\tilde{\M{u}}}_l) \rightarrow \NM{\M{0}}$.
Concurrently, the right-hand side (RHS) $\NM{\M g}^{\mathsc{nn}}_l$ for $l < L$ trends towards $\NM{\M F}^{\mathsc{nn}}_l(\Project_l\NM{\M{u}}_{l+1})$, emphasizing the alignment of the RHS computation with the projected solutions from higher levels.
These properties result in the formulation of the incremental multilevel collocation method, defined by the equations:
\begin{subequations}
  \label{eq:ml-cm:incremental}
  \begin{alignat}{2}
    &  \NM{\M F}^{\mathsc{nn}}_l(\NM{\M{u}}_l) 
    &&  = \NM{\M F}^{\mathsc{nn}}_l(\Project_l\NM{\M{u}}_{l+1})
        \,,\quad
        l = 1, \dots L-1,
  \\
    &  \NM{\M F}^{\mathsc{nn}}_L(\NM{\M{u}}_L) 
    && = \NM{\M 0}.
  \end{alignat}
\end{subequations}


\subsection{MLSDC transfer operators}
\label{sec:transferoperators}
The solution is interpolated from level $l-1$ to level $l$ using the \RII{coarse-to-fine} interpolation operator \RII{$\Interpolate_{l}$}
\begin{align}
\label{eq:interpolation}
\NM{\M u}_{l} = \Interpolate_{l}\NM{\M u}_{l-1}.
\end{align}
This operator relies on embedded interpolation techniques, where the operator $\Interpolate_{l}$ transfers the coarse solution accurately to the fine grid with no other than round-off errors.
Similarly, the solution can be projected from level \RII{$l+1$ to $l$} using the projection operator \RII{$\Project_{l}$}
\begin{align}
\label{eq:projection_int}
\NM{\M v}_{l} = \Project_{l}\NM{\M u}_{l+1} \,,
\end{align}
for which there is a choice between $L^2$-projection and embedded interpolation.
While $L^2$-projection is conservative, it is less accurate than embedded interpolation as it relies on the lower level basis functions.
In this paper, both methods for projection are implemented and investigated later on.
The residual restriction operator $\Restrict_{l}$ from level $l$ to $l-1$ is defined as the transpose of the \RII{coarse-to-fine} interpolation operator:
\begin{align}
\label{eq:restriction}
\Restrict_{l} \NM{\M r}_{l+1} = \transpose{\Interpolate_{1+1}} \NM{\M r}_{l+1}.
\end{align}


\subsection{Incremental multilevel SDC method}

Similarly to SDC, the multilevel SDC method is formulated as an iterative solver to the collocation method.
The corrector equations on level $l$ are given by
\begin{equation}
  \label{eq:ml-sdc:incremental}
  \M u^{k}_{l,m} 
    = \M u^{k}_{l,m-1} 
    + \Delta t \sum_{i=1}^{M_l} w^{\mathsc{nn}}_{l,mi} \M f(\M u^{k-1}_{l,i},t_{l,i})
    + \M g^{\mathsc{nn}}_{l,m}
    + \Delta \M H^k_{l,m} 
  \,,\quad
  m = 1,\dots M_l.
\end{equation}
The definitions of $\Delta \M H^k_{l,m}$ for the three different time integrators are given in equations \eqref{eq:DeltaH_EU} to \eqref{eq:DeltaH_SI2}.
\begin{remark}
    The $k$-th iteration $\NM {\M u}^k_l$ is always an approximation to the discrete solution $\NM {\M u}_l$ and therefore one specific $\tilde{\NM {\M u}}_l$. The tilde $\tilde{\NM {\M u}}_l$ indicator is left out in favor of readability.
\end{remark}
To express the MLSDC correction scheme in a compact form, its matrix form is introduced:
\begin{equation}
  \label{eq:ml-sdc:incremental:matrix}
  \NM{\M F}^{\mathsc{nn}}_l(\NM{\M{u}}_l^k) 
  = \NM{\M g}^{\mathsc{nn}}_l
  + \NM{\M C}^{\mathsc{nn}}_l \Delta \NM{\M H}^k_{l},
\end{equation}
where $\NM{\M C}^{\mathsc{nn}}_l$ is identity matrix and can be omited
\begin{equation}
  \NM{\M C}_l^{\RII{\mathsc{nn}}}
  = \begin{bmatrix}
      1      & 0      & \cdots  & 0      \\
      0      & \ddots &         & \vdots \\
      \vdots &        & \ddots  & 0      \\
      0      & \cdots &  0      & 1
    \end{bmatrix} \,.
\end{equation}
Upon convergence \, ${\Delta \NM{\M H}^k_{l} \rightarrow \NM{\M 0}}$ \,
and the incremental multilevel CM is recovered. The procedure for an incremental MLSDC correction sweep can be seen in algorithm \ref{alg:mlsdc:corrector}, the index $l$ dropped for readability.
\begin{algorithm}[ht]
\caption{Incremental MLSDC correction sweeps.}
\label{alg:mlsdc:corrector}
\begin{algorithmic}[1]
\Procedure{MLSDC\_Corrector}{}${(\NM{\M u},\NM{\M g}, N)}$
  \State
    ${\NM{\M u}^{0} \gets \NM{\M u}}$ 
  \For{$k = 1,N$}
    \State
      ${\M u_0^k \gets \M u_0}$
    \For{$m = 1,M$}
      \State
         ${\M u_m^k \gets 
           \textsc{Solve}
            \big(  \M u^{k}_{m}
                =  \M u^{k}_{m-1} 
                + \Delta t \sum_{i=1}^{M} w^{\mathsc{nn}}_{m,i} \M f(\M u^{k-1}_{i},t_{i})
                + \Delta \M H^k_{m} + g^{\mathsc{nn}}_{m}
            \big)
         }$
    \EndFor
  \EndFor
  \State
    ${\NM{\M u} \gets \NM{\M u}^{N}}$ 
\EndProcedure
\end{algorithmic}
\end{algorithm}

The incremental MLSDC residual reads

\begin{equation}
    \label{MLSDC:incremental:residual}
    \NM{\M r}^{\mathsc{nn}}_{l} = \NM{\M g}^{\mathsc{nn}}_{l} - \NM{\M F}^{\mathsc{nn}}_l(\NM{\M{u}}_l) \,.
\end{equation}


\subsection{Non-incremental multilevel collocation}

An analogous procedure to the non-incremental CM method yields the non-incremental multilevel collocation method, which is defined by the equations:
\begin{subequations}
  \label{eq:ml-cm:non-incremental}
  \begin{alignat}{2}
    &  \NM{\M F}^{\mathsc{0n}}_l(\NM{\tilde{\M{u}}}_l) 
    &&  = \NM{\M F}^{\mathsc{0n}}_l(\Project_l\NM{\M{u}}_{l+1})
        \,,\quad
        l = 1, \dots L-1 \,,
  \\
    &  \NM{\M F}^{\mathsc{0n}}_L(\NM{\tilde{\M{u}}}_L) 
    && = \NM{\M 0}.
  \end{alignat}
\end{subequations}


\subsection{Non-incremental multilevel SDC method}

The matrix form of the non-incremental MLSDC method reads
\begin{equation}
  \label{eq:ml-sdc:non-incremental:matrix}
  \NM{\M F}^{\mathsc{0n}}_l(\NM{\M{u}}_l^k) 
  = \NM{\M g}^{\mathsc{0n}}_l
  + \NM{\M C}^{\mathsc{0n}}_l \Delta \NM{\M H}^k_{l},
\end{equation}
with
\begin{equation}
  \NM{\M C}^{\mathsc{0n}}_l
  = \begin{bmatrix}
      1      & 0      & \cdots & 0      \\
      1      & 1      & \cdots & 0      \\
      \vdots & \vdots & \ddots & \vdots \\
      1      & 1      & \cdots & 1
    \end{bmatrix} \,.
\end{equation}
Again the non-incremental multilevel collocation method is recovered upon convergence.


\subsection{Non-equivalence of incremental and non-incremental multilevel-methods}

Equivalence is obvious for multilevel collocation, but not for multilevel SDC.
For single-level CM and SDC, the incremental and non-incremental formulations are equivalent via recursion, but for MLSDC equivalence requires
\begin{equation}
\label{eq:condition_equivalence}
  \NM{\M g}^{\mathsc{0n}}_{l} 
  = \NM{\M C}^{\mathsc{0n}}_l \,\NM{\M g}^{\mathsc{nn}}_{l}.
\end{equation}
For a two-level case the condition \eqref{eq:condition_equivalence} can be demonstrated. An illustrative example is provided in Appendix~\ref{sec:app:mlsdc:equivalence}.
Substituting the definition of the RHS and exploiting
${\NM{\M F}^{\mathsc{0n}}_l = \NM{\M C}^{\mathsc{0n}}_l \NM{\M F}^{\mathsc{nn}}_l}$
yields
\begin{equation}
  \Restrict_l 
    \bigl( \NM{\M g}^{\mathsc{0n}}_{l+1} 
         - \NM{\M F}^{\mathsc{0n}}_{l+1}(\NM{\M{u}}^k_{l+1})
    \bigr)
  = \NM{\M C}^{\mathsc{0n}}_l \,
       \Restrict_l 
         \bigl( \NM{\M g}^{\mathsc{nn}}_{l+1} 
              - \NM{\M F}^{\mathsc{nn}}_{l+1}(\NM{\M{u}}^k_{l+1})
         \bigr).
\end{equation}
For ${l = L-1}$ we have ${\NM{\M g}^{\mathsc{nn}}_{l+1} = \NM{\M g}^{\mathsc{0n}}_{l+1} = \NM{\M0}}$, which implies
\begin{equation}
  \Restrict_{L-1} 
          \NM{\M F}^{\mathsc{0n}}_{L}(\NM{\M{u}}^k_{L})
  = \NM{\M C}^{\mathsc{0n}}_{L-1} \,
       \Restrict_{L-1}
              \NM{\M F}^{\mathsc{nn}}_{L}(\NM{\M{u}}^k_{L}),
\end{equation}
or equivalently
\begin{equation}
  \Restrict_{L-1} 
      \NM{\M C}^{\mathsc{0n}}_{L} \NM{\M F}^{\mathsc{nn}}_{L}(\NM{\M{u}}^k_{L})
  = \NM{\M C}^{\mathsc{0n}}_{L-1} \,
       \Restrict_{L-1}
              \NM{\M F}^{\mathsc{nn}}_{L}(\NM{\M{u}}^k_{L}).
\end{equation}
As the construction of the transfer operators is independent of the mesh level it follows that
\begin{equation}
    \Restrict_{l} \NM{\M C}^{\mathsc{0n}}_{l+1}
  = \NM{\M C}^{\mathsc{0n}}_{l} \Restrict_{l}
\end{equation}
is a necessary condition for the equivalence of incremental and non-incremental MLSDC.
The above condition is not universally satisfied and can be refuted by a numerical example.
\RI{
\begin{remark}
The matrix $\NM{\M C}_l^{\mathsc{0n}}$ has dimensions $M_l \times M_l$, i.e.\ $\NM{\M C}_l^{\mathsc{0n}} \in \mathbb{R}^{M_l \times M_l}$. Matrix–vector products are to be understood as index contractions, e.g., $(\NM{\M C}_l^{\mathsc{0n}} \NM{\M g}^{\mathsc{nn}}_l)_m = \sum_{n=1}^{M_l} \M C_{l,mn}^{\mathsc{0n}} \M g^{\mathsc{nn}}_{l,n} $ 
\end{remark}
}
Consider $L = 2$ with $M_2 = 3$, $M_1 = 2$ and Radau Right points with $p$ coarsening in time. The operators for both cases result in:
\[
\NM{\M C}^{\mathsc{0n}}_{l} \Restrict_{l} = 
\begin{bmatrix}
1.22499 & 0.33594 & -0.07321 & 0.00000 \\
0.87674 & 1.19970 & 0.69907 & 0.00000 \\
1.00000 & 1.00000 & 1.00000 & 1.00000 \\
\end{bmatrix}
\]
and
\[
\Restrict_{l} \NM{\M C}^{\mathsc{0n}}_{l+1} = 
\begin{bmatrix}
1.48773 & 0.26274 & -0.07321 & 0.00000 \\
1.28778 & 1.63603 & 0.77227 & 0.00000 \\
1.22449 & 1.10123 & 1.30093 & 1.00000 \\
\end{bmatrix} \,.
\]
\begin{proposition}
    The incremental and non-incremental MLSDC formulations are not equivalent.
\end{proposition}
The outcome is similar for Lobatto and equidistant collocation points, which shows that incremental and non-incremental MLSDC are not equivalent, although they converge to the same solution.


\subsection{Multilevel algorithm}

\begin{figure}[H]
    \centering
    \begin{subfigure}[b]{0.4\textwidth}
        \includegraphics[width=\textwidth]{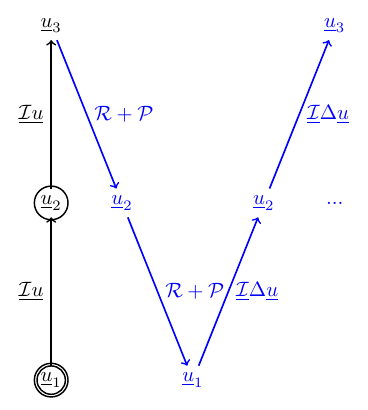}
        \caption{Cascadic starting strategy}
        \label{fig:cascade_sketch}
    \end{subfigure}
    \hfill
    \begin{subfigure}[b]{0.57\textwidth}
        \includegraphics[width=\textwidth]{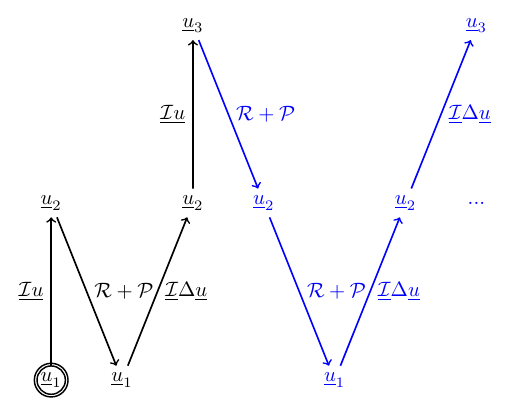}
        \caption{FMG$_1$ starting strategy}
        \label{fig_fmg1_sketch}
    \end{subfigure}
    \vspace{1em}
    \caption{Sketches for the starting strategies Cascade and FMG$_1$ \RII{, i.e. one cycle per added level,} \RII{each with} one attached \RII{full} V-cycle \RII{in blue}}
    \label{fig:start_strategies_sketches}
\end{figure}

In this subsection, we provide a detailed explanation of the MLSDC framework, which consists of a starting strategy combined with a V-Cycle. 
The algorithms for the V-cycle (algorithm \ref{alg:mlsdc:v-cycle}) and the various strategies for initiating an MLSDC algorithm  (algorithms \ref{alg:mlsdc:cascade} and \ref{alg:mlsdc:fmg}) are presented below. 
All of them are formally identical for incremental and non-incremental MLSDC, so the $\mathsc{nn}$ and $\mathsc{0n}$ indices are dropped.
\begin{algorithm}[ht]
\caption{MLSDC V-cycle.}
\label{alg:mlsdc:v-cycle}
\begin{algorithmic}[1]
\Procedure{Incremental MLSDC\_V\_Cycle}{}{$(\NM{\M u}, N_{c}, L)$}
  \State 
    ${\NM{\M g}_{L} \; \gets \; \NM{\M 0}}$
  \For{$l = L,2,-1$}
    \State 
      ${\textsc{SDC\_Corrector}(\NM{\M u}_l, \NM{\M g}_l, N_s = 1)}$ 
      \textcolor{Grey}{\Comment{{Pre-sweeping}}}
    \State
      \makebox[2.2em][l]{${\NM{\M r}_{l}}$}
        ${\gets \; \NM{\M g}_l -\NM{\M F}_{l}(\NM{\M u}_{l}) }$
        \textcolor{Grey}{\Comment{{Residual evaluation}}}
    \State
      \makebox[2.2em][l]{${\NM{\M v}_{l-1}}$}
        ${ \gets \; \Project_{l-1}\NM{\M u}_{l}}$
        \textcolor{Grey}{\Comment{{Solution restriction}}}
      \State
        \makebox[2.2em][l]{${\NM{\M g}_{l-1}}$}
          ${ \gets \; \NM{\M F}_{l-1}(\NM{\M v}_{l-1})
                    + \Restrict_{l-1}\NM{\M r}_{l} 
          }$
        \textcolor{Grey}{\Comment{{RHS composition}}}
  \EndFor
  \State 
    ${\textsc{SDC\_Corrector}(\NM{\M u}_1, \NM{\M g}_1, N_{c})}$ 
    \textcolor{Grey}{\Comment{{Coarse solution}}}
  \For{$l = 2,L$}
      \State
        \makebox[1.5em][l]{${\NM{\M u}_{l}}$}
          ${ \gets \; \NM{\M u}_{l}
                    + \Interpolate_{l}(\NM{\M u}_{l-1} - \NM{\M v}_{l-1})}$
          \textcolor{Grey}{\Comment{{Correction prolongation}}}
    \If{$l < L$}
      \State 
      ${\textsc{SDC\_Corrector}(\NM{\M u}_l, \NM{\M g}_l, N_s = 1)}$ 
      \textcolor{Grey}{\Comment{{Post-sweeping}}}
    \EndIf
  \EndFor
\EndProcedure
\end{algorithmic}
\end{algorithm}
The V-cycle is started at the fine level $l = L$ with a SDC sweep. Pre-sweeping is used to \RII{reduce} high-frequency error components on the finer grid, while low-frequency errors are removed on the coarse grid.
\RI{For a detailed investigation of the role of the SDC corrector as a smoother we refer to \cite{MG_Kremling2021}.}
\RII{The} solution and residuals are \RII{then} restricted to the coarser grids with the operators described in \ref{sec:transferoperators}.
The FAS right-hand side $\NM{\M g}_l$ is subsequently calculated with equation \eqref{eq:ml-cm:incremental:fas-mg:rhs}.
At the bottom level $l=1$, the problem is solved using $N_{c}$ coarse sweeps.
Finally, the correction is interpolated to the finer grids (line 11 in algorithm \ref{alg:mlsdc:v-cycle}). This correction is the key for boosting convergence.
A post-sweeping step is performed for levels $1 < l < L$, after that the cycle starts again with the current solution. Only on the last cycle a post-sweeping SDC sweep is optionally performed on the finest grid.
The number of post- and pre-sweeping steps $N_s$ could differ from $1$, conducted numerical studies do not suggest any benefit from that, so all tests were simply done with $N_s = 1$ sweep for pre- and post-sweeping, \RII{if not stated otherwise}.\\
Initial solutions on every grid are required with FAS. Therefore, certain starting strategies to produce $\NM{\M u}^{0}_L$ are needed.
One strategy is to spread the initial conditions as a constant through all nodes \cite{MG_Speck2014}.
The second strategy investigated in the present work is to start with the predictor solutions on all levels.
The third strategy, referred to as the Cascade (algorithm \ref{alg:mlsdc:cascade}) begins with a predictor on the coarse grid, followed by a single correction sweep \cite{MG_Huismann2016a}. The solution is then interpolated onto the next finer grid, where another correction is performed. 
On the finest grid, the V-cycle starts with the pre-sweeping step for the interpolated solution from the second finest grid.
\begin{algorithm}[ht]
\caption{MLSDC Cascade.}
\label{alg:mlsdc:cascade}
\begin{algorithmic}[1]
\Procedure{MLSDC\_Cascade}{}{$(\NM{\M u}, L)$}
  \State
   ${ \textsc{SDC\_Predictor}}(\NM{\M u}_{1})$ 
      \textcolor{Grey}{\Comment{{Predictor Solution Level 1}}}
  \For{$l = 1,L-1$}
    \State 
    ${ \textsc{SDC\_Corrector}(\NM{\M u}_{l}, \NM{\M g}_l = \NM{ \M {0}}, N_s = 1 })$ 
      \textcolor{Grey}{\Comment{{Corrector Sweep}}}
    
      \State
        \makebox[1.5em][l]{${\NM{\M u}_{l+1}}$}
          ${ \gets \; \Interpolate_{l+1}(\NM{\M u}_{l})}$
          \textcolor{Grey}{\Comment{{Solution interpolation}}}
  \EndFor
\EndProcedure
\end{algorithmic}
\end{algorithm}
The fourth and final approach is utilizing MLSDC as a Full Multigrid (FMG) method, as outlined in algorithm \ref{alg:mlsdc:fmg} (compare FMG in \cite[Fig. 1.2]{MG_Brandt2011}).
The algorithm \ref{alg:mlsdc:fmg} shows the FMG starting variant with the number of cycles \RII{$C_{\mathsc{fmg}}$} per added level. Rather than improving the initial solution with a corrector sweep, as done in the Cascade variant, we substitute the corrector with one or more V-Cycles on a growing subset of levels.
To ensure that each level converges to its solution before moving to the next, there is the possibility to perform multiple V-cycles per added level, contrary to just one per level like in \cite{MG_Brandt2011}).
Note that for 2 levels, the starting strategies from algorithm \ref{alg:mlsdc:cascade} and \ref{alg:mlsdc:fmg} are identical.
\begin{algorithm}[H]
\caption{MLSDC FMG.}
\label{alg:mlsdc:fmg}
\begin{algorithmic}[1]
\Procedure{MLSDC\_FMG}{}{$(\NM{\M u}, N_{c}, L, C_{\mathsc{fmg}})$}
\State
   ${ \textsc{SDC\_Predictor}}(\NM{\M u}_{1})$ 
      \textcolor{Grey}{\Comment{{Predictor Solution Level 1}}}

\State 
   ${ \textsc{SDC\_Corrector}(\NM{\M u}_{1}, \NM{\M g}_1 = \NM{ \M {0}}, N_s = 1)}$ 
      \textcolor{Grey}{\Comment{{Corrector Sweep Level 1}}}

\State
        \makebox[1.5em][l]{${\NM{\M u}_{2}}$}
          ${ \gets \; \Interpolate_{2}(\NM{\M u}_{1})}$
          \textcolor{Grey}{\Comment{{Solution interpolation}}}
          
\For{$l = 2,L-1$}
          
     \For{$c = 1,C_{\mathsc{fmg}}$} 
         \State 
            ${ \textsc{MLSDC\_V-cycle}(\NM{\M u}, N_{c}, l)}$ 
              \textcolor{Grey}{\Comment{{$C_{\mathsc{fmg}}$ MLSDC V-cycles}}}
     \EndFor

 \State
        \makebox[1.5em][l]{${\NM{\M u}_{l+1}}$}
          ${ \gets \; \Interpolate_{l+1}(\NM{\M u}_{l})}$
          \textcolor{Grey}{\Comment{{Solution interpolation}}}     
     
\EndFor 
\EndProcedure
\end{algorithmic}
\end{algorithm}
For illustration, the cascadic and FMG strategies are sketched in Figure~\ref{fig:start_strategies_sketches}. The coarse level, corresponding to $l = 1$, is located at the bottom of the sketch. 
Both strategies begin with the coarse solution at the bottom left, the notation $\mathcal{R} + \mathcal{P}$ represents fine-to-coarse operations (restriction and projection), while $\mathcal{I} \Delta u$ denotes coarse-to-fine prolongation of the correction. 
The \RI{circular markers in Figure~\ref{fig:start_strategies_sketches} denote} specific operations: a single circle \RI{indicates} a corrector step, \RI{while a double circle} corresponds to a predictor–corrector sequence. This notation \RI{is used only outside the V-cycles. Within the V-cycles, the SDC sweeps and are not explicitly marked.}
The first \RI{full} V-Cycle after the start is highlighted in blue.


%% file: stability.tex
\section{Temporal stability and accuracy}
\label{ch:stability}

In this section, the stability and accuracy properties of the MLSDC method are investigated. Starting point for the stability analysis is the convection-diffusion equation \eqref{eq:conv-diff} with a spatially periodic complex solution ${u(x,t) = \hat u(t) \mathrm{e}^{\iu k x}}$ and a constant velocity $v \in \mathbb R$, which yields
\begin{equation}
    \partial_t \hat u \mathrm{e}^{\iu k x} = -\d_x (v \hat u \mathrm{e}^{\iu k x}) + \d_x (\nu \d_x \hat u \mathrm{e}^{\iu k x}).
\end{equation}
With the definitions from Section \ref{ch:problem}, we note that 
\begin{equation}
    \begin{split}
        f_c & = v \hat u(t) \mathrm{e}^{\iu k x} = \hat f_c(t)\mathrm{e}^{\iu k x}, \\
        f_d & = \nu \d_x \hat u(t) \mathrm{e}^{\iu k x} = \hat f_d(t) \mathrm{e}^{\iu k x}.
    \end{split}
\end{equation}
It follows that $A_c = \hat A_c = v$ and $A_d = \hat A_d = \nu$. This can now simply be plugged into the predictor and corrector equations shown in Section \ref{ch:sdc} and Section \ref{ch:mlsdc}. If the spatial derivatives are evaluated, $\mathrm{e}^{\iu k x}$ can be cancelled out and a Dahlquist-type ODE \cite{TI_Dahlquist1963} in time is obtained,
\begin{equation}
  \label{eq:dahlquist}
  \D_t \hat u = \lambda \hat u.
\end{equation}
Here ${\lambda \in \mathbb C}$ with the real part ${\lambda_{\mathrm r} = -\nu k^2}$ and the imaginary part ${\lambda_{\mathrm i} = -v k}$.
\begin{remark}
Although the MLSDC methods were developed for PDEs with a real solution $u(x,t)$, they can be easily adopted for the present problem. After removing the common factor $\mathrm{e}^{\iu k x}$, the complex amplitudes $\hat H^k_m$, $\Delta \hat H^k_m$ and $\hat f(\hat u_i,t_i)$ take the role of  $H^k_m$, $\Delta H^k_m$ and \RI{$f(u_i,t_i)$} respectively.
\end{remark}

\subsection{Stability}

\begin{figure}[H]
    \centering
    \begin{subfigure}[b]{0.45\textwidth}
        \includegraphics[width=\textwidth]{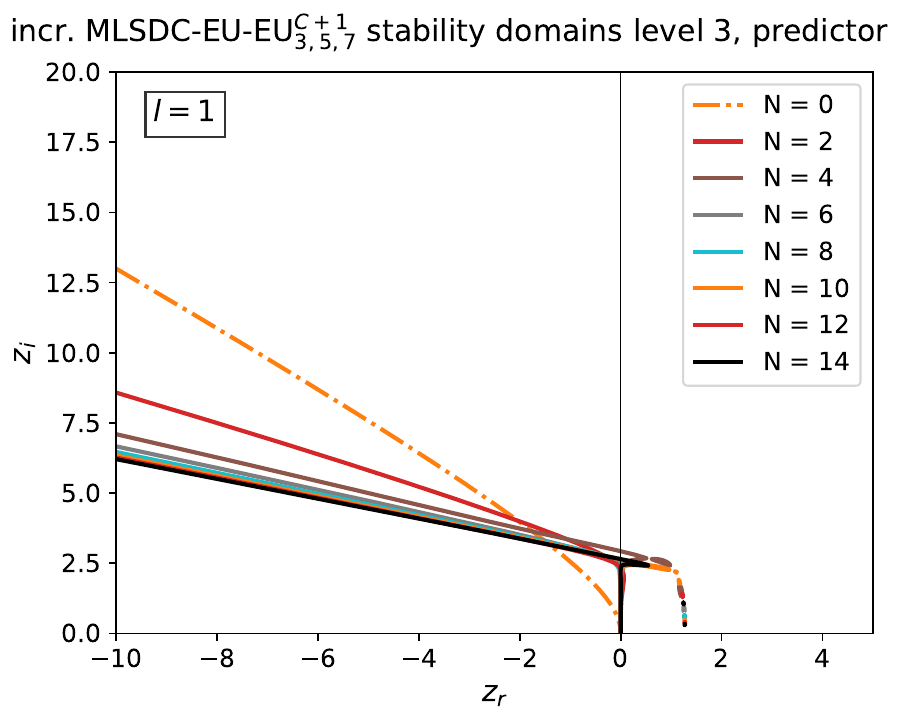}
        \caption{Bottom level $l = 1$ of MLSDC with IMEX Euler on all levels, predictor variant}
        \label{xxa}
    \end{subfigure}
    \hfill
    \begin{subfigure}[b]{0.45\textwidth}
        \includegraphics[width=\textwidth]{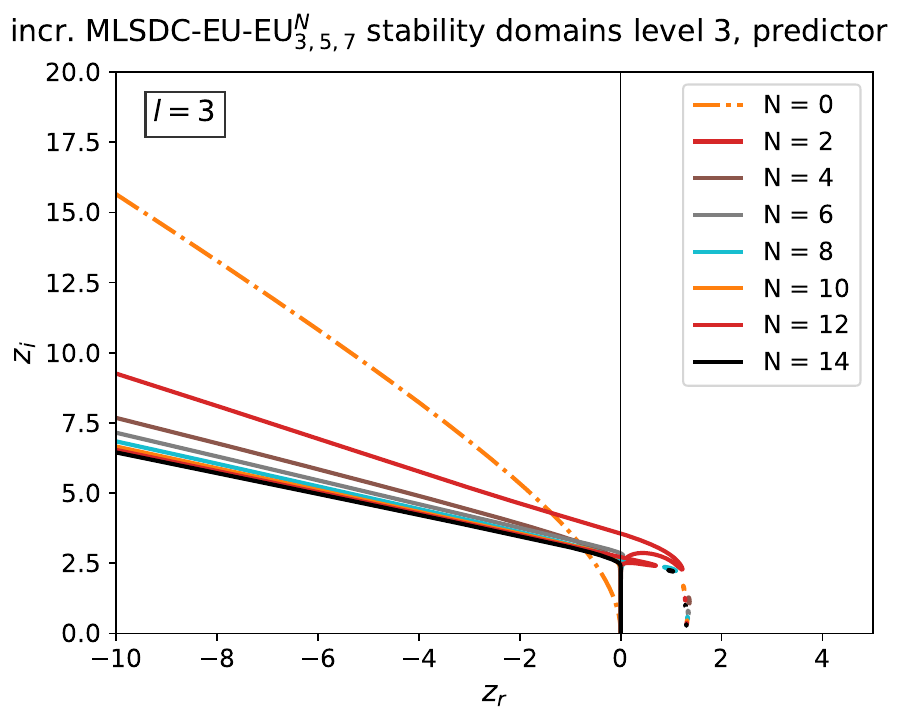}
        \caption{Top level $l = 3$ of MLSDC with IMEX Euler on all levels, predictor variant}
        \label{xxb}
    \end{subfigure}
    \vspace{1em} 
    \begin{subfigure}[b]{0.45\textwidth}
        \includegraphics[width=\textwidth]{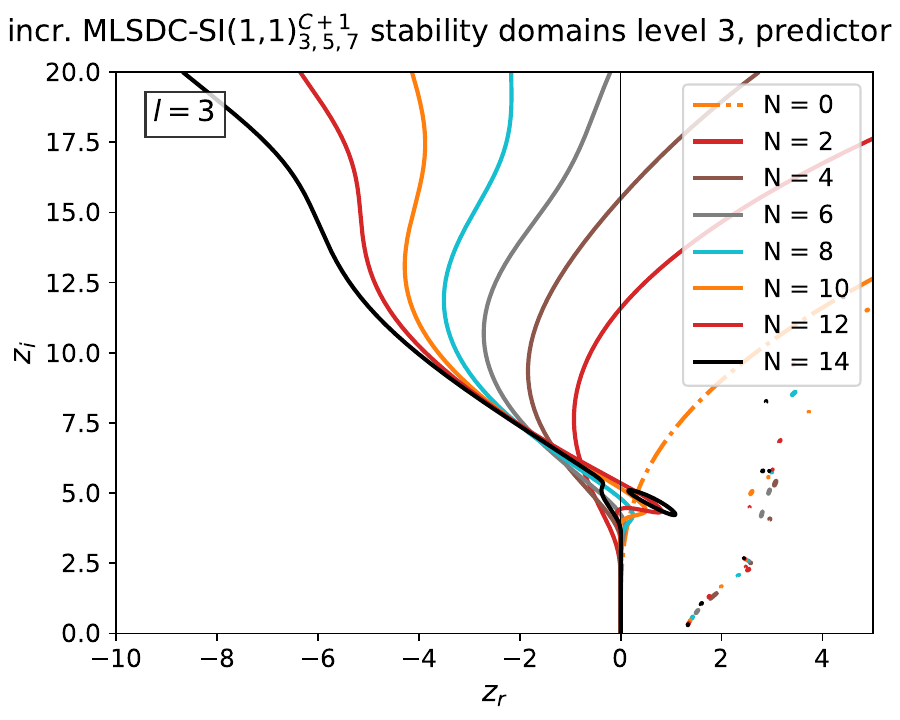}
        \caption{Top level $l = 3$ of MLSDC with SI(1) on all levels, predictor variant}
        \label{xxc}
    \end{subfigure}
    \hfill
    \begin{subfigure}[b]{0.45\textwidth}
        \includegraphics[width=\textwidth]{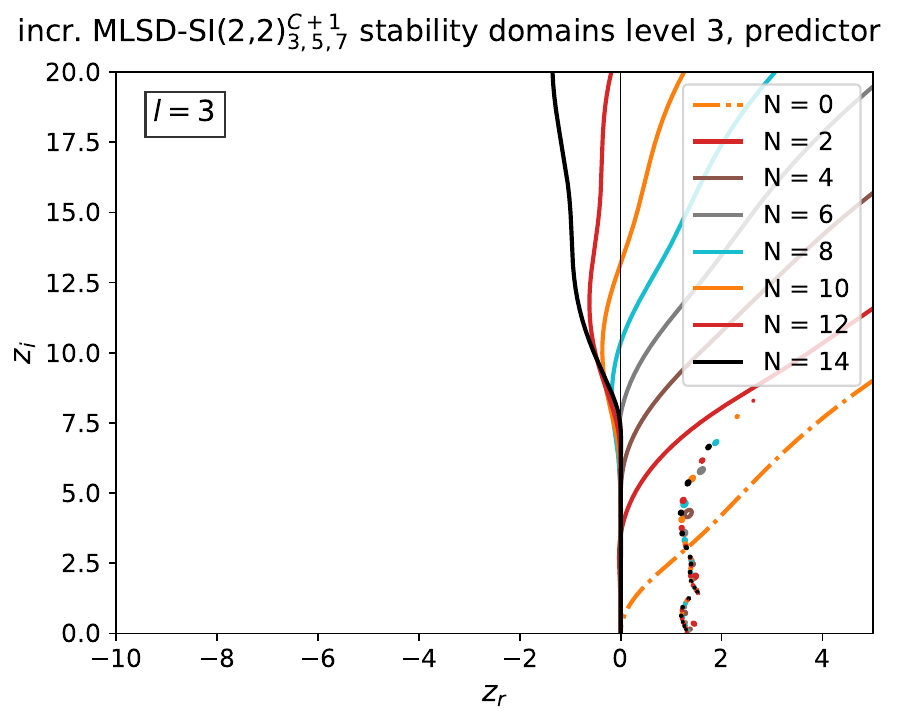}
        \caption{Top level $l = 3$ of MLSDC with SI(2) on all levels, predictor variant}
        \label{xxd}
    \end{subfigure}
    \vspace{1em} 
    \begin{subfigure}[b]{0.45\textwidth}
        \includegraphics[width=\textwidth]{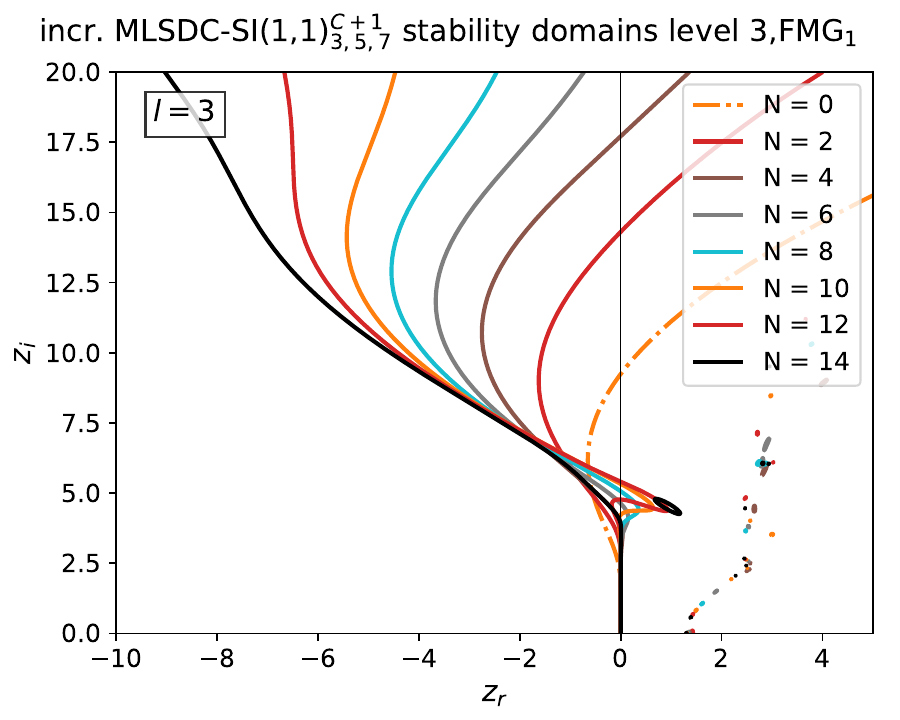}
        \caption{Top level $l = 3$ of MLSDC with SI(1) on all levels, FMG$_1$ variant}
        \label{xxe}
    \end{subfigure}
    \hfill
    \begin{subfigure}[b]{0.45\textwidth}
        \includegraphics[width=\textwidth]{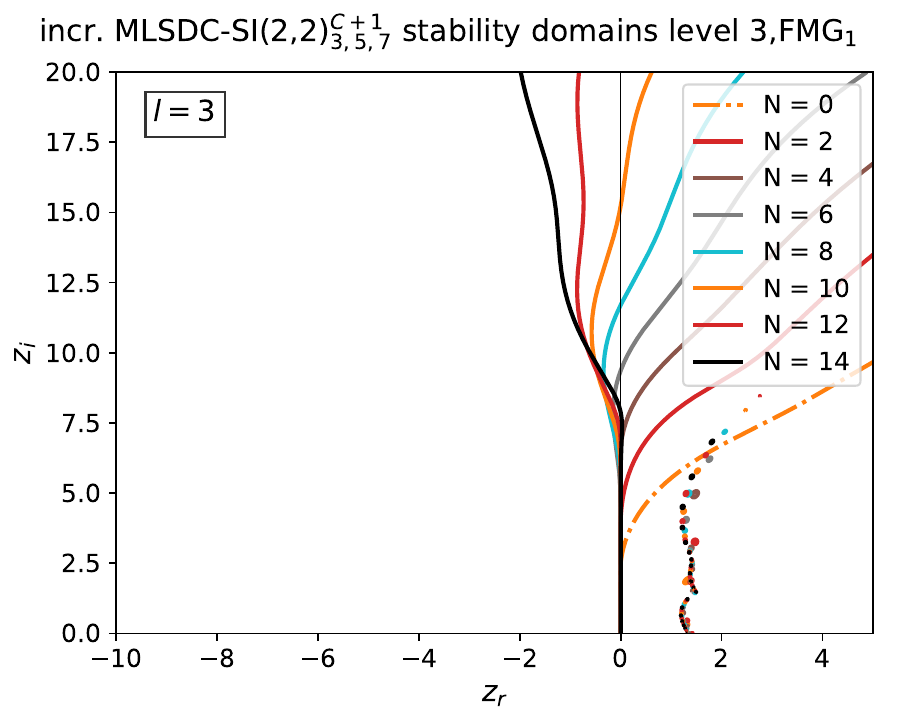}
        \caption{Top level $l = 3$ of MLSDC with SI(2) on all levels, FMG$_1$ variant}
        \label{xxf}
    \end{subfigure}
    \caption{Neutral stability curves for incremental MLSDC methods with Radau Right points and different time integrators shown after each odd V-cycle. \RII{FMG$_1$ denotes the full multigrid start strategy (algorithm \ref{alg:mlsdc:fmg}) with one cycle per added level, i.e., $C_{\mathsc{fmg}}=1$}}
    \label{fig:dahlquist_incr}
\end{figure}

\begin{figure}[H]
    \centering
    \begin{subfigure}[b]{0.45\textwidth}
        \includegraphics[width=\textwidth]{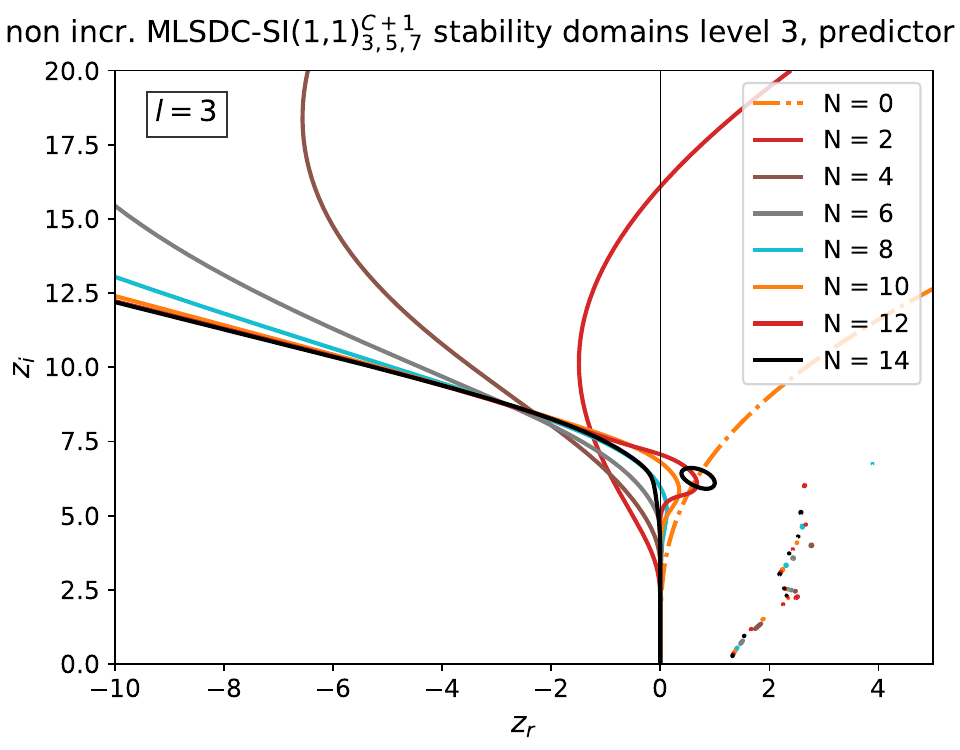}
        \caption{Top level $l = 3$ of MLSDC with SI(1) on all levels, predictor variant}
        \label{ni_xxa}
    \end{subfigure}
    \hfill
    \begin{subfigure}[b]{0.45\textwidth}
        \includegraphics[width=\textwidth]{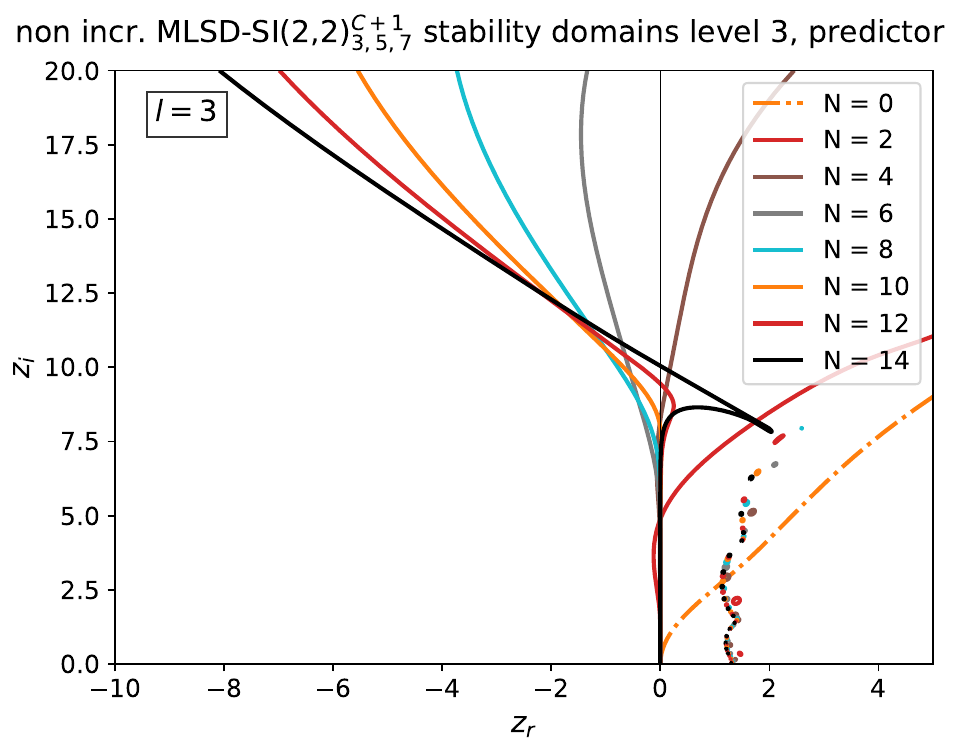}
        \caption{Top level $l = 3$ of MLSDC with SI(2) on all levels, predictor variant}
        \label{ni_xxb}
    \end{subfigure}
    \vspace{1em} 
    \caption{Neutral stability curves for non-incremental MLSDC methods with Radau Right points and different time integrators shown after each odd V-cycle}
    \label{fig:dahlquist_non-incr}
\end{figure}

The stability function can be obtained by computing the ratio of the amplitudes of two subsequent time steps, i.e., 
\begin{equation}
  R(z) = \frac{\hat u^{n+1}}{\hat u^n}
  \,,
\end{equation}
where ${z = \Delta t \lambda}$.
The set ${\{z \in \mathbb C; |R(z)| \le 1\}}$ is the stability domain of the method. 
The method is called $A$-stable if ${|R(z)| \le 1}$ for all $z$ with ${z_{\mathrm i} \le 0}$ and $L$-stable if, in addition, ${R = 0}$ for ${z \rightarrow \infty}$.
\begin{remark}
For conservation laws, 
${|z_i| = \Delta t |v| k}$ can be associated with the CFL number and
${-z_r = \Delta t \nu k^2}$ with the diffusion number 
of the fully discrete problem.
\end{remark}
\RII{The stability of the MLSDC method was investigated for different combinations of predictors and correctors introduced in Section ~\ref{ch:mlsdc}. In particular the IMEX Euler, SI(1) and SI(2) were explored.}
It is to be noted, that the predictor will also be used for the FMG starting strategy.
For the tests, one time step is performed with $\Delta t = 1$ and \RII{the} initial condition $\hat u(0) = 1$.
\RI{
For clarity, we briefly recall the parameters used throughout the stability analysis: 
$N$ denotes the number of fine-grid iterations, 
$C$ the number of cycles, 
$N_c$ the number of corrector sweeps on the coarse grid and 
$M_l$ the number of collocation nodes in time per level.
}
Figure \ref{fig:dahlquist_incr} shows representative neutral stability curves for MLSDC with the predictor and FMG$_1$ starting strategy using 3 levels with $M_l = (3,5,7)$. Radau Right collocation points were used and $C$ cycles performed. 
For example, MLSDC-SI(2,2)$^{C+1}_{3,5,7}$ means a MLSDC method with the SI(2) time integrator as a predictor and corrector, with 3 levels and a number of fine grid iterations $N = C+1$. This notation \RII{enables a direct comparison} to single-level SDC methods.
The number of coarse sweeps was set to $N_c = 2$.\\
Intuitively, one would only look at stability curves on the finest grid, but through FAS there are full solutions on all grids, which interact via the right-hand side.
One apprehension may therefore be that instabilities of the coarsest discretization will propagate to the fine grid.
An initial observation reveals that the stability of the multilevel SDC is slightly inferior to that of the single-level SDC \cite{TI_Stiller2024}, but the difference is marginal, and the multilevel method still performs remarkably well overall.
For a comparison to the single-level SDC methods, refer to Appendix \ref{sec:app:single-level-stability}.
It is further noticed, that the stability properties are rather similar on all levels, compare Figure \ref{xxa} with Figure \ref{xxb}. When comparing the SI(1) and SI(2) correctors to the single-level stability curves for $M = 3$ (Appendix \ref{sec:app:single-level-stability}, Figure \ref{fig:dahlquist_sdc}), the stability doesn't seem to be dictated by the stability of the coarsest level.
Compared to IMEX Euler, SI(1,1) yields a clear improvement in stability, which is further enhanced by the two-stage method SI(2,2).
The stability of the IMEX Euler method cannot keep up with the new SI time integrators and will no longer be considered from now on.
The predictor and FMG$_1$ starting variants hardly differ in stability with the predictor variant being slightly more stable possibly due to the interpolation within the FMG start, compare Figure \ref{xxc} and \ref{xxe} or \ref{xxd} and \ref{xxf}.
Increasing the number of cycles slowly erodes the stability, which complies with \cite{TI_Stiller2024}. 
A possible explanation may be given by the gradual removal of the artificial diffusion of SI(1) and SI(2) in both, single-level and multilevel SDC methods. This proceeds more rapidly for MLSDC, which may be the reason why the latter fails to preserve the $L$-stability of the related single-level method. 
A test with optimal SDC variants for each level, as identified in \cite[Table 1]{TI_Stiller2024}, \RII{did} not improve the stability here.
The comparison of the incremental and non-incremental MLSDC formulations for the SI time integrators indicates, that the incremental form is more stable, which underlines that they are not equivalent (compare Figure \ref{ni_xxa} with \ref{xxc} or Figure \ref{ni_xxb} with \ref{xxd}).
Moreover, the non-incremental form takes longer to compute because of the additional operations. Therefore, only the incremental formulation is considered in the following studies.
\begin{remark}
In application to space-time problems, a further improvement of stability is expected due to the coarsening in space, leading to lower $\CFL$ numbers or, respectively, imaginary parts $z_i$.
\end{remark}

\subsection{Accuracy}
\begin{figure}[H]
    \centering
    \begin{subfigure}[b]{0.45\textwidth}
        \includegraphics[width=\textwidth]{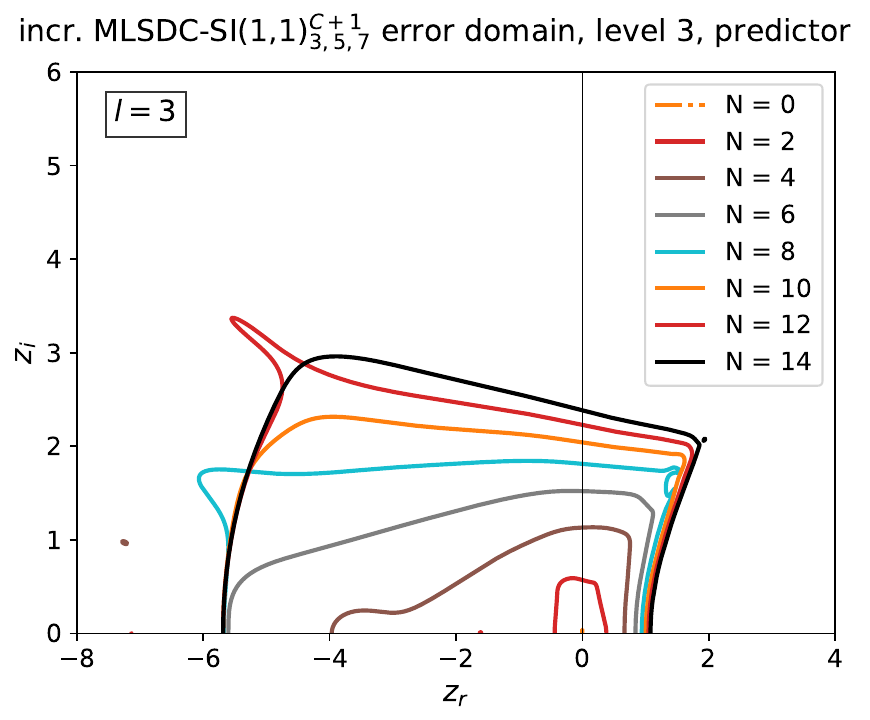}
        \caption{Top level $l = 3$ of incremental MLSDC with SI(1) fon all levels, predictor variant}
        \label{acc_xxa}
    \end{subfigure}
    \hfill
    \begin{subfigure}[b]{0.45\textwidth}
        \includegraphics[width=\textwidth]{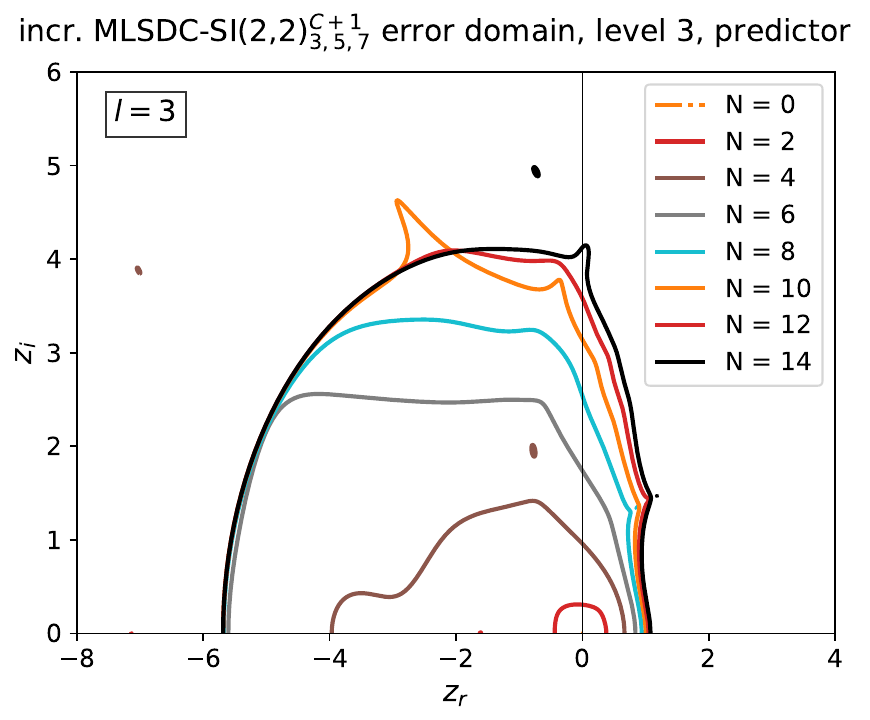}
        \caption{Top level $l = 3$ of incremental MLSDC with SI(2) on all levels, predictor variant}
        \label{acc_xxb}
    \end{subfigure}
    \vspace{1em}
    \begin{subfigure}[b]{0.45\textwidth}
        \includegraphics[width=\textwidth]{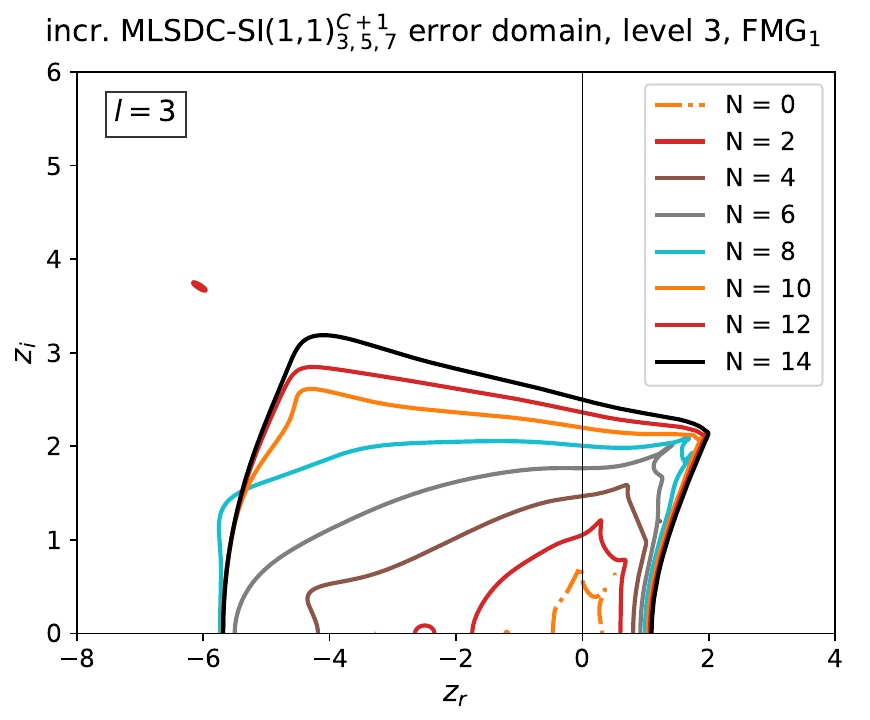}
        \caption{Top level $l = 3$ of incremental MLSDC with SI(1) on all levels, FMG$_1$ variant}
        \label{acc_xxc}
    \end{subfigure}
    \hfill
    \begin{subfigure}[b]{0.45\textwidth}
        \includegraphics[width=\textwidth]{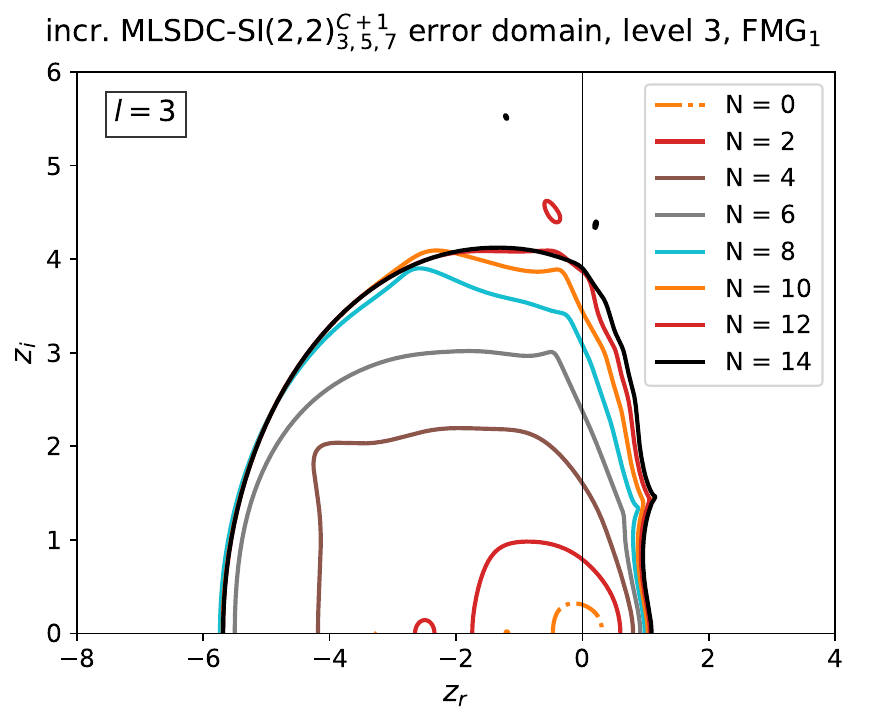}
        \caption{Top level $l = 3$ of incremental MLSDC with SI(2) on all levels, FMG$_1$ variant}
        \label{acc_xxd}
    \end{subfigure}
    \begin{subfigure}[b]{0.45\textwidth}
        \includegraphics[width=\textwidth]{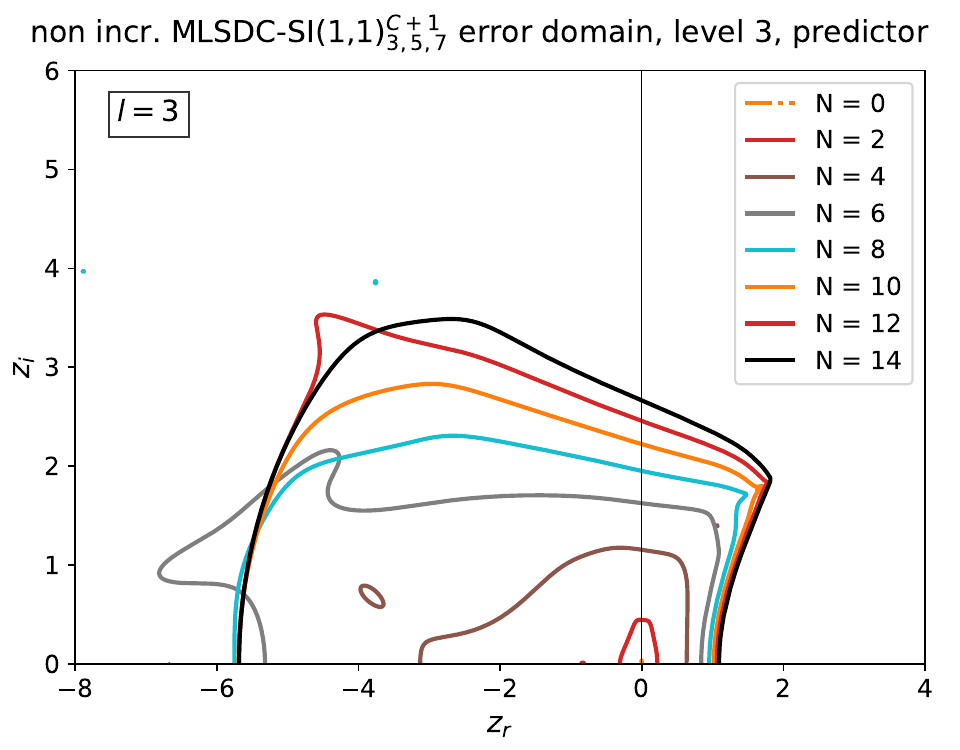}
        \caption{Top level $l = 3$ of non-incremental MLSDC with SI(1) on all levels, predictor variant}
        \label{acc_xxe}
    \end{subfigure}
    \hfill
    \begin{subfigure}[b]{0.45\textwidth}
        \includegraphics[width=\textwidth]{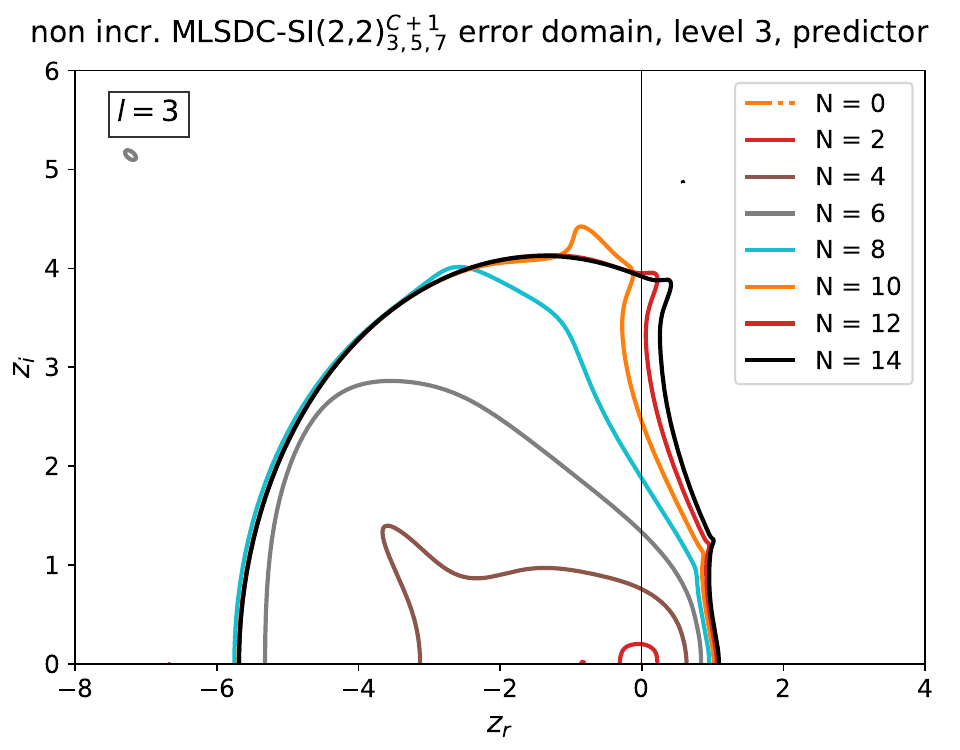}
        \caption{Top level $l = 3$ of non-incremental MLSDC with SI(2) on all levels, predictor variant}
        \label{acc_xxf}
    \end{subfigure}
    \caption{Accuracy curves ($\varepsilon = 10^{-6}$) for MLSDC methods with Radau Right points and different time integrators shown after each odd V-cycle for the incremental and non-incremental formulation. \RII{FMG$_1$ denotes the full multigrid start strategy (algorithm \ref{alg:mlsdc:fmg}) with one cycle per added level, i.e., $C_{\mathsc{fmg}}=1$}}
    \label{fig:dahlquist_incr_accuracy}
\end{figure}

Since increasing the number of cycles appears to adversely affect stability, a natural question arises: how many cycles does MLSDC require to achieve convergence?
Figure~\ref{fig:dahlquist_incr_accuracy} presents accuracy plots for both the incremental and non-incremental methods, showing only the SI(1) and SI(2) time integrators. These methods were selected due to their superior stability. The results correspond to an error tolerance of $\varepsilon = 10^{-6}$ applied uniformly across all levels.
Unlike the single-level SDC methods where lines between iterations are expected to be evenly spaced (Appendix \ref{sec:app:single-level-stability}, Figure \ref{fig:dahlquist_accuracy_sdc}), indicating a uniform increase in order by one per sweep, the plots show that the spacing between lines varies for MLSDC, with convergence initially occurring rapidly before slowing down.
Theoretically, \RII{a single-level} SDC \RII{method} with Radau Right points and $M=7$ would require 14 sweeps (Appendix \ref{sec:app:single-level-stability}, Figure \ref{fig:dahlquist_accuracy_sdc}), whereas MLSDC requires only 9 cycles to achieve the same accuracy for SI(2,2). 
A comparison of the incremental and non-incremental methods (see Figures \ref{acc_xxa} with \ref{acc_xxe} and Figures \ref{acc_xxb} with \ref{acc_xxf}) reveals only minor differences, with the incremental approach remaining preferable due to its stability advantages.
The FMG$_1$ variant is observed to perform slightly better compared to the predictor method, particularly in the initial phases of convergence (compare Figure \ref{acc_xxa} with \ref{acc_xxc} and Figure \ref{acc_xxb} with \ref{acc_xxd}), although both converge after the same number of iterations. 
The detailed comparison of the convergence between MLSDC and traditional SDC will be continued in subsequent sections.

%% file: spatial-discretization.tex
\section{Spatial discretization}
\label{ch:spatial_discretization}

As the investigations now encompass partial differential equations, it is essential to introduce spatial discretization.
\RII{In this study, a} spectral element method \RII{is used}, subdividing the domain $\Omega$ into a set of discrete elements \RII{$\Omega^e$. The} subdivision yields a computational domain denoted by ${\Omega_h = \cup\Omega^e}$.
Within this domain, we define a set of element boundaries $\Gamma_h$. These are further divided into $\Gamma_h^{\mathsc i}$, representing the internal boundaries, and $\Gamma_h^{\d}$, indicating the external ones.
The approximate solution $\M u_h(x,t)$ is sought within a specific function space at a given time $t$.
This function space is defined as
\begin{align}
\mathbb U_h = \left\{
    \M u_h \in \bigl[\mathbb L^2(\Omega)\bigr]^d  : 
    \M u_h|_{\Omega^e} \in  \bigl[\mathbb P_{\! P}(\Omega^e)\bigr]^d 
    \quad 
    \forall \Omega^e \in \Omega
  \right\}.
\end{align}
Here, $\mathbb L^2$ denotes the space of square-integrable functions, and $\mathbb P_{\! P}$ refers to the space of polynomials with a degree of at most $P$.
To formulate the problem, weighting functions $\M w_h$ are introduced. Mathematically, these weighting functions belong to the same function space, ${\M w_h \in \mathbb U_h}$. For physical consistency, $\M w_h$ must have dimensions that are reciprocal to those of $\M u_h$.
At the boundary between adjacent elements, the average and jump operators are given by
\begin{align}
\avg{\M u_h} = \frac{\M u_h^- + \M u_h^+}{2} \,, \quad \\
\quad \jmp{\M u_h} = n^- \M u_h^- + n^+ \M u_h^+.
\end{align}
Here, $\M u_h^{\pm}$ are the traces of the corresponding element solutions and $n^{\pm}$ the normals. 
The discontinuous Galerkin formulation is composed of discrete functionals $\mathcal M$ for projection to weighting functions,
$\Fc$ for convection,
$\Fd$ for diffusion and
$\Fs$ for sources.
\RII{
The discrete form of the collocation method \eqref{eq:cm:incremental} becomes
\begin{equation}
  \label{eq:collocation:radau-iia:full}
  \PO(\M w_h,\M u_{h,m}) = \PO(\M w_h,\M u_{h,m-1})
      + \Delta t \sum_{i=1}^{M} w^{\mathsc{nn}}_{m,i} \F(\M w_h, \M u_{h,i}, t_m),
  \quad
  m = 1,\dots M
  \,,
\end{equation}
where} the right hand side can be written in the form
\begin{equation}
\mathcal F(\M w_h, \M u_h, \M A_{\mathrm d,h})
= \Fc(\M w_h, \M u_h)
+ \Fd(\M w_h, \M u_h, \M A_{\mathrm d,h})
+ \Fs(\M w_h, \M u_h, t) \,.
\end{equation}
For details it is referred to \cite{TI_Stiller2024}.
The diffusion matrix $\M A_{\mathrm d,h}$ can be tailored to fit the Euler, SI(1) and SI(2) time integrators introduced above.
Moreover, it can be supplemented by an artificial viscosity in order to capture shocks and other discontinuities.
In some of the numerical experiments discussed below, we used the artificial viscosity proposed by \citet{SE_Persson2006a} (see also Section 4.3 in \cite{TI_Stiller2024}).

\subsection{Multilevel discretization.}
\RII{The previously introduced spectral element discretization is extended to a multilevel framework by coupling it with the MLSDC time integration scheme.}
Consider a space-time slab (or a space-time element) for a one-level discretization with a domain described by ${Q_h = \Omega_h \times T_h}$.
The temporal domain $T_h$ is similarly divided into intervals $T^n$, each of length $\Delta t$ and containing $M$ subintervals as described in the previous sections. As a result, the combined space-time elements are represented by 
\begin{align}
Q^{e,n} = \Omega^e \times T^n.
\end{align}
In a multilevel discretization, this structure is extended to multiple levels and denoted by ${\M Q_h = [ Q_l ]}$, where ${Q_l = \Omega_l \times T_l}$ is the domain for the levels ${l = 1 \dots L}$.
At each level, the spatial domain is represented as ${\Omega_l = \cup \Omega^e_l}$, consisting of elements with length $\Delta x_l$ and polynomial degree $P_l$.
\RII{The multilevel spatial domain can be expressed as ${\M \Omega = [ \Omega_l ]}$.}
The temporal domain $T_l$ contains intervals $T^n_l$ of length $\Delta t_l$ and degree $M_l$.
This results in space-time elements at each level denoted by 
\begin{align}
Q^{e,n}_l = \Omega^e_l \times T^n_l.
\end{align}

\subsection{MLSDC transfer operators.}
The discrete solution of level $l$ is denoted as $\M u_l(x,t)$, where the subscript $h$ has been dropped for better readability.
For transferring information between subsequent levels, the operators 
$\mathcal I_l$, $\mathcal P_l$ and $\mathcal R_l$ need to be extended from time to space-time operands.
This is achieved by defining corresponding spatial operators and applying them either before or after their temporal counterparts.
As the spatial approximation is composed of polynomial basis functions in each element, the design principles are largely identical to the temporal case.
One exception appears in the case of 2:1 $h$-coarsening, where two adjacent fine elements are fused into one on the next coarser level.
This coarsening technique may include optional adjustments for discontinuities, such as removing jumps via linear blending (default) or coefficient averaging at interfaces.
If embedded interpolation is used for projection $\mathcal P_l$, each fine element interpolates to the corresponding part of the coarse element. 
Contributions to the center of the latter are averaged.
With $L^2$-projection, both fine elements contribute to the whole coarse element.
The spatial restriction operator $\mathcal R_l$ is again defined as the transpose of the interpolation $\mathcal I^t_{l+1}$.
\RI{
In summary, the framework provides the possibility of employing either $h$– or $p$–refinement (or a combination thereof) in space, while in time only $p$–refinement is considered across the levels.
}

%% file: mlsdc-basics.tex
\section{Preliminary studies}
\label{ch:mlsdc_basics}

\begin{figure}[H]
    \centering
    \begin{subfigure}[b]{0.45\textwidth}
        \centering
        \includegraphics[width=\textwidth]{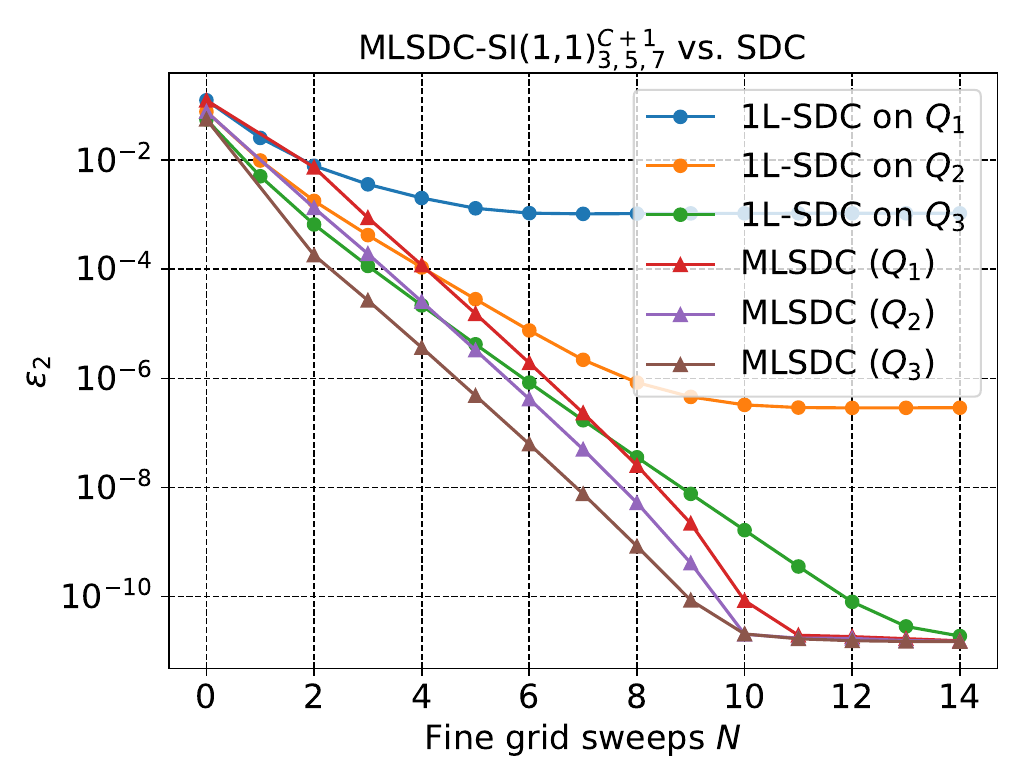}
        \caption{Test 1 - embedded interpolation}
        \label{fig:test_1_int}
    \end{subfigure}
    \hfill
    \begin{subfigure}[b]{0.45\textwidth}
        \centering
        \includegraphics[width=\textwidth]{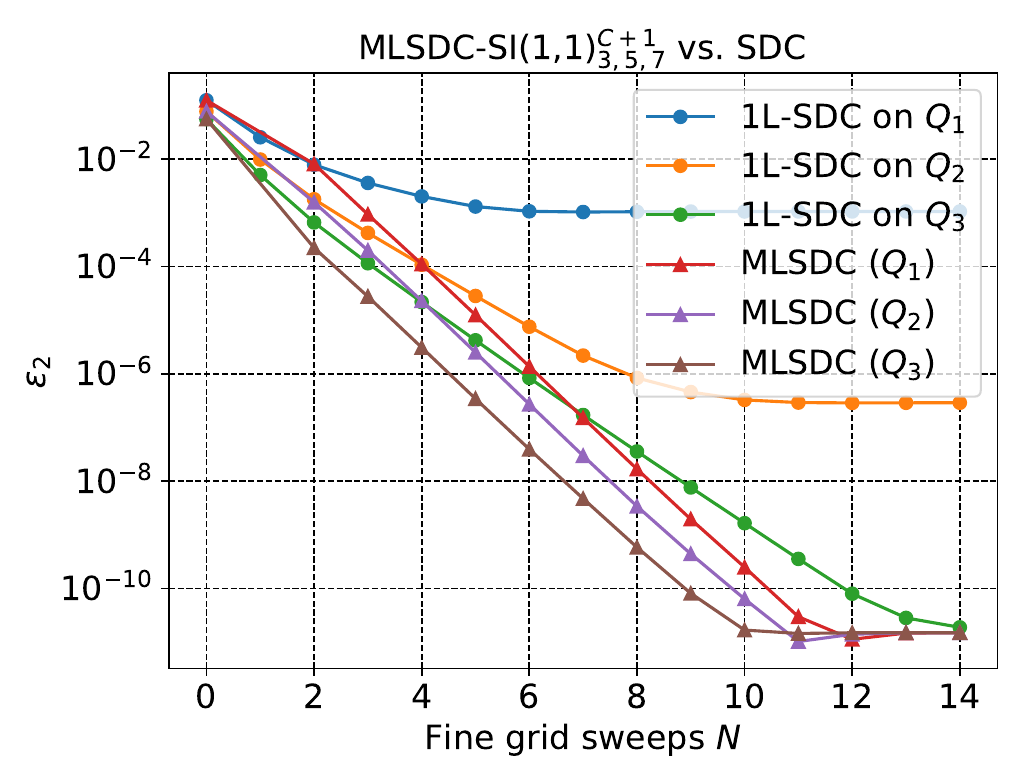}
        \caption{Test 1 - $L^2$-projection}
        \label{fig:test_1_l2p}
    \end{subfigure}
    \caption{MLSDC V-cycle vs. single-level SDC corrector, number of coarse sweeps $N_{c} = 2$, predictor variant}
    \label{fig:test_1}
\end{figure}

\begin{figure}[H]
        \centering
    \begin{subfigure}[b]{0.45\textwidth}
        \centering
        \includegraphics[width=\textwidth]{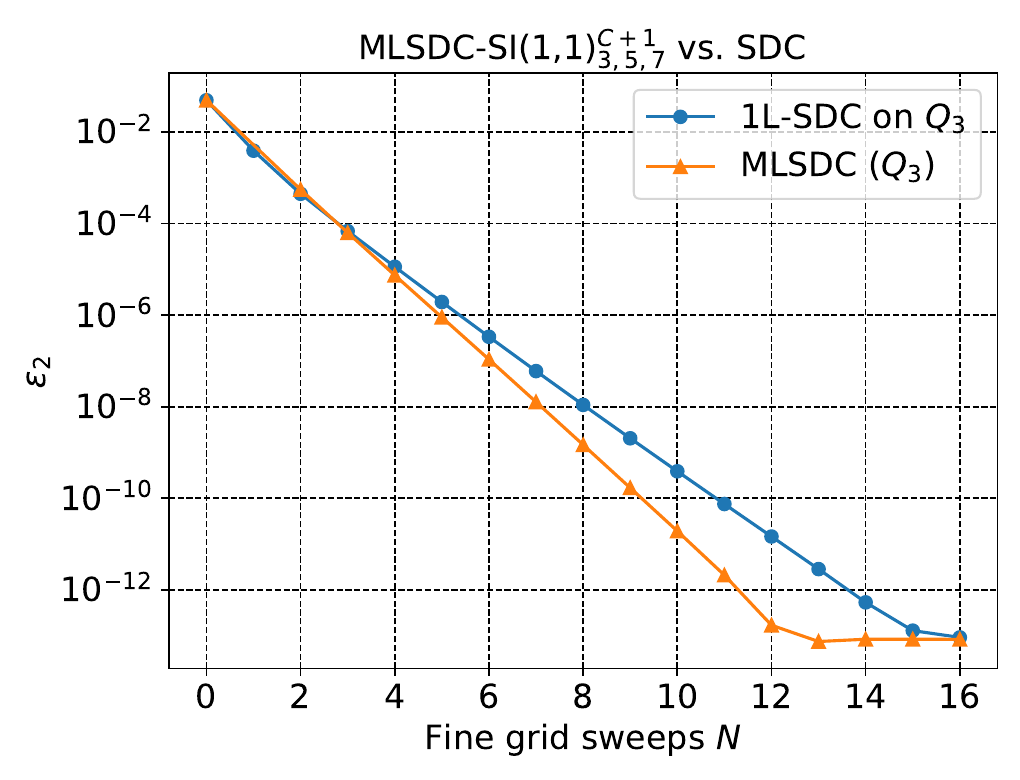}
        \caption{Test 2 - number of coarse sweeps $N_{c} = 2$, stronger coarsening}
        \label{fig:test2-248}
    \end{subfigure}
    \hfill
    \begin{subfigure}[b]{0.45\textwidth}
        \centering
        \includegraphics[width=\textwidth]{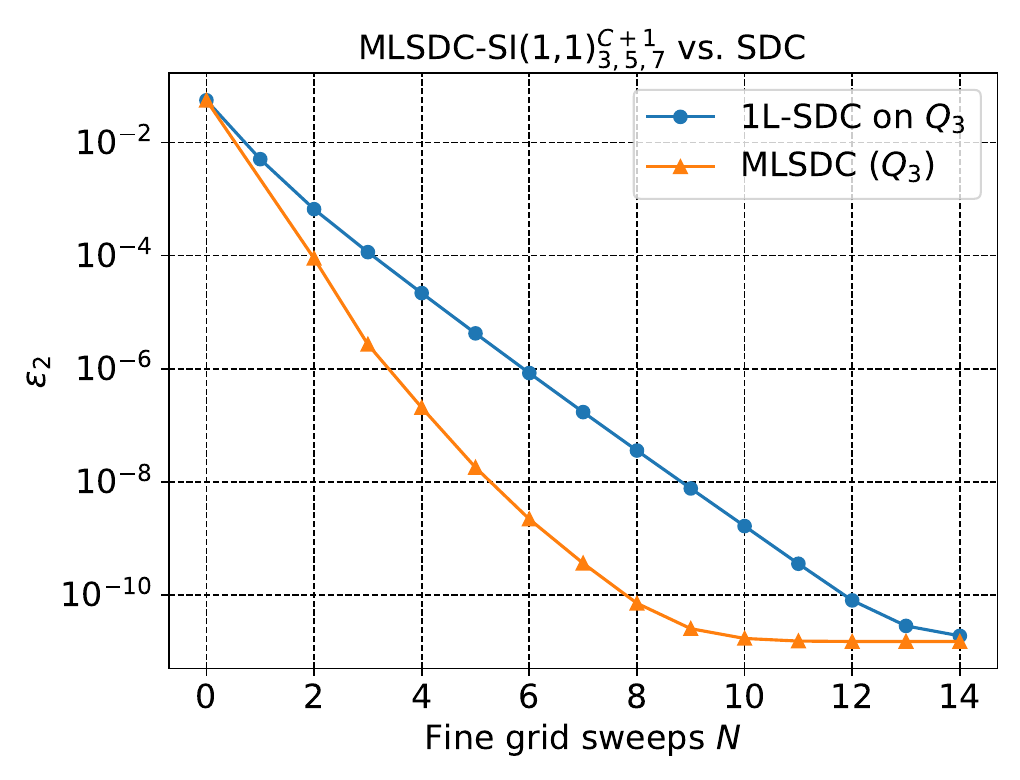}
        \caption{Test 3 - number of coarse sweeps $N_{c} = 3$}
        \label{fig:test3-nc3}
    \end{subfigure}
    \caption{MLSDC V-cycle vs. single-level SDC corrector, predictor variant, embedded interpolation, top level $l = 3$}
    \label{fig:test_23}
\end{figure}

\begin{figure}[H]
    \centering
  \includegraphics[width=0.7\textwidth]{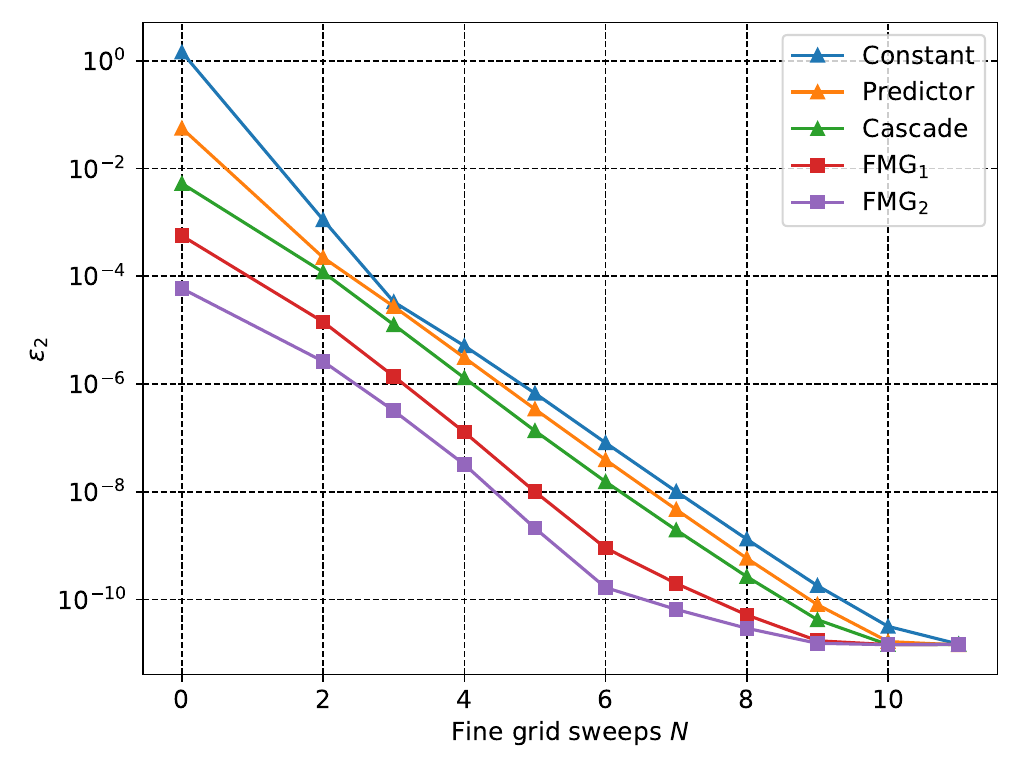}
    \caption{Test 4 - errors for different starting strategies to the MLSDC V-cycle, top level $l = 3$ }
    \label{fig:test_4}
\end{figure}

In this section, the basic functionality and key parameters as well as \RII{the} performance of the MLSDC algorithm \RII{are} investigated by first examining it with a simple convection-diffusion example \eqref{eq:conv-diff}.
To further simplify the initial analysis, only temporal coarsening is considered and one time step performed in $T =[0,0.1]$ with the time step $\Delta t = 10^{-2}$ being equivalent to the length of the time slab.
Periodic boundary conditions are used for the problem. 
The spatial domain is $x \in \Omega = [0, 1]$ on all levels and the parameters $\nu = 2 \cdot 10^{-2}$ and $v = 1$ are set.
The exact solution to the problem is given by:
\begin{equation}
  \label{eq:wave-packet}
  u(x,t) = \sum_{i=1}^7 a_i \sin\big( \kappa_i(x - s_i - vt) \big)
                        e^{-\kappa_i^2 \nu t}
  .
\end{equation}
Equation~\eqref{eq:wave-packet} defines a wave packet comprising various wave numbers $\kappa_i$, amplitudes $a_i$, and phase shifts $s_i$, as specified in Table \ref{tab:wave-package}.
\begin{table}[ht]
  \begin{tabular}{lccccccc} 
  \toprule
  $i$        & $1$    & $2$    & $3$     & $4$     & $5$     & $6$     & $7$     \\\midrule
  $\kappa_i$ & $2\pi$ & $6\pi$ & $10\pi$ & $14\pi$ & $18\pi$ & $24\pi$ & $30\pi$ \\
  $a_i$      & $1.00$ & $1.50$ & $1.80$  & $1.70$  & $1.50$  & $1.30$  & $1.15$  \\
  $s_i$      & $0.00$ & $0.05$ & $0.10$  & $0.15$  & $0.20$  & $0.30$  & $0.18$  \\
  \bottomrule
  \end{tabular}
    \caption{Wave package coefficients}
    \label{tab:wave-package}
\end{table}
The discretization parameters for space and time are recorded in Table \ref{tab:test_124} with $N_e$ being the number of elements and $n_t$ the number of time steps.
The error should primarily be influenced by the time discretization, ensuring that spatial discretization effects remain minimal.
For the first test, a single-level SDC corrector with three different discretizations is compared with a 3 level MLSDC V-cycle with the number of sweeps for the coarse solution $N_{c} = 2$.
\begin{table}
\centering
\hspace{-0.5cm}
\begin{minipage}{.55\linewidth}
\centering
\begin{tabular}{cccccc}
\toprule
 & \multicolumn{2}{c}{Space} & \multicolumn{2}{c}{Time} \\
\cmidrule(r){2-3} \cmidrule(r){4-5}
 & $N_{e,l}$ & $P_l$ & $n_{t,l}$ & $M_l$ \\
\midrule
Level 1 \RI{($Q_1$)} & 32 & 15 & 1 & 3 \\
Level 2 \RI{($Q_2$)} & 32 & 15 & 1 & 5  \\
Level 3 \RI{($Q_3$)} & 32 & 15 & 1 & 7  \\
\bottomrule
\end{tabular}
\caption{Discretization parameters for tests 1,3 and 4}
\label{tab:test_124}
\end{minipage}%
\hspace{-1.5cm}
\begin{minipage}{.55\linewidth}
\centering
\begin{tabular}{cccccc}
\toprule
 & \multicolumn{2}{c}{Space} & \multicolumn{2}{c}{Time} \\
\cmidrule(r){2-3} \cmidrule(r){4-5}
 & $N_{e,l}$ & $P_l$ & $n_{t,l}$ & $M_l$ \\
\midrule
Level 1 \RI{($Q_1$)} & 32 & 15 & 1 & 2 \\
Level 2 \RI{($Q_2$)} & 32 & 15 & 1 & 4  \\
Level 3 \RI{($Q_3$)} & 32 & 15 & 1 & 8  \\
\bottomrule
\end{tabular}
\caption{Discretization parameters for test 2}
\label{tab:test_3}
\end{minipage}
\end{table}
Figure \ref{fig:test_1} shows the $L^2$ error of the variable $\NM{u}_l$ for each level and after each iteration.
The test is done using embedded interpolation (Figure \ref{fig:test_1_int}) for the projection operator or $L^2$-projection (Figure \ref{fig:test_1_l2p}).
The iteration count is equivalent to the number of fine grid SDC sweeps \RII{($N=C+1$)}. So for the single-level SDC it is the number of sweeps and for the MLSDC method it is the number of cycles + 1, to \RII{account for} the post-sweeping. Zero iterations show the error before applying the SDC corrector or entering the V-cycle.
Both methods start at the same point with the predictor solution.
The SI(1) method described in Section \ref{ch:sdc} was used for the predictor and corrector.
All tests were performed using Radau-Right and Lobatto points as collocation points in time. Overall, Radau Right points yielded better accuracy than Lobatto points throughout all studies, and therefore all results in this document are shown for Radau Right points.
Firstly, it can be observed that the MLSDC algorithm converges with fewer iterations and saves SDC sweeps on the finest grid, which was expected.
Another observation is that all levels in the V-cycle converge towards the solution on the finest grid, which means that the FAS-RHS works as intended.
Note that the SDC Corrector alone for the coarser discretization converges to a higher error.
\RII{Further,} MLSDC achieves a faster convergence rate, especially in the first cycle. In contrast, standalone SDC maintains a constant convergence rate of one order per sweep. This behaviour aligns with the accuracy analysis presented in Section \ref{ch:stability}.
In this test, the MLSDC method converges after 9 cycles (10 fine grid sweeps), whereas the single-level SDC needs 14 sweeps.
The difference between the $L^2$-projection and the embedded interpolation is hardly visible for this test case. Since the embedded interpolation still introduces smaller projection errors, it will be the default choice for upcoming tests. 
Using the restriction operator \eqref{eq:restriction} for variable restriction produced inferior results and is therefore not considered.\\
Test 2 repeats the first test but with a stronger coarsening in time, compare Table \ref{tab:test_3}. The results for the top levels are shown in Figure \ref{fig:test2-248}.
The single-level SDC requires 16 sweeps to reach convergence, while the MLSDC V-cycle converges after 12 cycles, corresponding to 13 fine grid sweeps. 
Although MLSDC still converges more efficiently, its advantage over the single-level SDC is noticeably less pronounced in this case, indicating that the coarsening in the time discretization between the levels cannot be too permissive.
Following the insights of \cite{MG_Speck2014}, a more restrictive temporal coarsening of $M_{l-1} = M_l/2 + 1$ will be approached from now.\\
The third test repeats test 1 with a different number of coarse sweeps $N_c = 3$.
The parameter $N_{c}$ significantly impacts the convergence rate of the MLSDC method. As $N_{c}$ increases, the convergence accelerates (compare e.g. the top levels in Figure \ref{fig:test_1} and Figure \ref{fig:test3-nc3}).
However, it also boosts computational cost and may generate instabilities, necessitating an optimal balance for an efficient runtime.
Note that the performance gain of MLSDC is limited by using the predictor-based start. However, this choice was necessary to ensure a fair comparison with single-level SDC.
Test 1 is again repeated for the starting strategies presented in Section \ref{ch:mlsdc}, see Figure \ref{fig:test_4}. Only errors on the finest level are displayed after each iteration.
The constant-value initialization approach yields the worst results, converging one iteration later than the predictor variant.
The Cascade variant performs only slightly better than the predictor start, converging after 9 cycles as well.
The FMG variants converge the fastest with both FMG$_1$ and FMG$_2$ variants almost finishing after 7 cycles.
Runtime savings with the FMG$_1$ variants over the Cascade can be anticipated, as the performed first incomplete V-cycles of the FMG starting variants occur solely between the two coarsest levels. Consequently, the reduction in full V-cycles should compensate for the additional computational cost incurred during the startup phase.
A runtime analysis would be relevant at this point; however, measurements taken over a single time step may be neither reliable nor representative and the code is not trimmed for best performance. Thus, a more extensive runtime investigation will be conducted in future work with multiple time steps.
For comparable accuracy, the Cascade and constant approaches are expected to be slower than FMG.
Based on the observations and analysis, a set of default parameters is established for the subsequent tests. The number of coarse correction sweeps is set to  $N_c = 2 $, as this configuration provides good convergence properties while maintaining stability. 
For projection, embedded interpolation is chosen over $L^2$-projection due to its slightly favourable accuracy. 
Regarding the MLSDC starting strategy, the FMG$_1$ approach is selected for two-level and multilevel configurations. Note that the FMG$_1$ approach is equivalent to the Cascade strategy for a two-level configuration.
The incremental formulation is preferred over the non-incremental method, as it proved to be more stable. The time integration method is fixed to SI(2), based on insights from the Dahlquist-type analysis.
Thus, unless explicitly stated otherwise, all subsequent tests are conducted using $N_c = 2$, embedded interpolation as the projection method, the Cascade strategy or FMG$_1$ for initialization, and the incremental formulation.
Furthermore, post-sweeping is enabled by default.

%% file: numerical-studies.tex
\section{Numerical Experiments}
\label{ch:numerical_experiments}

The accuracy and stability of the MLSDC method are extensively evaluated through numerical experiments on linear convection-diffusion, Burgers, and compressible flow problems in this section.
Two central questions arise in the context of MLSDC: First, does the method retain the same stability properties as single-level SDC when applied to partial differential equations - or potentially even improves stability due to the combined space-time coarsening?
Second, how rapidly does MLSDC converge in practice for PDEs with space-time coarsening, and to what extent does it reduce the number of iterations on the fine grid compared to standard SDC?
For the upcoming tests, the $\CFL$ number is defined as
\begin{equation}
\label{eq:cfl}
  \CFL = \frac{\Delta t\,\lambda_{\mathrm c, \max}}{\Delta x}
  \,,
\end{equation}
\RI{where $\lambda_{\mathrm c, \max}$ is the maximum magnitude of the eigenvalues of $\Ac$ in $\Omega_h$}. The characteristic length scale $\Delta x$ is following the approach in \cite[Sec. 7.3.3]{SE_Canuto2011a} and chosen as
\begin{equation}
  \Delta x = \frac{\Delta x^e}{2 \delta} 
  \,.
\end{equation}
There, ${\delta \sim P^2}$ denotes the largest eigenvalue of the convection problem with unit velocity and one-sided Dirichlet boundary conditions in the standard element ${[-1,1]}$ \RI{and $\Delta x^e$ denotes the element width}.
The $\CFL$ number is determined from the initial conditions and on the finest grid, unless otherwise stated.

\subsection{Convection-diffusion}
\label{sec:num_convdiff}
\begin{figure}[H]
    \centering
        \includegraphics[width=0.8\textwidth]{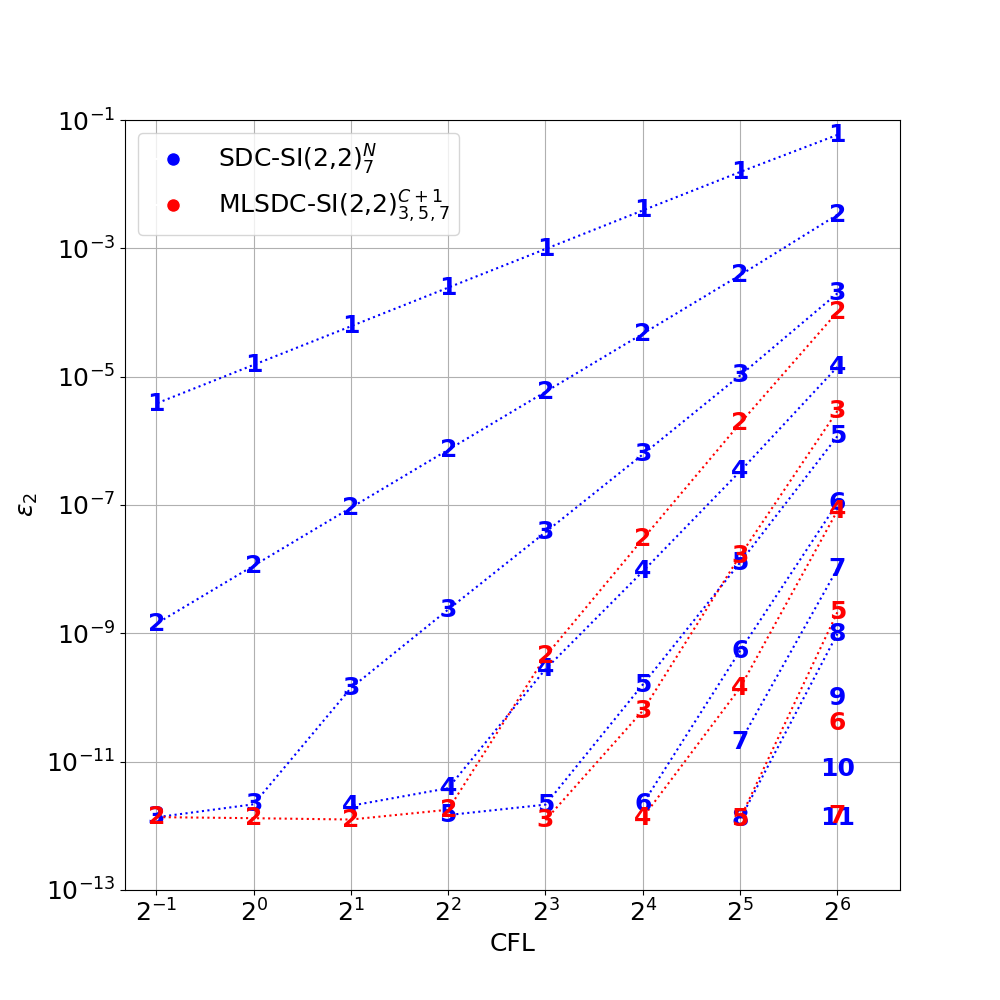}
        \label{cd_nu1e-3_xxa}
    \caption{$L_2$ errors over $\CFL$ numbers SDC$_7^N$ vs. MLSDC$_{3,5,7}^{C+1}$ for convection-diffusion and $\nu=10^{-3}$, top level, markers indicate the number of fine grid sweeps}
    \label{fig:convdiff_nu1e-3}
\end{figure}

The numerical experiments in this section build upon the tests presented in \cite[Sec. 6]{TI_Stiller2024}.
The initial set of test problems focuses on the convection-diffusion equation \eqref{eq:conv-diff}, with the exact solution given by the same wave packet as discussed in Section \ref{ch:mlsdc_basics}.
The tests are performed in the periodic domain ${\Omega = [0,1]}$ and a time interval ${[0,5]}$, i.e. five periods of the longest wave.
The discretization parameters for the first test can be taken from Table \ref{tab:num_experiments_convdiff_1}, the number of coarse sweeps $N_c = 2$. We first consider a constant velocity of $v=1$ and diffusivity $\nu = 10^{-3}$.
\begin{table}[H]
\centering
\centering
\begin{tabular}{cccc}
\toprule
 & \multicolumn{2}{c}{Space} & \multicolumn{1}{c}{Time} \\
\cmidrule(r){2-3} \cmidrule(r){4-4}
 & $N_{e,l}$ & $P_l$ & $M_l$ \\
\midrule
Level 1 & 16 & 15 & 3 \\
Level 2 & 32 & 15 & 5 \\
Level 3 & 64 & 15 & 7 \\
\bottomrule
\end{tabular}
\caption{Discretization parameters for the convection-diffusion tests}
\label{tab:num_experiments_convdiff_1}
\end{table}
Figure~\ref{fig:convdiff_nu1e-3} compares the single-level SDC method using the fine grid discretization (level 3) with MLSDC. 
The $L_2$ error of the fine grid is shown after each cycle of the MLSDC method (red markers) and after each sweep for SDC (blue markers) for different $\CFL$ numbers. The markers specify the number of fine grid sweeps. The $C+1$ notation indicates $C$ V-cycles + 1 fine grid sweep to account for the post-sweeping.
The error is displayed after each iteration, since there is no theoretical maximum number of iterations required for the MLSDC method to reach convergence. 
The $\CFL$ number on the finest grid is initially set to $64$ and is subsequently reduced by a factor of 2 down to $1/2$ \RII{, i.e. the time step size is reduced}.
As expected, the errors decrease after each iteration, whether it corresponds to a sweep in SDC or a cycle in the MLSDC method. 
There are clear savings in fine grid iterations for MLSDC. For $\CFL$ numbers of $0.5$ to $4$, it requires only one full cycle to converge and hence two fine grid sweeps. For $\CFL = 64$, MLSDC converges after 6 cycles, while the SDC method requires 11 sweeps. 
The advantage of MLSDC becomes more pronounced as the number of needed iterations to converge increases. 
As a result, the benefits of MLSDC are more apparent for higher $\CFL$ numbers. No stability issues were encountered in any of the cases tested with $\nu = 10^{-3}$.

\begin{figure}
    \centering
        \includegraphics[width=0.8\textwidth]{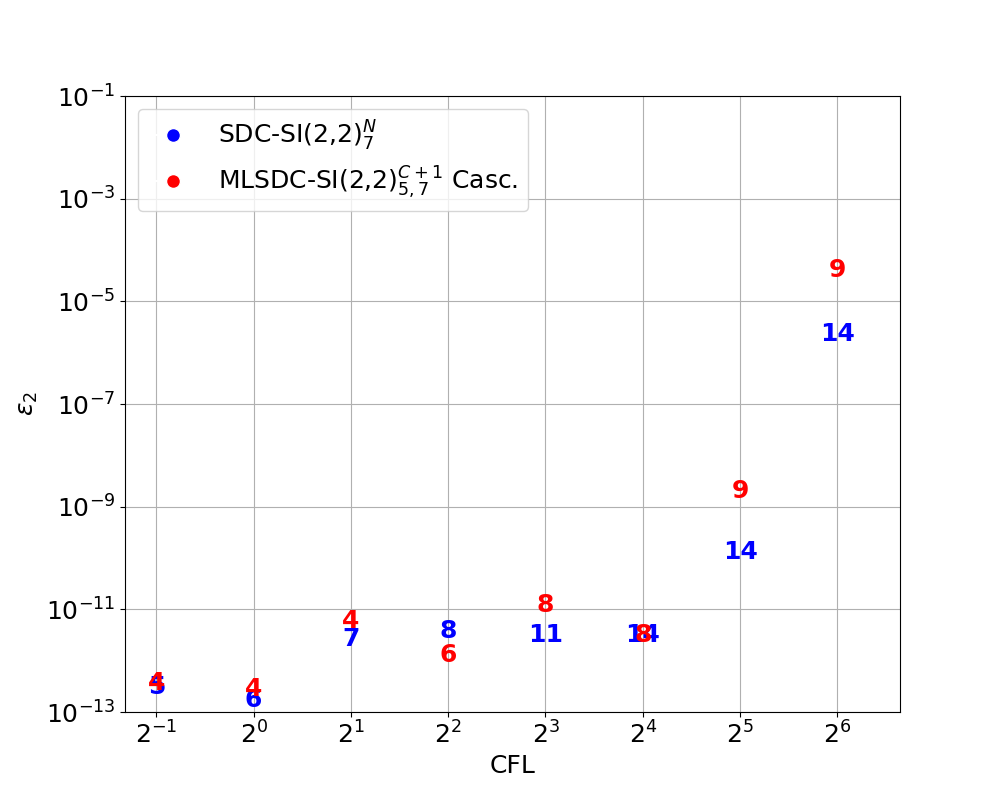}
        \caption{$L_2$ errors over $\CFL$ numbers SDC$_7^N$ vs. MLSDC$_{5,7}^{C+1}$ for convection-duffusion and $\nu=0$, top level, markers indicate the number of fine grid sweeps}
        \label{fig:convdiff_nu0}
\end{figure}

The second test considers the case of pure convection $\nu = 0$ in the time interval ${[0,5]}$. 
The discretization parameters are identical to the previous test.
Figure \ref{fig:convdiff_nu0} compares the single-level SDC method to a two-level method (level 2 and 3) with the Cascade start.
For clarity, only the iteration count at which the method converges for a given $\CFL$ number is shown. A test case is considered converged when the relative change in the $L_2$ error between two iterations falls below $10$ percent. This criterion is applied consistently across all subsequent experiments.
The general results for $\nu = 0$ are consistent with those observed for $\nu = 10^{-3}$, although more iterations are generally required to achieve convergence.
Across all $\CFL$ numbers, MLSDC exhibits faster convergence. For instance, at $\CFL = 32$ and $\CFL = 64$, it converges after 8 cycles compared to 14 sweeps required by the single-level SDC method. However, for larger time steps, the method tends to converge to a higher terminal error, despite the reduced number of iterations.
Stability issues are observed for both methods, and these become more pronounced in MLSDC as the number of coarse-grid levels, $N_c$, increases.
Applying coarsening in both time and space mitigates these issues to some extent, as the coarse-grid $\CFL$ number is reduced.
Up to a $\CFL$ number of 16, the method remains stable. Beyond this threshold, it becomes unstable if more than 9 cycles are performed for this problem.
Another key observation is that the resolution gap between levels must not be too large. Otherwise, transfer errors between levels increase significantly, adversely affecting both stability and accuracy. 
In particular, the coarsest level is critical for determining the overall performance of the MLSDC method.
Using a 3-level FMG$_1$ variant with the discretization parameters from Table~\ref{tab:num_experiments_convdiff_1} results in faster convergence for smaller $\CFL$ numbers. However, for higher $\CFL$ values, stability and accuracy degrade. This issue will be investigated in more detail in the subsequent test cases.

\subsection{Burgers equation}

\begin{figure}[H]
    \centering
    \includegraphics[width=0.8\textwidth]{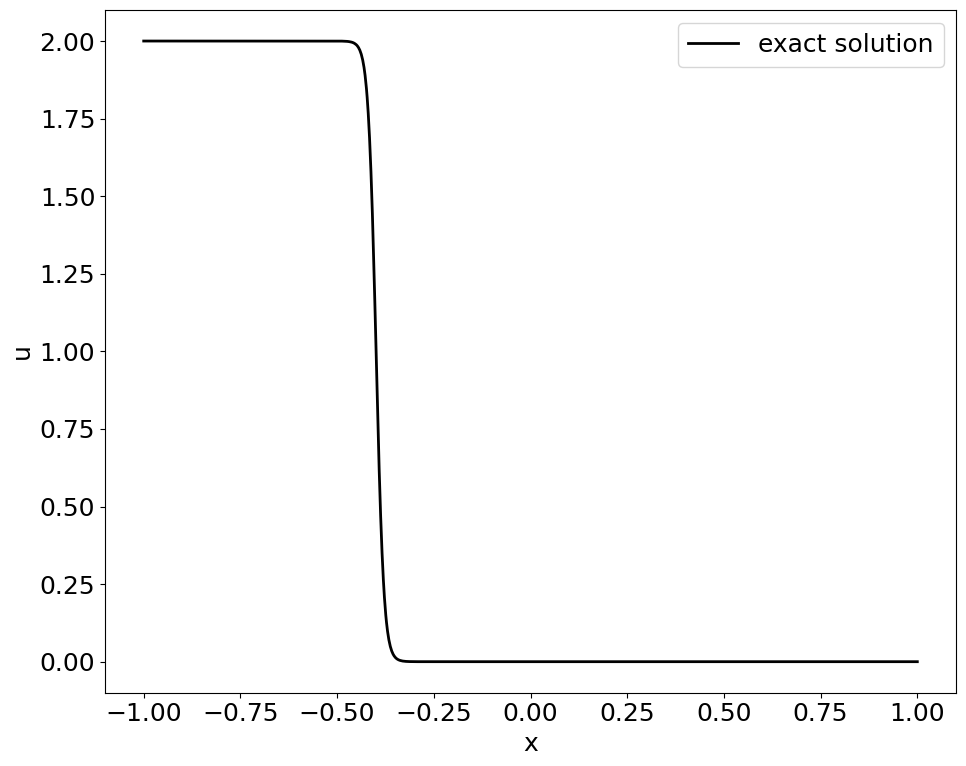}
    \caption{Burgers moving front with $\nu = 10^{-2}$ at t=0.1}
    \label{fig:burgers_mf_exact}
\end{figure}

\begin{figure}[H]
    \centering
    \begin{subfigure}[t]{0.47\textwidth}
        \includegraphics[width=\textwidth]{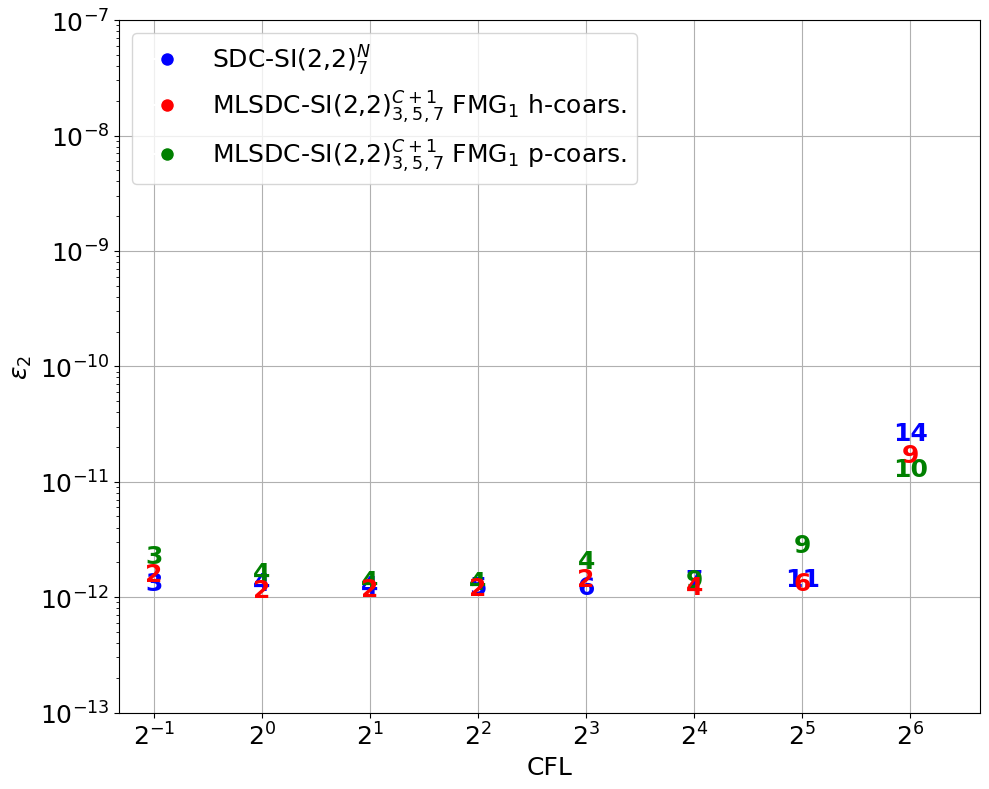}
        \caption{high spatial resolution (1)}
        \label{fig:burgers_nu1e-2_hires}
    \end{subfigure}
    \hfill
    \begin{subfigure}[t]{0.47\textwidth}
        \includegraphics[width=\textwidth]{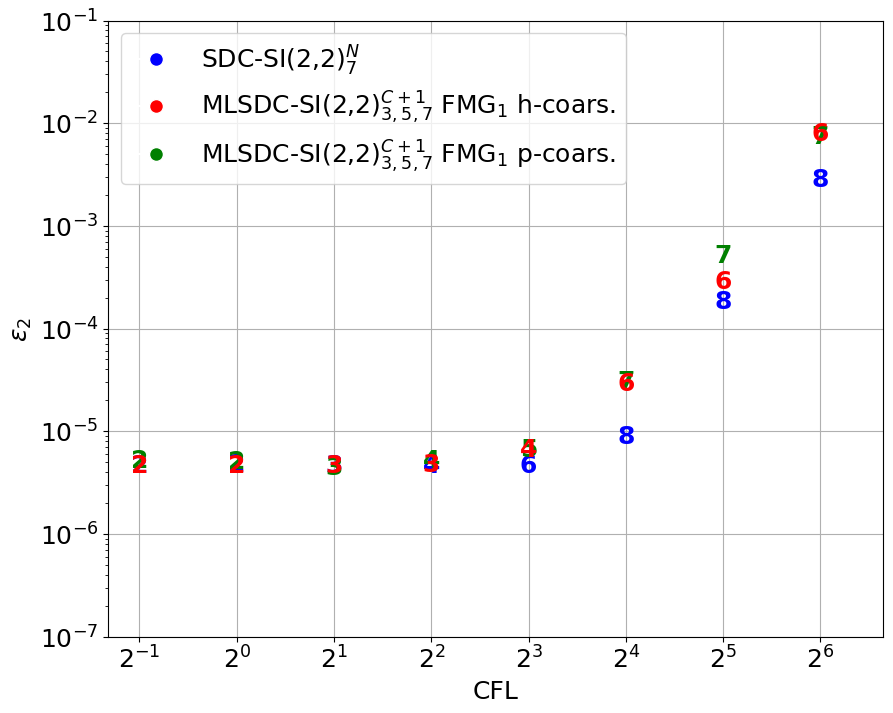}
        \caption{low spatial resolution (2)}
        \label{fig:burgers_nu1e-2_lores}
    \end{subfigure}
    \caption{$L_2$ errors over $\CFL$ numbers SDC$_7^N$ vs. MLSDC$_{3,5,7}^{C+1}$ for Burgers moving front and $\nu=10^{-2}$, top level, markers indicate the number of fine grid sweeps}
    \label{fig:burgers_mf}
\end{figure}

The Burgers equation \eqref{eq:burgers} presents a numerical challenge due to its non-linear dynamics, which can lead to the formation of shock waves and discontinuities.
Additionally, variations in the convective Jacobian matrix cause changes in diffusivity within the semi-implicit methods devised in this study.
The test case considered here relates to the moving front solution
\begin{equation}
  u(x,t) = 1 - \tanh\left(\frac{x + 0.5 - t}{2 \nu}\right)
\end{equation}
for ${x \in [-1,1]}$ and $t \in [0,0.1]$, see \cite{SE_Hesthaven2008a}.
Assuming a diffusivity of $\nu = 10^{-2}$ yields a thin advancing front (Figure \ref{fig:burgers_mf_exact}), which is expensive to resolve numerically.
We consider two cases: (1) the front is spatially well-resolved on all levels, including the coarsest; (2) the front is spatially under-resolved on all levels, including the finest. 
To capture the front and to avoid unphysical oscillations in the underresolved case, the artificial diffusivity of Persson and Peraire is used with parameters $\ds = 0.5$ and $\cs = 0.5$, see \cite{SE_Persson2006a} and \cite{TI_Stiller2024} for details.
Figure~\ref{fig:burgers_mf} displays the $L_2$ errors achieved with the different methods for the problem. The markers in the figures indicate the number of fine grid sweeps, as listed in Table~\ref{tab:burgers_iterations_cfl_hires} and \ref{tab:burgers_iterations_cfl_lores}.\\
We first consider the highly resolved test case. Table~\ref{tab:burgers_iterations_cfl_hires} presents the number of fine grid operations required to achieve convergence across various $\CFL$ numbers (compare Figure \ref{fig:burgers_nu1e-2_hires}). The fine grid $\CFL$ number is kept constant across the two different spatial resolutions. As a result, for a given $\CFL$ number, the corresponding time step/slab size differs between the test cases. 
\begin{table}[h!]
    \centering
    \begin{tabular}{cccc|cccccccc}
    \hline
    \multicolumn{4}{c|}{Method} & \multicolumn{8}{c}{$\CFL$} \\ 
      & $N_e$ & $P$ & Start  & 0.5 & 1 & 2 & 4 & 8 & 16 & 32 & 64 \\ \hline
    SDC-SI(2,2)$_7^N$  & $ 200$ & $15$ & - & 3 & 4 & 4 & 5 & 6 & 7 & 11 & 14 \\ \hline
            & $\{200,200,200\}$ & $\{15,15,15\}$ & FMG$_1$ & 2 & 2 & 2 & 2 & 2 & 4 & 6 & 9 \\   
             & $\{50,100,200\}$ & $\{15,15,15\}$ & FMG$_2$ & 2 & 2 & 2 & 2 & 2 & 3 & 6 & 9 \\
    MLSDC-SI(2,2)$_{3,5,7}^{C+1}$ & $\{50,100,200\} $ & $\{15,15,15\}$ & FMG$_1$ & 2 & 2 & 2 & 2 & 2 & 4 & 6 & 9 \\
              & $\{50,100,200\}$ & $\{15,15,15\}$ & Casc. & 2 & 2 & 2 & 2 & 3 & 5 & 7 & 10 \\
             & $\{200,200,200\} $ & $\{5,10,15\}$ & FMG$_1$ & 3 & 4 & 4 & 4 & 4 & 9 & 9 & 10 \\
    \end{tabular}
    \caption{Number of needed fine grid iterations until convergence for different $\CFL$ numbers at high spatial resolution}
    \label{tab:burgers_iterations_cfl_hires}
\end{table}
In this setting, both spatial $h$-coarsening and $p$-coarsening perform without any issues. However, $h$-coarsening proves to be more effective, as the $L_2$ error on the coarsest level is lower compared to that obtained with $p$-coarsening in this test.
The iteration savings achieved by MLSDC are significant. For certain values of the $\CFL$ number, it requires less than half as many fine grid sweeps as single-level SDC. At lower $\CFL$ numbers, MLSDC often fully converges with just one cycle, which is equivalent to only two fine grid iterations.
The results for the starting strategies described in Section~\ref{ch:mlsdc_basics} are reproduced as well.
It should also be mentioned that, in the cases with space-time coarsening, the embedded interpolation method generally outperformed the $L_2$-projection as a variable projection method. In some instances, this led to significantly lower $L_2$ errors after a fixed number of iterations, sometimes by an order of magnitude.\\
If all levels suffer from poor resolution (Figure \ref{fig:burgers_nu1e-2_lores} and Table \ref{tab:burgers_iterations_cfl_lores}), the MLSDC still performs stably, although the benefits it provides become less pronounced.
\begin{table}[h!]
    \centering
    \begin{tabular}{cccc|cccccccc}
    \hline
    \multicolumn{4}{c|}{Method} & \multicolumn{8}{c}{$\CFL$} \\
      & $N_e$ & $P$ & Start  & 0.5 & 1 & 2 & 4 & 8 & 16 & 32 & 64 \\ \hline
    SDC-SI(2,2)$_7^N$  & $ 40$ & $8$ & - & 2 & 2 & 3 & 4 & 6 & 8 & 8 & 8 \\ \hline
            & $\{40,40,40\}$ & $\{8,8,8\}$ & FMG$_1$ & 2 & 2 & 2 & 2 & 3 & 6 & 6 & 6 \\
    MLSDC-SI(2,2)$_{3,5,7}^{C+1}$ & $\{10,20,40\} $ & $\{8,8,8\}$ & FMG$_1$ & 2 & 2 & 3 & 3 & 4 & 6 & 6 & 6 \\
              & $\{40,40,40\}$ & $\{2,4,8\}$ & FMG$_1$ & 2 & 2 & 3 & 4 & 5 & 7 & 7 & 7 \\
    \end{tabular}
     \caption{Number of needed fine grid iterations until convergence for different $\CFL$ numbers at low spatial resolution}
     \label{tab:burgers_iterations_cfl_lores}
 \end{table}
The solution projected from the finer levels cannot be properly represented on the poorly resolved coarser levels. This also seems to affect the FAS RHS and the coarse-to-fine correction becomes inefficient.
For MLSDC with $h$- or $p$-coarsening, this leads to slower convergence and higher terminal $L_2$ errors compared to the single-level SDC solution for high $\CFL$ numbers.
Moreover, it was observed that an under-resolved solution on the coarse level can introduce jumps or oscillations in the numerical solution. In the context of the Cascade and FMG starting strategies, these artifacts are subsequently interpolated to the next finer level, resulting in unfavourable initial values for the MLSDC method. 
This, in turn, can delay convergence or degrade the overall efficiency of the method. As a practical workaround, the starting values are smoothed prior to interpolation, e.g. by averaging out jumps across element interfaces on the coarse level. This procedure is applied in all Burgers moving front test cases.
The issues associated with an unresolved coarse level are known to the authors. Potential remedies, such as the use of filters and smoothing of the starting solutions, are already under investigation, but beyond the scope of this study.
These limitations are less restrictive in cases with higher diffusion, as the solution is smoother in such scenarios.

\subsection{Euler and Navier-Stokes equations.}

\begin{figure}[H]
  \subcaptionbox{SDC-SI(2,2)$_{7}^{11}$ density distribution
    \label{fig:cns:wave:sdc}}
    {\includegraphics[scale=0.33]{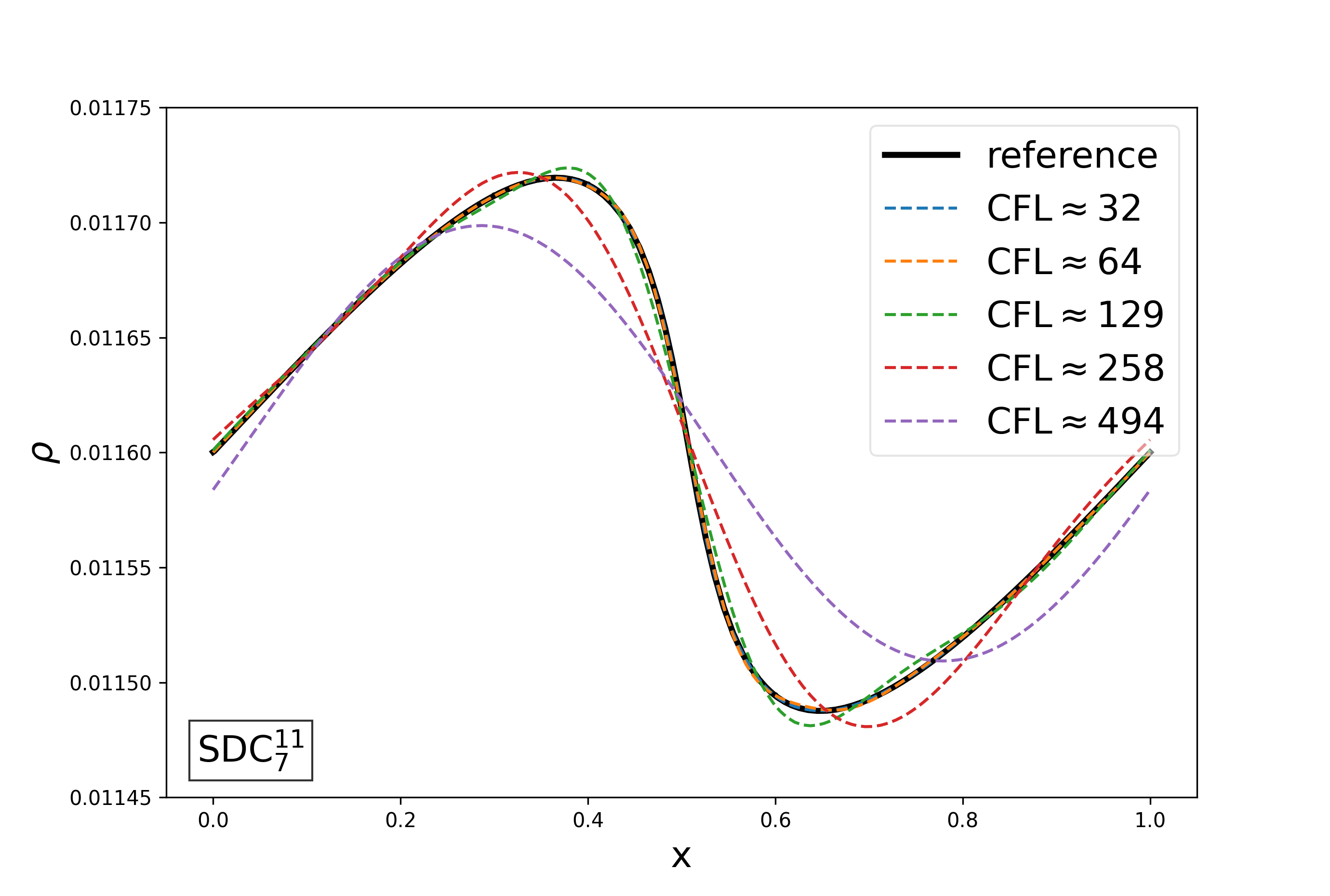}}
  \subcaptionbox{\RI{SDC-SI(2,2)$_{7}^{11}$ deviation from reference
    \label{fig:cns:wave:sdc_err}}}
    {\includegraphics[scale=0.33]{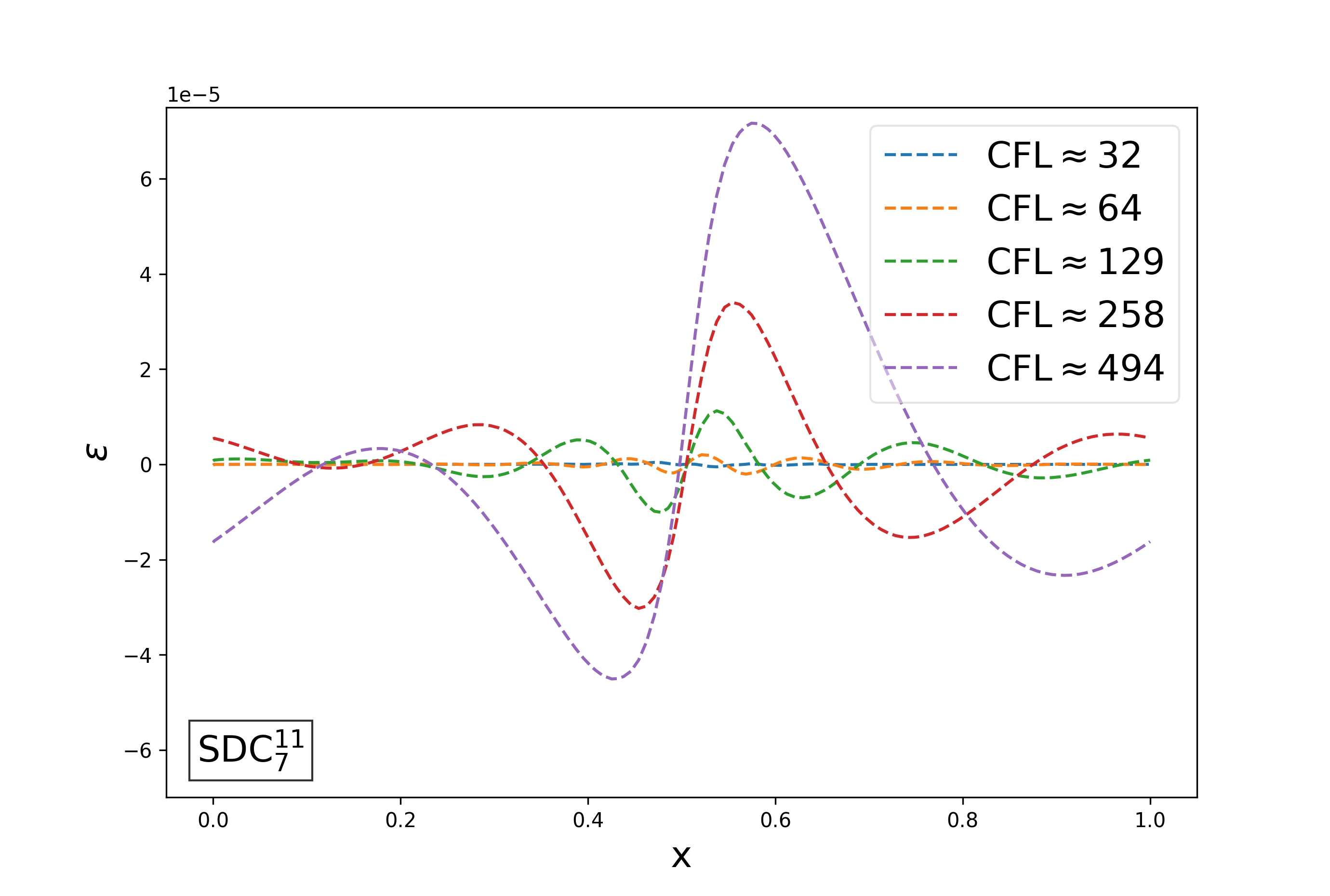}} \\[1ex]
  \subcaptionbox{MLSDC-SI(2,2)$_{5,7}^8$ density distribution
    \label{fig:cns:wave:mlsdc}}
    {\includegraphics[scale=0.33]{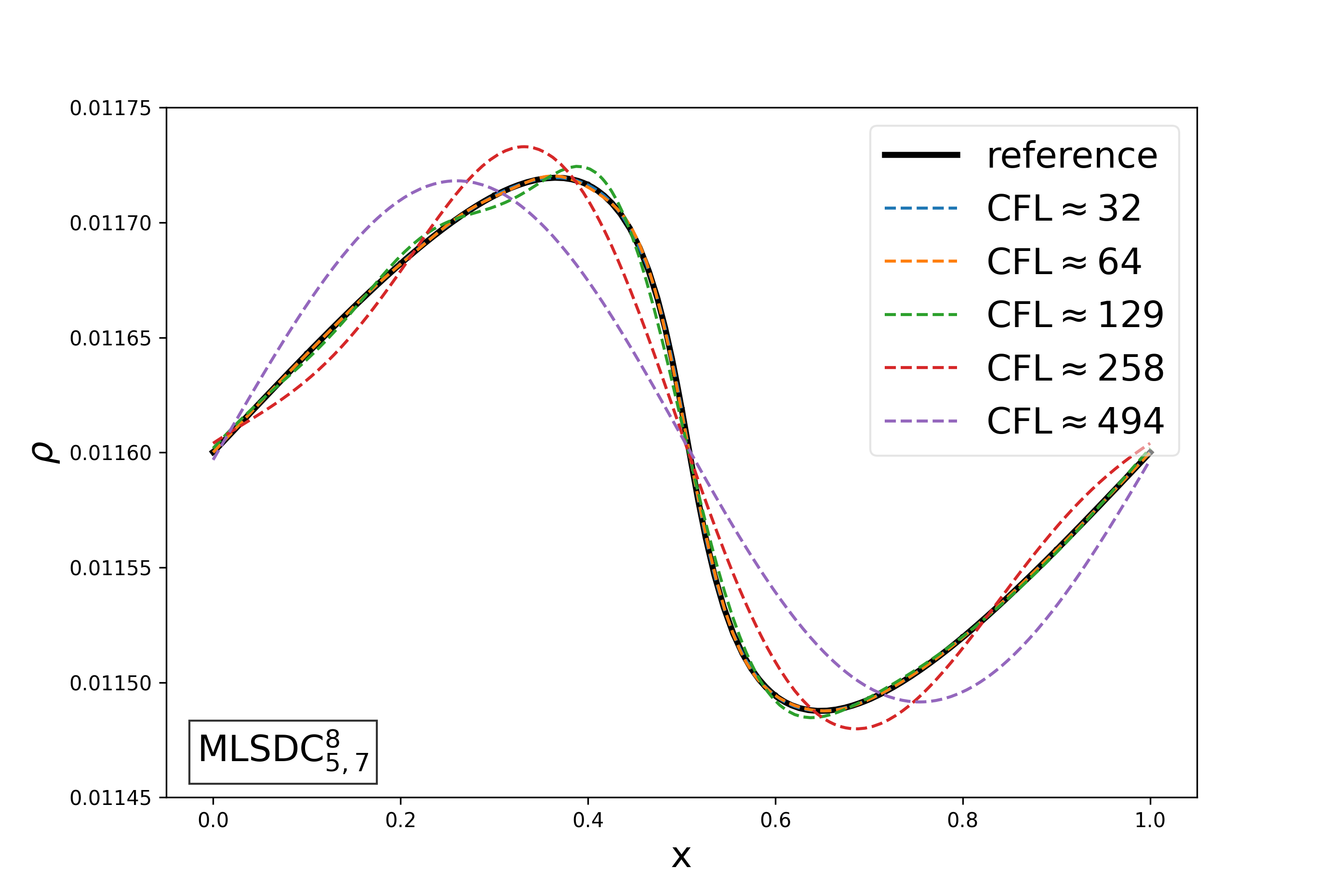}}
  \subcaptionbox{\RI{MLSDC-SI(2,2)$_{5,7}^8$ deviation from reference
    \label{fig:cns:wave:mlsdc_err}}}
    {\includegraphics[scale=0.33]{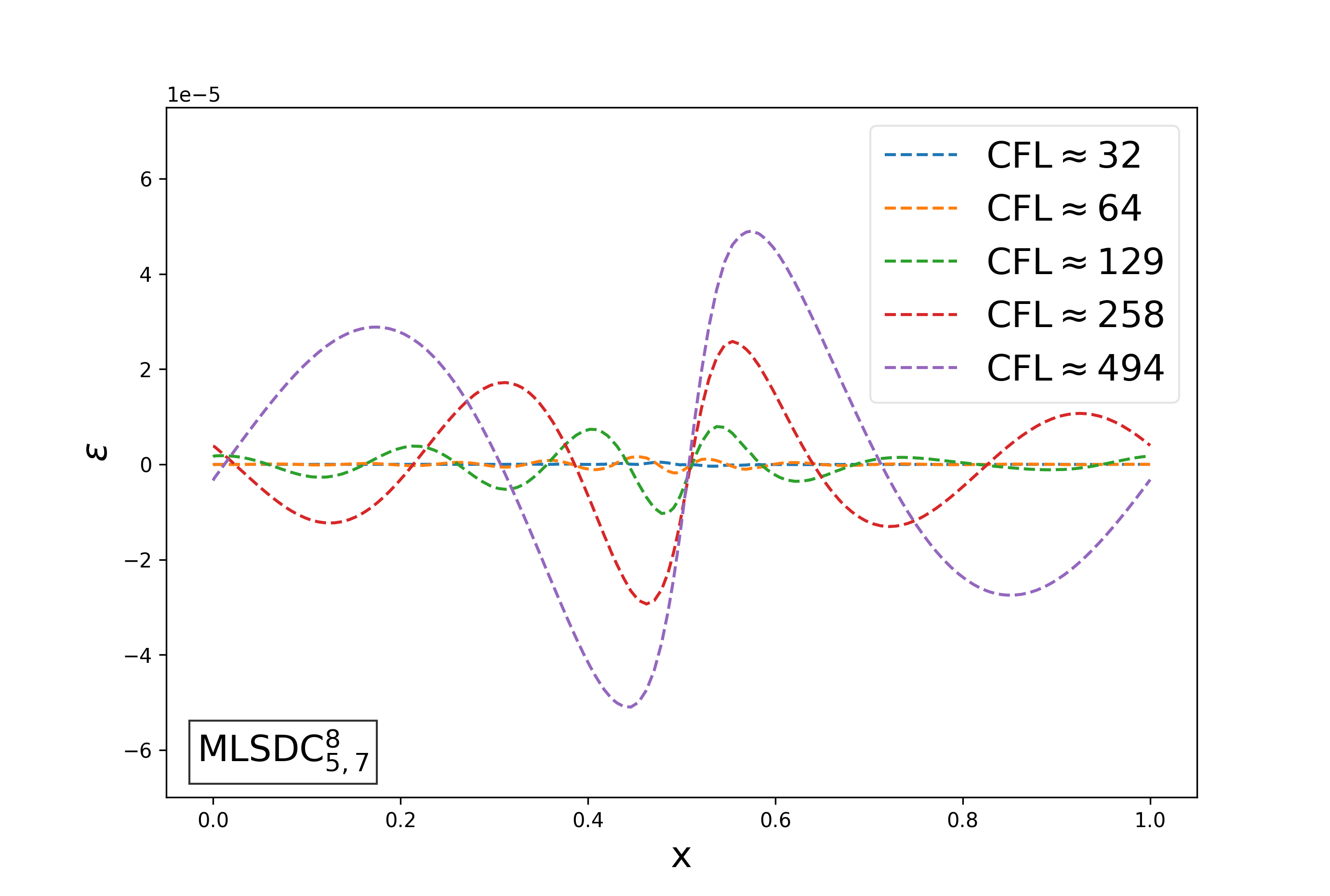}} 
  \caption{
    Acoustic wave density distributions obtained with SDC-SI(2,2)$_{7}^{11}$ vs MLSDC-SI(2,2)$_{5,7}^8$ \RI{and their deviations from a reference solution}.
    \label{fig:cns:wave_2l}}
\end{figure}

\begin{figure}[H]
    \includegraphics[scale=0.56]{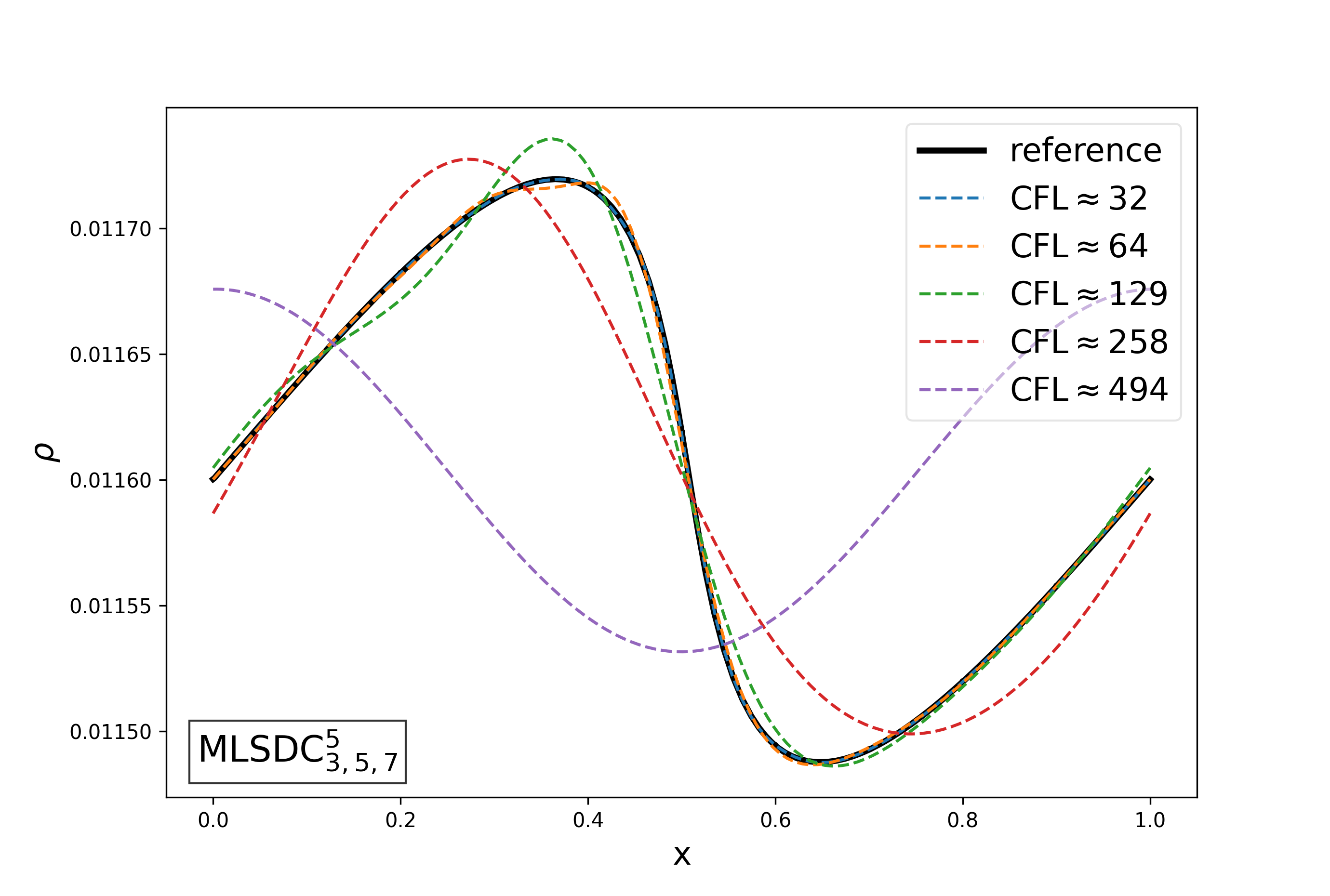} 
  \caption{Acoustic wave density distributions 
    obtained with MLSDC-SI(2,2)$_{3,5,7}^5$ FMG$_1$ variant
    \label{fig:cns:wave_3l}}
\end{figure}

Now  some compressible flow problems are investigated. 
The fluid is modeled as a perfect gas and exhibits properties corresponding to air, with a specific gas constant $R = 287.28$ and a heat capacity ratio $\gamma = 1.4$, under standard conditions and assuming SI units. 
The first problem describes an acoustic wave traveling downstream in a uniform flow characterized by a 
Mach number $\mathit{Ma} = 0.1$,
ambient pressure $p_{\infty} = 1000$, 
temperature $T_{\infty} = 300$, 
dynamic viscosity $\eta = 10^{-5}$ and  
Prandtl number $\mathit{Pr} = 0.75$. 
This wave is initialized by superimposing the periodic velocity perturbation, defined mathematically by the equation:
\begin{equation}
  v'(x) = \hat{v} \sin(2 \pi x),
\end{equation}
with an amplitude of 
${\hat v = 10^{-2} a_{\infty}}$, where 
${a_{\infty} \approx 347.36}$ is the speed of sound in the undisturbed fluid.
Assuming an isotropic initial state, the inital temperature is given by
\begin{equation}
  T_0(x) = \frac{1}{\gamma R}\left( a_{\infty} +\frac{\gamma-1}{2}v' \right)^2
  \,,
\end{equation}
the pressure by
\begin{equation}
  p_0(x) = p_{\infty} \left(\frac{T}{T_{\infty}}\right)^{\gamma/(\gamma-1)}
  \,,
\end{equation}
and the velocity by
\begin{equation}
  v_0(x) = \mathit{Ma} \, a_{\infty} + v'(x)
  \,.
\end{equation}
It is to be noted, that the amplitude of the wave is sufficiently large to trigger non-linear steepening effects in the wave dynamics. 
Conversely, the impact of viscous damping is considered negligible owing to the high Reynolds number $\mathit{Re} = \rho_{\infty} a_{\infty}/\ \eta \approx 4 \times 10^5$.
The progression of the wave is studied within a periodic domain $ \Omega = [0,1]$, from $t = 0$ to $t_{\mathrm{end}} = 2.619 \times 10^{-2}$. 
This time frame approximately covers 10 \RI{wave passages}.
The discretization parameters employed in the study are detailed in Table \ref{tab:num_experiments_cns_wave}.
\begin{table}[H]
\centering
\begin{tabular}{cccccc}
\toprule
 & \multicolumn{2}{c}{Space} & \multicolumn{2}{c}{Time} \\
\cmidrule(r){2-3} \cmidrule(r){4-5}
 & $N_{e,l}$ & $P_l$ & $M_l$ \\
\midrule
Level 1 & 6 & 8 & 3 \\
Level 2 & 12 & 8 & 5 \\
Level 3 & 12 & 15 & 7 \\
\bottomrule
\end{tabular}
\caption{Discretization parameters for the acoustic wave test}
\label{tab:num_experiments_cns_wave}
\end{table}
Since no exact solution is available, the TVD-RK3 method of Shu and Osher \cite{TI_Shu1988b} with ${\Delta t \approx 10^{-7}}$ or ${\CFL \approx 0.019}$ serves as a reference.
The first test corresponds to a 2-level MLSDC method (level 2 and level 3 from Table \ref{tab:num_experiments_cns_wave}) with a Cascade start vs. a single-level SDC method.
The results can be seen in Figure \ref{fig:cns:wave_2l}. The MLSDC method can stably solve the CNS acoustic wave problem for $\CFL$ numbers up to 494.
After 7 cycles (8 fine grid sweeps), it already shows better agreement with the reference than single-level SDC after 11 sweeps, especially for high $\CFL$ numbers (494 and 258).
The second test considers a three-level discretization using the FMG$_1$ starting variant, see Figure~\ref{fig:cns:wave_3l}.
Compared to the two-level method, this setup achieves even greater iteration savings for lower $\CFL$ numbers.
However, the method encounters difficulties in accurately solving the problem for greater $\CFL$ numbers, particularly at $\CFL = 494$. 
It may be suspected that this issue is related to an insufficient temporal resolution on the coarse level, which worsens when the $\CFL$ number increases.
Increasing the number of cycles does not improve the solution; instead, it introduces slight instabilities.
The three-level method still successfully solves the problem for $\CFL$ numbers up to 258, which remains a notable result.
In summary, the three-level MLSDC offers clear advantages at small $\CFL$ numbers. At larger $\CFL$ values, however, the limitations due to poor resolution on the coarsest grid become apparent.

\begin{figure}[H]
  \subcaptionbox{SDC-SI(2,2)$_7^N$
    \label{fig:cns:st_12_sdc}}
    {\includegraphics[scale=0.3]{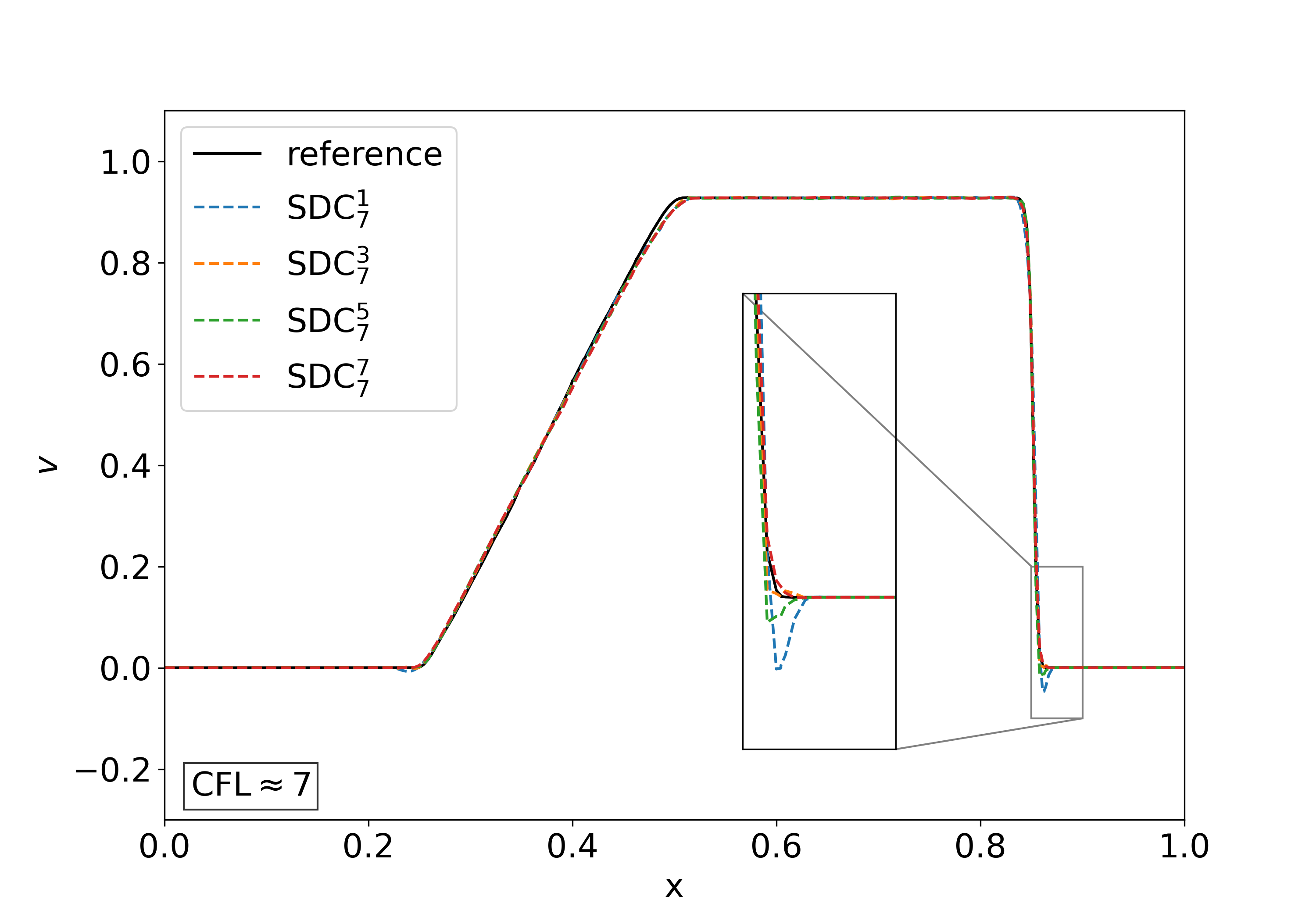}}
  \subcaptionbox{MLSDC-SI(2,2)$_{5,7}^C$
    \label{fig:cns:st_12_mlsdc}}
    {\includegraphics[scale=0.3]{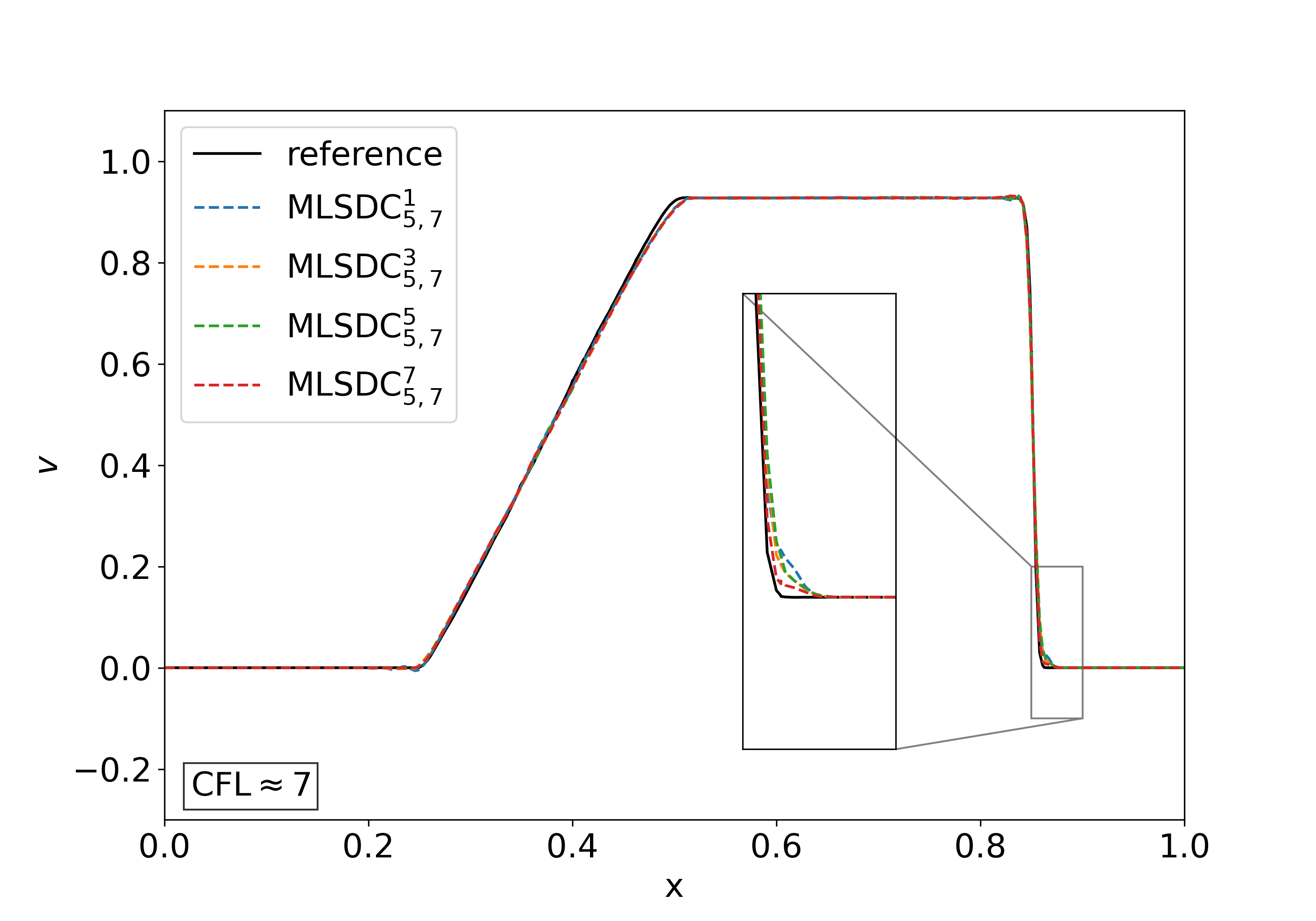}}
  \caption{%
    Shock Tube velocity distributions 
    obtained with SDC-SI(2,2)$_{7}^{N}$ vs MLSDC-SI(2,2)$_{5,7}^C$ at $\CFL \approx 7$
    \label{fig:cns:st_7}}
\end{figure}

The second example is Sod's shock tube problem for the Euler equations with initial conditions \cite{BM_Sod1978a}
\begin{subequations}
  \begin{alignat}{4}
    &\rho = 1,     &\quad& p = 1,   &\quad& v = 0 && \qquad\text{if} \quad x <   0.5 \\
    &\rho = 0.125, &\quad& p = 0.1, &\quad& v = 0 && \qquad\text{if} \quad x \ge 0.5
  \end{alignat}
\end{subequations}
to solve on ${ \Omega = [0,1]}$ until ${t = 0.2}$.
For the numerical studies, the domain is divided into 80 elements of degree 5 on the finest level.
The discretization parameters for the test can be taken from Table \ref{tab:num_experiments_cns_st}, the number of coarse sweeps $N_c = 2$.
\begin{table}[H]
\centering
\begin{tabular}{cccccc}
\toprule
 & \multicolumn{2}{c}{Space} & \multicolumn{2}{c}{Time} \\
\cmidrule(r){2-3} \cmidrule(r){4-5}
 & $N_{e,l}$ & $P_l$ & $n_{t,l}$ & $M_l$ \\
\midrule
Level 1 & 80 & 3 & 32 & 5 \\
Level 2 & 80 & 5 & 32 & 7 \\
\bottomrule
\end{tabular}
\caption{Discretization parameters for the shock-tube test}
\label{tab:num_experiments_cns_st}
\end{table}
The artificial diffusivity of Persson and Peraire \cite{SE_Persson2006a} is used with $\ds = 2$ and $\cs = 0.4$ to capture the evolving shock and the contact discontinuity.
Figure \ref{fig:cns:st_7} shows the velocity $v$ obtained with SDC-SI(2,2)$_7^{N}$ and MLSDC-SI(2,2)$_{5,7}^{C}$ using $32$ time steps ($\CFL \approx 7$) after different numbers of iterations. The MLSDC used the Cascade starting strategy.
The TVD-RK3 with $10^5$ steps or ${\Delta t = 2\times 10^{-6}}$ serves as a reference.
The shock tube problem could only be solved stably with spatial coarsening, using $P_l=5$ on all levels led to instability, likely due to an excessively large coarse grid $\CFL$ number. 
With 32 time steps, the agreement with the reference solution is already very good.
Although both SDC and MLSDC appear to resolve the wave front around $x \approx 0.85$ equally well, there are overshooting oscillations in the lower part of the domain in the SDC solution, which can not be observed in the MLSDC method.
In addition, the MLSDC method converges toward the reference solution more rapidly.

\subsection{Runtime analysis}

\RIV{
This paper is focused on stability, convergence and accuracy of the SI-MLSDC method.
In our comparative analysis, we specifically examine the iterations executed on the fine grid for the single-level SDC and the MLSDC.
At this point, a small runtime analysis shall be provided.} \\
The computational model utilized in this small analysis is intentionally simplified.
These simplifications include the absence of factors such as communication overhead and the exclusion of additional runtime overhead from the transfer operators in MLSDC.
In addition, all implicit equations are currently being iterated to convergence at all levels (e.g. within the diffusion solver).
However, this streamlined model allows for a clearer examination of the core computational elements without the confounding effects of more complex system interactions.
For estimating the runtime of single-level SDC implementations, we suggest the following:
\begin{equation}
    t_{\text{Run,SDC}} \approx c_{R} \times (P+1) \times N_e \times M \times n_t \times N \,,
\end{equation}
where $c_R$ is a coefficient representing the runtime associated with a single operation within the SDC framework, encompassing activities such as the solution of implicit parts and spatial discretization. 
Additionally, $N$ denotes the number of sweeps performed.
This formulation provides a foundational estimate, albeit approximate.
The computational efficiency of the MLSDC method can be captured through the following modeled runtime estimate:
\begin{equation}
\begin{split}
    t_{\text{Run,MLSDC}} & \approx \sum_{l=2}^L c_{R} \times (P_l +1 ) \times N_{e,l} \times M_l \times n_{t,l}  \times (C+1)  \\
    & + c_{R} \times (P_1+1) \times N_{e,1} \times M_1 \times n_{t,1}  \times N_c \times C \,.
\end{split}
\end{equation}
In scenarios where only $p$-refinement in space and time is applied, and the standalone SDC utilizes the fine grid discretization, a two-level MLSDC would offer a time-saving advantage under the condition:
\begin{equation}
\begin{split}
    (P_2+1) \times M_2 \times N & > (P_2 +1 ) \times M_2 \times (C+1)  \\
    & + (P_1+1) \times M_1  \times N_c \times C \,.
\end{split}
\end{equation}
For a simplified assumption, $N_c = 2$, $P_1$ and $P_2$ are high polynomial orders and $P_1 = P_2/2$ and $M_1 = M_2/2$, the MLSDC approach becomes beneficial if it converges after fewer than $C < 2/3 \times (N-1)$ cycles.\\
\RIV{
This model should be reinforced with a simple numerical example.
For the convection–diffusion problem described in Section~\ref{sec:num_convdiff}, the runtime of selected MLSDC methods is compared against single-level SDC as well as widely used competitive time integration methods, namely TVD-RK3 \cite{TI_Shu1988b} and the third order IMEX-RK ARS(4,4,3) \cite{TI_Ascher1997a}. 
As a representative test case, we consider the configuration given in Table~\ref{tab:num_experiments_convdiff_1} and a diffusivity of $\nu = 10^{-3}$.
For the case without diffusivity, explicit schemes such as TVD-RK3 clearly outperform semi-implicit MLSDC or SDC approaches; see, for example, Section~6.4 in~\cite{TI_Stiller2024}.
Once diffusivity is introduced, however, semi-implicit SDC and MLSDC schemes demonstrate distinct advantages.
Both two-level and three-level MLSDC configurations are examined for this problem.
All simulations are performed using the HiSPEET framework, compiled with the Intel Fortran Compiler 2021.10.0 (build date 20230609). Runs are executed in single-threaded mode (no parallelization) on an Intel Xeon Platinum 8470 CPU (52 cores) operating at 2.00\,GHz.
Each runtime experiment is repeated ten times, and the reported runtime corresponds to the average across these runs.\\
\begin{table}[h!]
\centering
\renewcommand{\arraystretch}{1.2}
\setlength{\tabcolsep}{4pt}
\small
\begin{tabular}{lcccccc}
\toprule
& \multicolumn{2}{c}{Level 3} & \multicolumn{2}{c}{Level 2} & \multicolumn{2}{c}{Level 1} \\
\cmidrule(lr){2-3}\cmidrule(lr){4-5}\cmidrule(lr){6-7}
Case & Rel. cost [\%] & Ops & Rel. cost [\%] & Ops & Rel. cost [\%] & Ops \\
\midrule
1D MLSDC$_{5,7}$ & 100 & $C$ & 36 & $2 \cdot C$ & -- & --  \\
1D MLSDC$_{3,5,7}$ & 100 & $C$ & 36 & $N_c \cdot C$ & 11 & $N_c \cdot C$ \\
3D MLSDC$_{5,7}$  & 100 & $C$ & 9 & $2 \cdot C$ & -- & -- \\
3D MLSDC$_{3,5,7}$ & 100 & $C$ & 9 & $N_c \cdot C$ & 0.6 & $N_c \cdot C$ \\
\bottomrule
\end{tabular}
\caption{Operation count and comp. cost per level relative to the finest discretization}
\label{tab:level-costs}
\end{table}
Including a post-smoothing step implies that the finest level of the MLSDC method would perform $C+1$ fine-grid iterations. 
However, for the runtime analysis some optional MLSDC parameters were optimized for better performance.
The numerical experiments revealed that omitting the post-smoothing step in the MLSDC method improved the overall time-to-solution for the convection–diffusion problem, as it significantly reduced the computational cost while only marginally increasing the error.
It shall be noted that this observation may not necessarily generalize to other test cases.\\
Following the estimates in Table~\ref{tab:level-costs}, for the numerical example with our discretization parameters, the three-level 1D MLSDC with $N_c=2$ requires $1.94C$ fine-grid operations.
In addition, there are overheads for transfer operators and due to the starting strategy. The FMG$_2$ start for this method requires addionally $1.27$ fine grid operations.
A MLSDC$_{3,5,7}^C$ method without post-smoothing and using an FMG$_2$ starting strategy is expected to be approximately as computationally expensive as an SDC$_7^9$ method with four cycles, excluding the transfer operator overhead. Consequently, if the same error can be achieved with an MLSDC$_{3,5,7}^3$ scheme, a theoretical runtime reduction is obtained. 
This comparison is illustrated in Figure~\ref{fig:convdiff_runtime1}.
\begin{figure}[H]
  \subcaptionbox{SDC and MLSDC with lesser number of sweeps/cycles
    \label{fig:convdiff_runtime1}}
    {\includegraphics[scale=0.48]{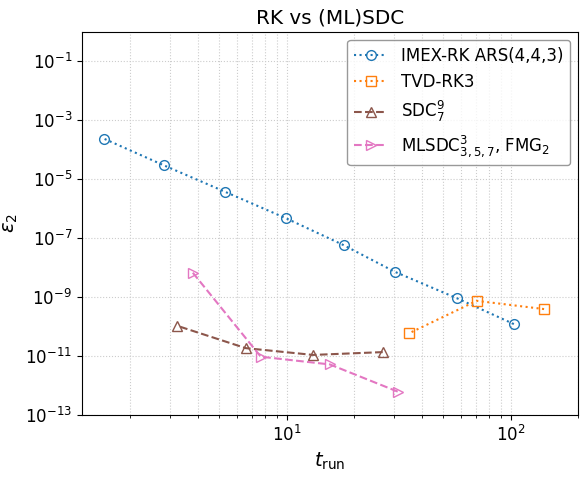}}
  \subcaptionbox{SDC and MLSDC with higher number of sweeps/cycles
    \label{fig:convdiff_runtime2}}
    {\includegraphics[scale=0.48]{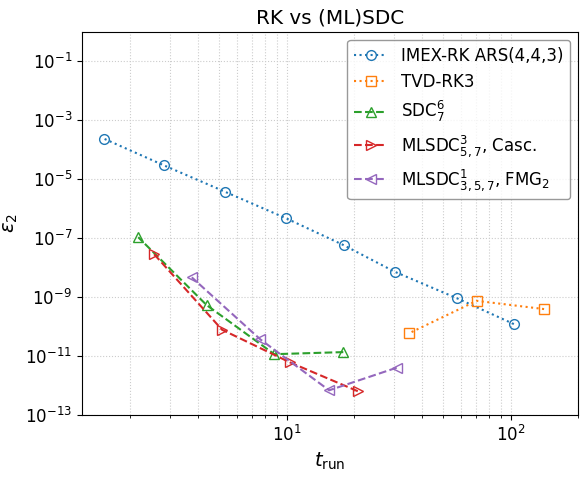}}
  \caption{%
    Error-runtime comparison between Runge-Kutta, MLSDC and SDC methods using either the preconditioned conjugate gradient method (PCG) or the fast direct solver (FDS)
    \label{fig:convdiff_runtime}}
\end{figure}
The estimation agrees well with measured runtimes, except when runtimes become very small, i.e. at large $\CFL$ numbers, where the overhead of the transfer operators becomes more pronounced. 
Both SDC and MLSDC demonstrate competitive performance compared to TVD-RK3 and IMEX-RK(4,4,3). The TVD-RK3 method quickly exhibits instability when the time step size exceeds a certain threshold.
This behavior was already observed in \cite{TI_Stiller2024} for the case with diffusion.
The same test was repeated with fewer sweeps or cycles to not give the impression that an unreasonable high number of iterations had been performed using the SDC method unnecessarily.
The results are shown in Figure ~\ref{fig:convdiff_runtime2}.
The investigation also includes a two-level MLSDC scheme, which requires $1.72C$ fine-grid operations in addition to $0.36$ for the Cascade start. 
The results shown in Figure~\ref{fig:convdiff_runtime1} can be reproduced reasonably well, indicating that the runtime estimates for MLSDC are consistent. 
For global refinement in one dimension, only minor runtime savings are observed compared to the single-level SDC method in this test case. However, achieving such savings was not the primary objective of this work. 
The advantages of the MLSDC approach are twofold: 
\begin{itemize}
    \item It provides a natural pathway for local adaptivity, allowing finer levels to be employed only where needed rather than globally. This is the focus of ongoing research and promises substantially greater runtime savings. 
    \item The runtime estimation in one dimension shows good agreement with the numerical experiment. Based on the estimates in Table~\ref{tab:level-costs}, even larger savings are expected in three dimensions, since the relative cost of the coarse levels decreases compared to the fine level. This aspect is also the subject of current investigations. 
\end{itemize}
Other factors naturally influence performance, such as the number of iterations required by the implicit solvers. 
In principle, the ML code could be further optimized for better runtime, for example by incompletely solving implicit equations on intermediate levels, exploring alternative starting strategies, optimizing implementation details, or applying $h/p$-coarsening in space. 
However, such efforts require a delicate balance between error and runtime and are of limited relevance, since we expect more significant improvements from adaptivity and three-dimensional simulations than from fine-tuning the one-dimensional ML algorithm. 
From this perspective, the present paper serves as a pilot study, aiming to demonstrate a proof of concept for the proposed method. 
It is therefore acceptable that globally refined one-dimensional MLSDC does not yet yield substantial runtime savings. 
We have shown that the number of iterations can be reduced significantly, and considerable runtime savings are expected in future work. 
}

\subsection{Convergence analysis}
\RIV{
The convection–diffusion problem described in Section~\ref{sec:num_convdiff} is revisited here to perform a brief convergence analysis.
For this convergence study, the spatial and temporal orders are aligned on the finest level, requiring slight adjustments to the discretization parameters from Section \ref{sec:num_convdiff}: on the finest grid, $M=7$ collocation nodes with Radau Right points yield a temporal order of $\mathcal{O}(t) = 2M-1=13$, corresponding to a spatial polynomial degree of $P=12$ (i.e., order $\mathcal{O}(x)= P+1=13$).
The spatial order remains identical across all levels and for each mesh, while the number of temporal collocation nodes varies as $M_l = \{3, 5, 7\}$. The number of elements is doubled on each higher level. Consequently, the refinement corresponds to an $h$-refinement in space and a $p$-refinement in time.
The errors are slightly higher than in Section~\ref{sec:num_convdiff} since the polynomial degree $P_l$ is lower on all levels in this study.
The tested method is a MLSDC$_{3,5,7}^{10}$ method, with sufficient cycles to converge to the minimum error.
\begin{table}[H]
\centering
\renewcommand{\arraystretch}{1.2}
\setlength{\tabcolsep}{6pt}
\small
\begin{tabular}{lcccc}
\toprule
    ML Mesh $\{N_{e,l}\}$ & $\CFL=32$ & $\CFL=64$ & $\CFL=128$ & $\CFL=256$ \\
    \midrule
    $\M \Omega_{(1)}\,\,\{5,10,20\}$    & $3.519\times10^{-9}$  & $5.555\times10^{-7}$  & $3.905\times10^{-3}$  &  $2.443\times10^{-1}$ \\
    $\M \Omega_{(2)}\,\,\{10,20,40\}$   & $1.385\times10^{-11}$ & $1.693\times10^{-10}$ & $5.556\times10^{-7}$  & $3.905\times10^{-3}$ \\
    $\M \Omega_{(3)}\,\,\{20,40,80\}$   & $3.556\times10^{-11}$ & $4.864\times10^{-11}$ & $2.064\times10^{-10}$ & $5.556\times10^{-7}$ \\
    $\M \Omega_{(4)}\,\, \{40,80,160\}$  & $9.507\times10^{-11}$ & $1.371\times10^{-10}$ & $2.425\times10^{-10}$ & $3.348\times10^{-10}$ \\
\bottomrule
\end{tabular}
\caption{$L_2$ errors $\varepsilon_2$ for each spatial discretization and CFL. Different spatial meshes $\M \Omega_{(k)} = \{\Omega_1, \Omega_2,\Omega_3\}_{(k)}$ are shown as the number of elements on all levels}
\label{tab:convergence_errors-all}
\end{table}
The corresponding results are presented in Figure~\ref{fig:convdiff_nu1e-3_convergence} and Table~\ref{tab:convergence_errors-all}.
The study proceeds as follows: for a fixed $\CFL$ number, the number of elements is successively doubled on all levels to increase the resolution, with the time step adjusted accordingly. This process is then repeated for different $\CFL$ numbers, i.e. time steps.
Results can be seen in Fugure~\ref{fig:convdiff_nu1e-3_convergence}.
\begin{figure}[H]
    \centering
        \includegraphics[width=0.8\textwidth]{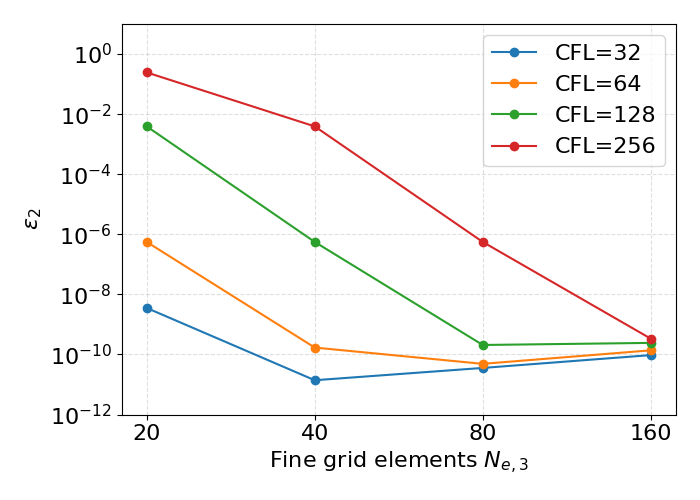}
        \label{cd_nu1e-3_xxa}
    \caption{$L_2$ errors over different spatial discretizations for various $\CFL$ numbers, MLSDC$_{3,5,7}^{10}$ for convection-diffusion and $\nu=10^{-3}$, top level}
    \label{fig:convdiff_nu1e-3_convergence}
\end{figure}
The observed convergence orders are determined (i) for a fixed $\CFL$ number between two successive spatial discretizations, and (ii) for the same spatial discretization between two different $\CFL$ numbers. 
The observed order reflects the combined accuracy in space and time (i.e., the minimum of the spatial and temporal orders), which were aligned through the chosen discretization parameters.
As for the other tests, the $\CFL$ number is determined from the finest grid. 
\begin{table}[H]
\centering
\renewcommand{\arraystretch}{1.2}
\setlength{\tabcolsep}{6pt}
\small
\begin{subtable}{0.48\textwidth}
\centering
\begin{tabular}{l|cccc}
\toprule
 & \multicolumn{4}{c}{$\CFL$} \\
\cmidrule(l){2-5}
ML mesh & $256$ & $128$ & $64$ & $32$ \\
\midrule
$\M\Omega_{(2)}$ & 6.0 & 12.8 & 11.7 & 8.0 \\
$\M\Omega_{(3)}$ & 12.8 & 11.4 & - & - \\
$\M\Omega_{(4)}$ & 10.7 & - & - & - \\
\bottomrule
\end{tabular}
\caption{Measured order at fixed $\CFL$ number}
\label{tab:convergence_a}
\end{subtable}
\hfill
\begin{subtable}{0.48\textwidth}
\centering
\begin{tabular}{l|ccc}
\toprule
 & \multicolumn{3}{c}{ML mesh } \\
\cmidrule(l){2-4}
$\CFL$ & $\M\Omega_{(1)}$ & $\M\Omega_{(2)}$ & $\M\Omega_{(3)}$ \\
\midrule
128 & 6.0 & 12.8 & 11.4 \\
64  & 12.8 & 11.7 & - \\
32  & 7.3 & - & - \\
\bottomrule
\end{tabular}
\caption{Measured order at fixed multilevel mesh $\M\Omega_{(k)}$}
\label{tab:convergence_b}
\end{subtable}
\caption{Calculated convergence orders from observed errors: (a) fixed $\CFL$ number, (b) fixed multilevel mesh $\M\Omega_{(k)}$}
\label{tab:orders-space-time}
\end{table}
It is difficult to determine the spatial and temporal orders of accuracy separately, as both components contribute to the overall error. Consequently, the observed convergence rates correspond to the minimum of the spatial and temporal orders.
The observed order was based on error reduction with a refinement ratio of~2, achieved by doubling the element count (fixed $\CFL$ number) or halving the time step, i.e. halving the $\CFL$ number (fixed spatial mesh).
For spatial refinement (fixed $\CFL$, Table~\ref{tab:convergence_a}):
\begin{equation}
\mathcal{O}_{\mathrm{\M \Omega}} 
= \log_2\!\left(\frac{\varepsilon_{\M \Omega_{(k-1)}}}{\varepsilon_{\M\Omega_{(k)}}}\right).
\end{equation}
For temporal refinement (fixed spatial mesh, Table~\ref{tab:convergence_b}):
\begin{equation}
\mathcal{O}_{\CFL} 
= \log_2\!\left(\frac{\varepsilon_{\CFL_b}}{\varepsilon_{\CFL_a}}\right).
\end{equation}
Convergence plots are shown for $\CFL$ numbers up to 256, as the MLSDC method achieves very high order in both space and time (as seen in Figure~\ref{fig:convdiff_nu1e-3}, convergence to the spatial error limit is already observed for $\CFL = 64$ at moderate spatial resolution).
The observed convergence orders are generally close to the expected value, around $13$, in space and time. As the solution approaches the asymptotic error limit, either at low $\CFL$ numbers or high spatial resolution, the measured order gradually decreases.
The study demonstrates that the observed convergence rates in both space and time closely approach the corresponding theoretical orders.
}

%% file: conclusion.tex
\section{Conclusions}
\label{ch:conclusion}

This paper introduces a new class of robust semi-implicit multilevel spectral deferred correction (SI-MLSDC) methods, characterized by excellent stability properties. 
The core foundation of these methods is the use of SI time integrators, which employ a semi-implicit splitting of the convection term inspired by the Lax-Wendroff method while treating diffusion and source terms implicitly, following the approach in \cite{TI_Stiller2024}.
The stability and accuracy of the proposed methods are thoroughly analyzed in Section~\ref{ch:stability} using a Dahlquist-type problem. 
The SI(2,2) time integrator demonstrated the best stability and accuracy properties. In contrast, the IMEX Euler method exhibited the poorest performance and was therefore not considered further in this study.
Moreover, the non-equivalence between the incremental and non-incremental MLSDC formulations is shown, and the former exhibits favourable stability properties. 
Among the aspects considered in this work is a detailed presentation and discussion of transfer operators, particularly the choice between $L^2$-projection and embedded interpolation as a projection method.
To further assess their numerical performance, the SI-MLSDC method was applied to high-order discontinuous Galerkin formulations of one-dimensional conservation laws, demonstrating high-order accuracy in space and time.
In the numerical studies, the convergence properties of the proposed SI-MLSDC methods were examined, along with key influencing factors such as different starting strategies (constant, cascadic, predictor, FMG) and the number of coarse sweeps, $N_c$. 
Tests using the convection-diffusion equation highlight the excellent stability and accuracy of SI-MLSDC, remaining stable for $\CFL$ numbers up to at least 64.
It is demonstrated, that the proposed SI-MLSDC method reduces the number of iterations required on the fine grid, with these savings becoming even more pronounced for the Burgers moving front problem. \RII{The order ideally increases by one per sweep in single-level SDC, but by significantly more than one per cycle in the present MLSDC approach.}

For instance, for certain $\CFL$ numbers, MLSDC-SI(2,2)$_{3,5,7}^{C+1}$ FMG required less than half the number of fine-grid sweeps compared to SDC-SI(2,2)$_7^N$. 
Similar iteration savings on the fine grid and faster convergence were also observed for other non-linear problems, including the Euler and Navier-Stokes equations.
In some cases, stabilization techniques were necessary near discontinuities. This included artificial diffusion methods as well as a combined $p$-coarsening strategy in both space and time, rather than $p$-coarsening in time alone, to reduce the $\CFL$ number on the coarse grid. 
Examples such as the Shock Tube test case illustrate this need.
\RIII{Future work should also focus on extending the MLSDC framework toward shock-capturing applications by incorporating monotonicity-preserving detection and correction mechanisms, similar to those employed in MOOD or TVD–IMEX schemes, to enhance robustness and stability at large $\CFL$ numbers.}
Issues related to an under-resolved problem on the coarsest grid were identified, and potential solution strategies were discussed.
These observations may align with the insights regarding an overly permissive coarsening strategy. A solution to this limitation is under investigation.
While the current implementation does not always yield direct runtime savings, partly due to the code not being optimized for performance, future improvements are expected with the introduction of local adaptivity in combination with reliable error estimators and the application to three-dimensional problems. 
At this stage, MLSDC provides the necessary infrastructure for further advancements and already displays considerable fine grid operation savings while remaining stable for a wide variety of problems.
The next steps involve local refinement in space-time, following the ideas of \citet{MG_Emmett2019} but with arbitrary high orders, more refined error estimators, improved stability, and greater flexibility.

%% file: appendix.tex



\section{Equivalence of incremental and non-incremental two-level Spectral Deferred Correction}
\label{sec:app:mlsdc:equivalence}

Consider ${M_l = 2}$, denote 
${\M f^{k-1}_{l,i} = \M f(\M u^{k-1}_{l,i},t_{l,i})}$.
Non-incremental MLSDC becomes 
\begin{subequations}
 	\begin{alignat}{2}
    & \M u^{k}_{l,1} 
    &&  = \M u^{k}_{l,0} 
        + \Delta t \sum_{i=1}^{M_l} w^{\mathsc{0n}}_{l,i1} \M f^{k-1}_{l,i}
        + \M g^{\mathsc{0n}}_{l,1}
        + \M H^k_{l,1} - \M H^{k-1}_{l,1}
  \\
    & \M u^{k}_{l,2} 
    &&  = \M u^{k}_{l,0} 
        + \Delta t \sum_{i=1}^{M_l} w^{\mathsc{0n}}_{l,i2} \M f^{k-1}_{l,i}
        + \M g^{\mathsc{0n}}_{l,2}
        + \bigl( \M H^k_{l,1} + \M H^k_{l,2} \bigr)
        - \bigl( \M H^{k-1}_{l,1} + \M H^{k-1}_{l,2} \bigr)
  \end{alignat}
\end{subequations}
and incremental MLSDC
\begin{subequations}
 	\begin{alignat}{2}
    & \M u^{k}_{l,1} 
    && = \M u^{k}_{l,0} 
       + \Delta t \sum_{i=1}^{M_l} w^{\mathsc{nn}}_{l,i1} \M f^{k-1}_{l,i}
       + \M g^{\mathsc{nn}}_{l,1}
       + \M H^k_{l,1} - \M H^{k-1}_{l,1}
  \\
    & \M u^{k}_{l,2} 
    && = \M u^{k}_{l,1} 
       + \Delta t \sum_{i=1}^{M_l} w^{\mathsc{nn}}_{l,i2} \M f^{k-1}_{l,i}
       + \M g^{\mathsc{nn}}_{l,2}
       + \M H^k_{l,2} - \M H^{k-1}_{l,2}.
  \end{alignat}
\end{subequations}
Using the relation between $\NM w^{\mathsc{nn}}_l$ and $\NM w^{\mathsc{0n}}_l$ the latter can be rewritten to
\begin{subequations}
 	\begin{alignat}{2}
    & \M u^{k}_{l,1} 
    && = \M u^{k}_{l,0} 
       + \Delta t \sum_{i=1}^{M_l} w^{\mathsc{0n}}_{l,i1} \M f^{k-1}_{l,i}
       + \M g^{\mathsc{nn}}_{l,1}
       + \M H^k_{l,1} - \M H^{k-1}_{l,1}
  \\
    & \M u^{k}_{l,2} 
    &&  = \M u^{k}_{l,0} 
       + \Delta t \sum_{i=1}^{M_l} w^{\mathsc{0n}}_{l,i2} \M f^{k-1}_{l,i}
       + \bigl( \M g^{\mathsc{nn}}_{l,1}
              + \M g^{\mathsc{nn}}_{l,2} \bigr)
       + \bigl( \M H^k_{l,1} + \M H^k_{l,2} \bigr)
       - \bigl( \M H^{k-1}_{l,1} + \M H^{k-1}_{l,2} \bigr).
  \end{alignat}
\end{subequations}
Hence, equivalence requires
\begin{subequations}
 	\begin{alignat}{2}
	  & \M g^{\mathsc{0n}}_{l,1} 
	  && = \M g^{\mathsc{nn}}_{l,1} 
	\\
	  & \M g^{\mathsc{0n}}_{l,2}
	  && = \M g^{\mathsc{nn}}_{l,1} + \M g^{\mathsc{nn}}_{l,2}
  \end{alignat}
\end{subequations}
or, more generally
\begin{equation}
  \NM{\M g}^{\mathsc{0n}}_{l} 
  = \NM{\M C}^{\mathsc{0n}}_l \,\NM{\M g}^{\mathsc{nn}}_{l}.
\end{equation}


\section{Single-level stability and accuracy}
\label{sec:app:single-level-stability}

\begin{figure}[H]
    \centering
    \begin{subfigure}[b]{0.45\textwidth}
        \includegraphics[width=\textwidth]{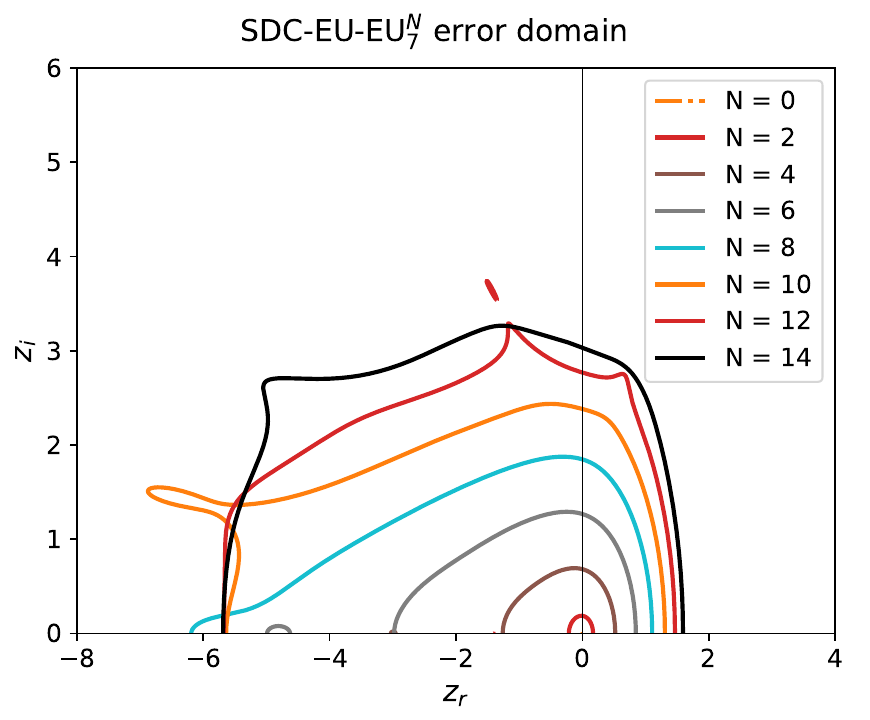}
        \caption{Single-level SDC$_7^N$ with IMEX Euler (corresponds to level 3 MLSDC)}
        \label{acc_xxasdc}
    \end{subfigure}
    \hfill
    \begin{subfigure}[b]{0.45\textwidth}
        \includegraphics[width=\textwidth]{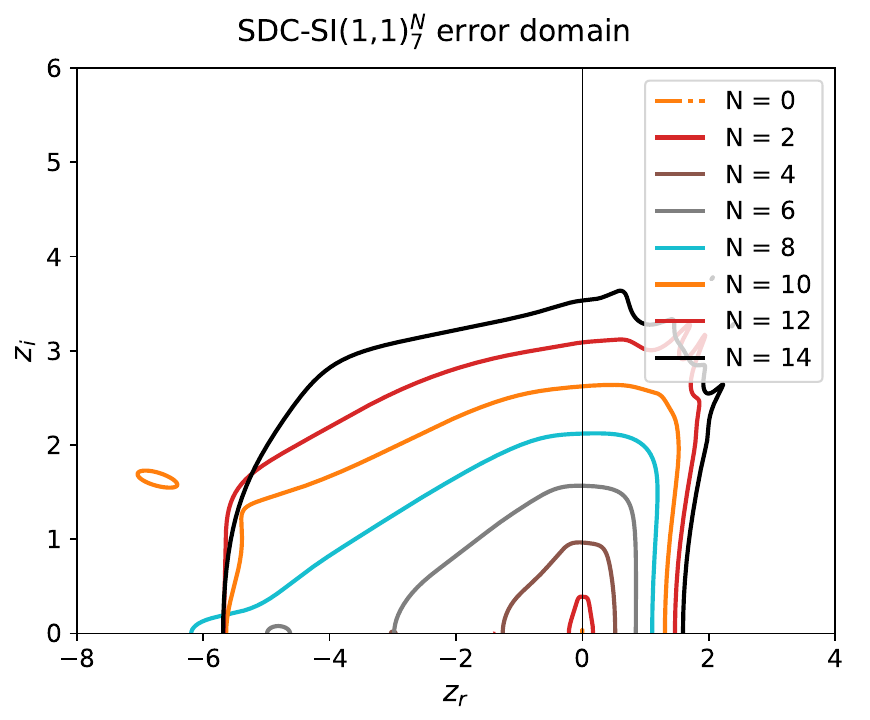}
        \caption{Single-level SDC$_7^N$ with SI(1,1) (corresponds to level 3 MLSDC)}
        \label{acc_xxbsdc}
    \end{subfigure}
    \vspace{1em}
    \begin{subfigure}[b]{0.45\textwidth}
        \includegraphics[width=\textwidth]{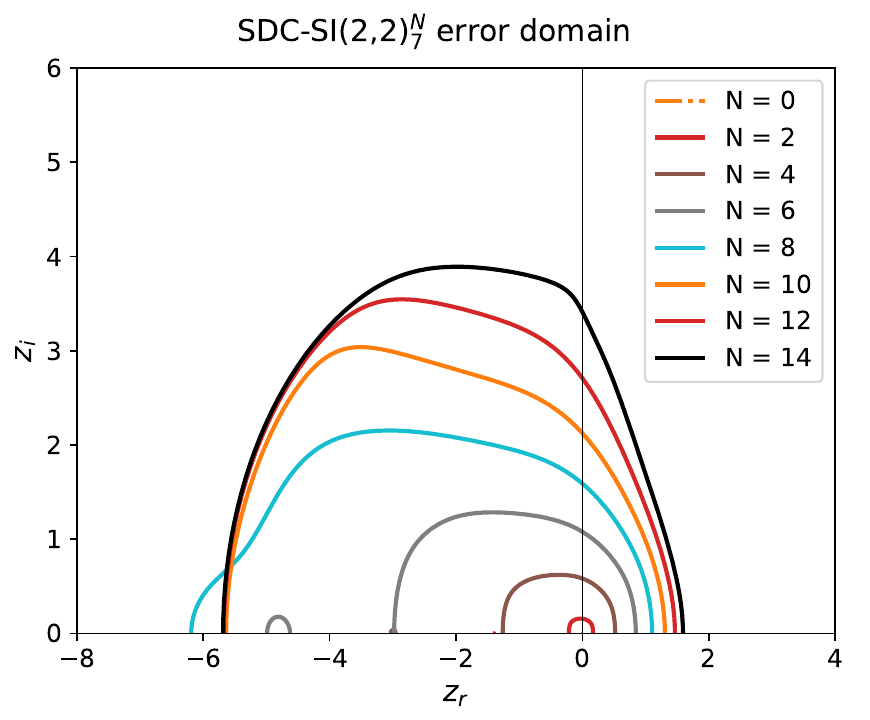}
        \caption{Single-level SDC$_7^N$ with SI(2,2) (corresponds to level 3 MLSDC)}
        \label{acc_xxcsdc}
    \end{subfigure}
    \hfill
    \caption{Accuracy curves ($\varepsilon = 10^{-6}$) for single-level SDC methods with RR points and different time integrators shown after each sweep}
    \label{fig:dahlquist_accuracy_sdc}
\end{figure}

\begin{figure}
    \centering
    \begin{subfigure}[b]{0.45\textwidth}
        \includegraphics[width=\textwidth]{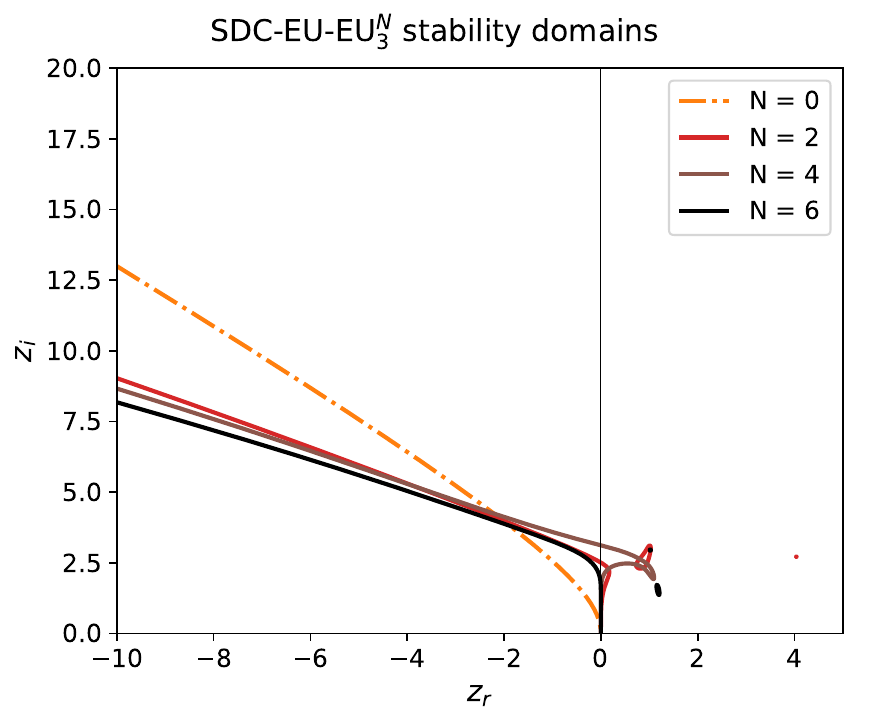}
        \caption{Single-level SDC$_3^N$ with IMEX Euler (corresponds to level 1 MLSDC)}
        \label{xxasdc}
    \end{subfigure}
    \hfill
    \begin{subfigure}[b]{0.45\textwidth}
        \includegraphics[width=\textwidth]{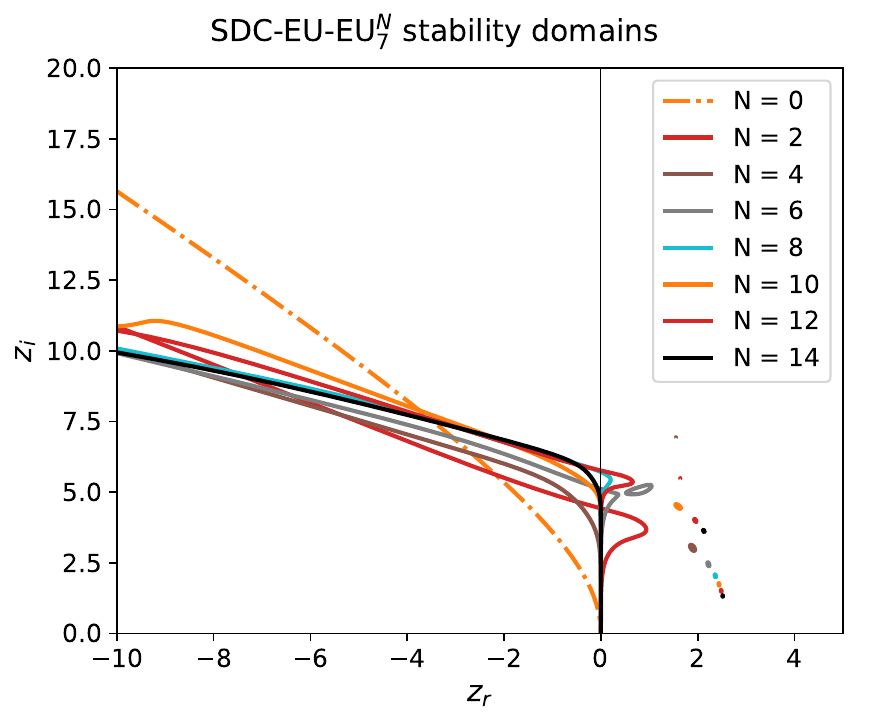}
        \caption{Single-level SDC$_7^N$ with IMEX Euler (corresponds to level 3 MLSDC)}
        \label{xxbsdc}
    \end{subfigure}
    \vspace{1em}
    \begin{subfigure}[b]{0.45\textwidth}
        \includegraphics[width=\textwidth]{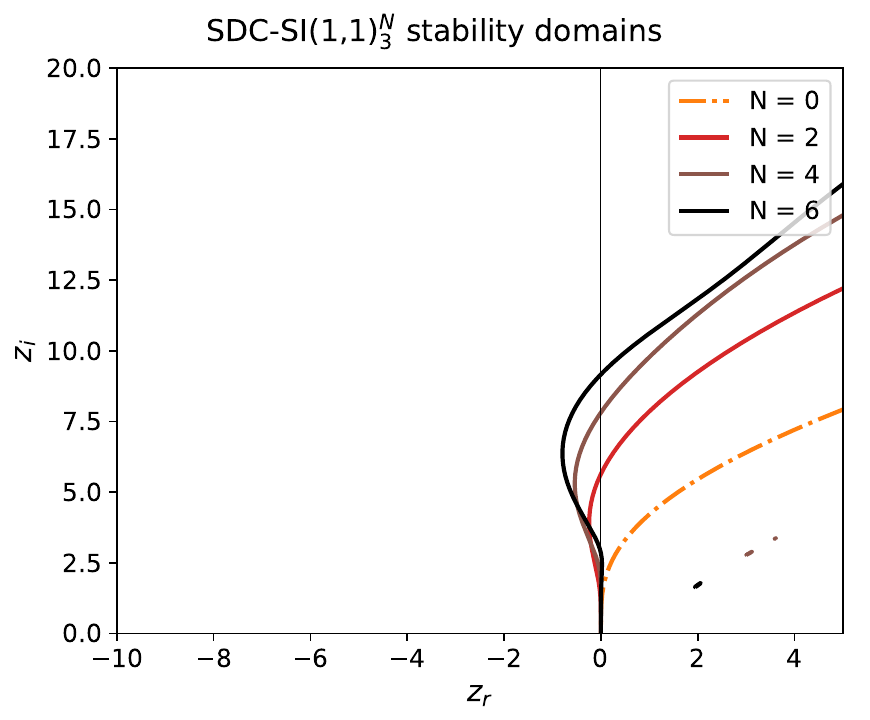}
        \caption{Single-level SDC$_3^N$ with SI(1,1) (corresponds to level 1 MLSDC)}
        \label{xxcsdc}
    \end{subfigure}
    \hfill
    \begin{subfigure}[b]{0.45\textwidth}
        \includegraphics[width=\textwidth]{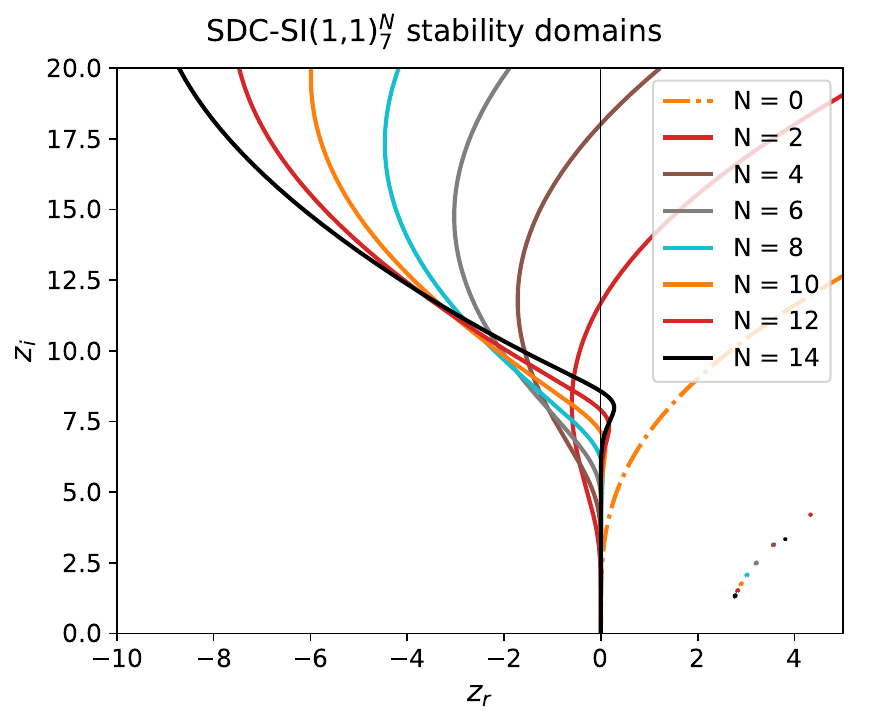}
        \caption{Single-level SDC$_7^N$ with SI(1,1) (corresponds to level 3 MLSDC)}
        \label{xxdsdc}
    \end{subfigure}
    \vspace{1em}
    \begin{subfigure}[b]{0.45\textwidth}
        \includegraphics[width=\textwidth]{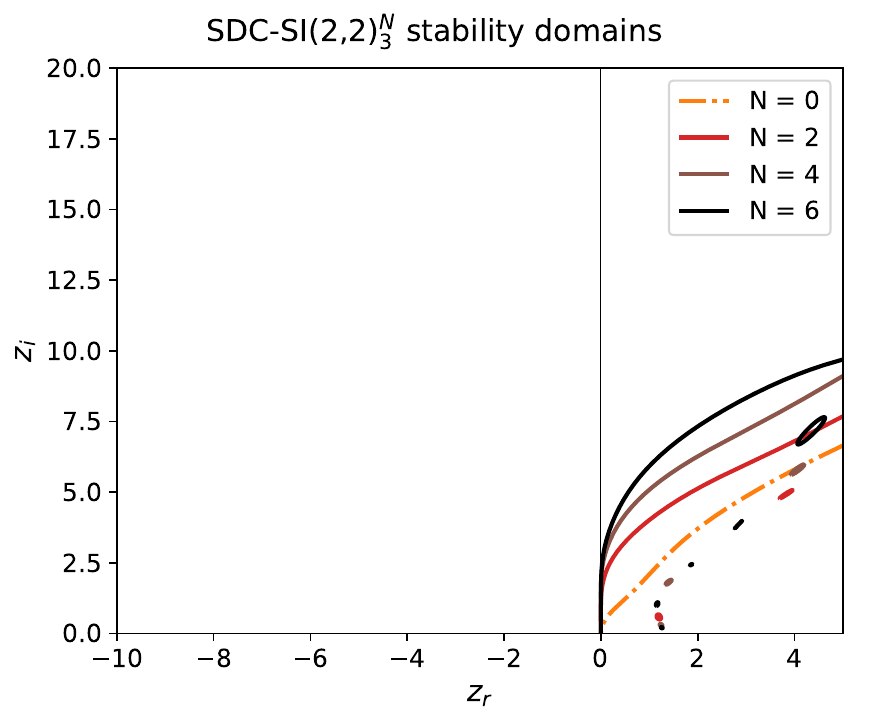}
        \caption{Single-level SDC$_3^N$ with SI(2,2) (corresponds to level 1 MLSDC)}
        \label{xxesdc}
    \end{subfigure}
    \hfill
    \begin{subfigure}[b]{0.45\textwidth}
        \includegraphics[width=\textwidth]{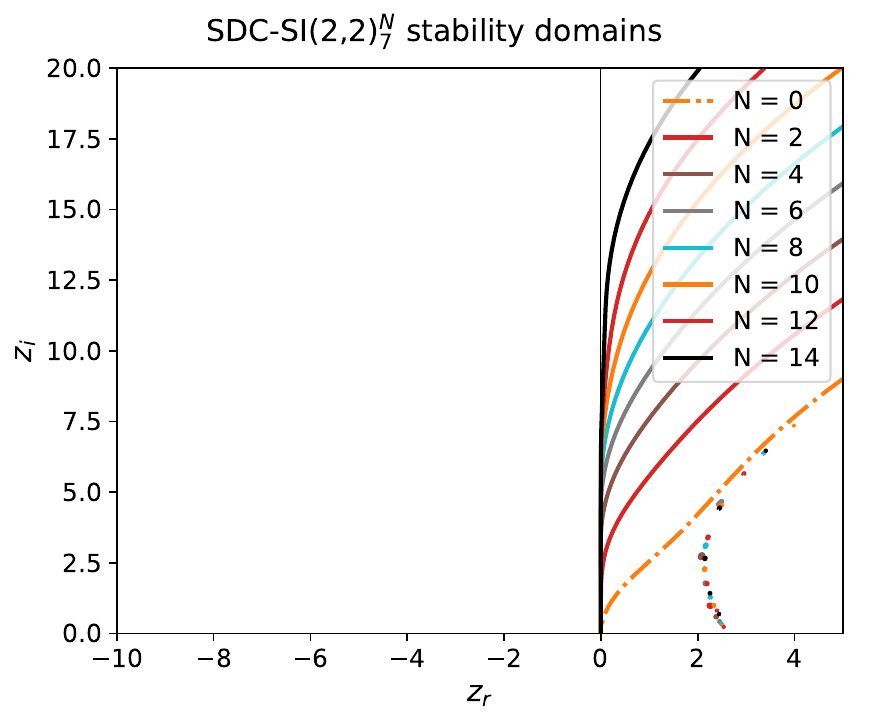}
        \caption{Single-level SDC$_7^N$ with SI(2,2) (corresponds to level 3 MLSDC)}
        \label{xxfsdc}
    \end{subfigure}
    \caption{Neutral stability curves for single-level SDC methods with RR points and different time integrators shown after each second sweep}
    \label{fig:dahlquist_sdc}
\end{figure}